\documentclass[11pt]{article}
\usepackage[margin = 1in]{geometry}

\usepackage{lipsum}
\usepackage{amsfonts}
\usepackage{graphicx}
\graphicspath{ {images/} }
\usepackage{epstopdf}
\usepackage{dsfont}
\usepackage{mathtools}
\usepackage{empheq}
\usepackage{amsfonts}
\usepackage{amsgen,amsthm,amsmath,amstext,amsbsy,amsopn,amssymb,stmaryrd,mathabx,mathrsfs}
\usepackage{comment}
\usepackage[font=footnotesize,labelfont=bf]{caption}
\usepackage{subfigure}

\usepackage{blkarray}

\usepackage{url}
\usepackage{longtable}
\usepackage{mathtools}
\usepackage{multirow}
\usepackage{array}

\usepackage{enumitem}
\usepackage{hyperref}
\usepackage{natbib}
\usepackage{booktabs}
\usepackage{algorithm}
\usepackage{algpseudocode}
\usepackage{xcolor}
\usepackage[utf8]{inputenc} 
\usepackage[T1]{fontenc}

\ifpdf
  \DeclareGraphicsExtensions{.eps,.pdf,.png,.jpg}
\else
  \DeclareGraphicsExtensions{.eps}
\fi

\newcommand{\cB}{\mathcal{B}}

\newcommand{\cP}{\mathcal{P}}
\newcommand{\cX}{\mathcal{X}}

\newcommand{\cD}{\mathcal{D}}

\newcommand{\bbP}{\mathbb{P}}
\newcommand{\bbE}{\mathbb{E}}
\newcommand{\bbR}{\mathbb{R}}
\newcommand{\bbZ}{\mathbb{Z}}

\newcommand{\bbB}{\mathbb{B}}

\renewcommand{\exp}{{\rm{exp}}}
\newcommand{\TV}{{\sf TV}}

\newcommand{\argmax}{\mathop{\rm arg\max}}
\newcommand{\Ham}{{\rm Ham}}

\newcommand{\var}{{\textrm Var}}

\newcommand{\indi}{{\mathds{1}}}

\newcommand{\wh}{\widehat}
\newcommand{\wt}{\widetilde}
\newcommand{\wb}{\widebar}

\newcommand{\DP}{{\textrm{DP}}}

\newtheorem{Definition}{Definition}
\newtheorem{Theorem}{Theorem}

\newtheorem{Lemma}{Lemma}

\newtheorem{Corollary}{Corollary}

\newtheorem{Proposition}{Proposition}

\DeclareMathAlphabet\mathbfcal{OMS}{cmsy}{b}{n}

\hypersetup{
    colorlinks=true,
    linkcolor=blue,
    filecolor=blue,      
    urlcolor=cyan,
    citecolor=blue
}

\makeatletter
\newcommand*{\rom}[1]{\expandafter\@slowromancap\romannumeral #1@}
\makeatother

\allowdisplaybreaks

\begin{document}
	\title{Stability and Accuracy Trade-offs in Statistical Estimation}
			
\author{Abhinav Chakraborty$^{*1}$ ~ Yuetian Luo$^{*2}$~ and ~ Rina Foygel Barber$^3$}
{\let\thefootnote\relax\footnotetext{* The first two authors contributed equally and are listed in alphabetical order.}}

\footnotetext[1]{Department of Statistics, Columbia University}
\footnotetext[2]{Department of Statistics, Rutgers University}
\footnotetext[3]{Department of Statistics, University of Chicago}

	\date{}
	\maketitle

\begin{abstract}
Algorithmic stability is a central concept in statistics and learning theory that measures how sensitive an algorithm's output is to small changes in the training data. Stability plays a crucial role in understanding generalization, robustness, and replicability, and a variety of stability notions have been proposed in different learning settings. However, while stability entails desirable properties, it is typically not sufficient on its own for statistical learning---and indeed, it may be at odds with accuracy, since an algorithm that always outputs a constant function is perfectly stable but statistically meaningless. Thus, it is essential to understand the potential statistical cost of stability. In this work, we address this question by adopting a statistical decision-theoretic perspective, treating stability as a constraint in estimation. Focusing on two representative notions—worst-case stability and average-case stability—we first establish general lower bounds on the achievable estimation accuracy under each type of stability constraint. { We obtain sharp upper bounds for several mean-estimation settings and near-sharp bounds, up to logarithmic factors, for nonparametric regression.} Together, these results characterize the near optimal trade-offs between stability and accuracy across these tasks. Our findings formalize the intuition that average-case stability imposes a qualitatively weaker restriction than worst-case stability, and they further reveal that the gap between these two can vary substantially across different estimation problems.
\end{abstract}


\begin{sloppypar}
\section{Introduction} \label{sec: introduction}
Over the past few decades, stability has emerged as a central principle in modern data science. At a high level, algorithmic stability expresses the idea that the output of an algorithm should not change massively under a slight perturbation to the input data. Intuitively, it is a desirable property for an algorithm on its own---but it also leads to many downstream implications, as well.
For instance, one of the most important and early implications of stability is that stable algorithms generalize \citep{rogers1978finite,devroye1979distribution,kearns1997algorithmic,bousquet2002stability,poggio2004general,shalev2010learnability}. Stability also plays a critical role in a number of other problems, such as model selection \citep{meinshausen2010stability,yu2013stability,shah2013variable,ren2023derandomizing}, differential privacy \citep{dwork2009differential,dwork2014algorithmic}, predictive inference \citep{steinberger2018conditional,barber2021predictive,liang2023algorithmic}, adaptive and selective inference \citep{dwork2015reusable,zrnic2023post} and inference for cross-validation estimate of model and algorithm risks \citep{austern2020asymptotics,bayle2020cross,kissel2022high,luo2024limits}. Stability is also one of three foundational principles of the PCS framework---Predictability, Computability, and Stability---introduced in \cite{yu2020veridical} for veridical data analysis. Finally, stability is closely related to many other notions in learning theory and statistics, such as learnability, replicability and sensitivity \citep{shalev2010learnability,impagliazzo2022reproducibility,kalavasis2023statistical,trillos2025wasserstein}.

Tailored to different downstream tasks, numerous notions of stability have been proposed in the literature. We refer readers to \cite{bousquet2002stability,shalev2010learnability,kim2021black} and the references therein for a few important notions. There are also a few algorithms which are known to satisfy various stability notions. A few canonical examples include $k$-nearest neighbors (kNN) \citep{devroye1979distribution}, ridge regression \citep{bousquet2002stability}, and stochastic gradient descent \citep{hardt2016train,lei2020fine}. Despite its long history in the literature, there is still a lack of a systematic understanding of the relation between different stability notions. This leads to our first motivating question:
\begin{quote}
	 {\it What is the connection between different stability notions, and can we make a quantitative comparison between them?}
\end{quote}

At the same time, when we consider using stable algorithms for our tasks, another important question we need to address is what level of stability we want. Is more stability necessarily better? In fact, while stability has a number of appealing features, it is typically not the whole story in a statistical learning context. For example, in the seminal work of \cite{bousquet2002stability}, stability was leveraged as an important tool to bound the generalization error of a predictor. Typically, as the predictor becomes more stable, it can also have smaller generalization error; however, this does not necessarily imply it has smaller test error or excess risk \citep{chen2018stability}, which is the ultimate goal we care about in the prediction task. A naive but extremely stable predictor is the constant predictor, which ignores the whole dataset and could have a huge test error. Similarly, in the predictive inference context, algorithmic stability enables a jackknife-based predictive interval to have distribution-free coverage guarantee \citep{barber2021predictive}. This is a surprising result as, without stability, jackknife-based predictive intervals can have severe undercoverage. On the other hand, in most practical settings, we also want the predictive interval to be short so that it can be practically useful. Unfortunately, it is not guaranteed that boosting stability would yield a shorter predictive interval. In fact, the reverse implication is more likely to be true, i.e., as the algorithm becomes more stable, the predictor is forced to become more and more like a constant predictor, and as a result, the corresponding statistical efficiency will also degrade. These considerations lead to our second key question:
\begin{quote}
	{\it Is there any trade-off between stability and statistical utility, and if yes, can we characterize it?}
\end{quote}
The main goal of the current paper is to provide a sharp understanding of these trade-offs in the setting of different notions of stability.

\paragraph{Trade-offs and constraints in statistical decision theory.} Statistical decision theory has served as a classical approach to assessing the utility of a statistical estimator in the frequentist sense. Classical statistical decision theory aims to look for an estimator with the smallest risk among all possible estimators. However, this classical approach can suffer from a significant deficiency nowadays, as modern data analysis is complex and often has many constraints. Motivated by that, people have started to incorporate different additional desiderata in statistical estimation, other than a small risk. Encouragingly, over the past decade or so, a flurry of research has been devoted to developing statistical decision theory under various practical constraints, and many of them have even been employed in practice. This includes but is not limited to computational constraints \citep{berthet2013complexity,zhang2014lower,ma2015computational,brennan2020reducibility}, communication constraints \citep{zhang2013information,garg2014communication,braverman2016communication}, and privacy constraints \citep{bassily2014private,dwork2015robust,duchi2018minimax,barber2014privacy,cai2021cost}. In this work, we introduce another desideratum, stability, into this family. Specifically, we take a statistical decision theory point of view, treating stability as a type of constraint in statistical inference tasks. Our focus lies on characterizing the effect of stability in statistical estimation and the fundamental trade-off between stability and accuracy.

\subsection{Our Contributions}\label{sec:contribution}
Building our work on the classical minimax theory, in this paper, we provide a formal framework to compare different notions of stability and characterize the effect of stability as a constraint in statistical estimation. We focus on two general notions of stability considered in the literature: worst-case stability (or uniform stability) and average-case stability. We first provide a couple of general techniques for deriving minimax lower bounds under these types of stability constraints. These bounds are useful in that they not only characterize the ``statistical cost of stability'', but they also allow us to compare different stability notions and stable estimators. 

{ Based on our lower-bound results, we construct stable estimators for a number of canonical estimation problems, including bounded mean estimation, heavy-tailed mean estimation, sparse mean estimation, and nonparametric regression function estimation. For several mean-estimation settings, we provide sharp upper bounds, while for nonparametric regression, we obtain upper bounds that are sharp up to logarithmic factors.}

\begin{figure}
	\centering
	\subfigure[Bounded Distribution Class]{\includegraphics[width=0.46\textwidth]{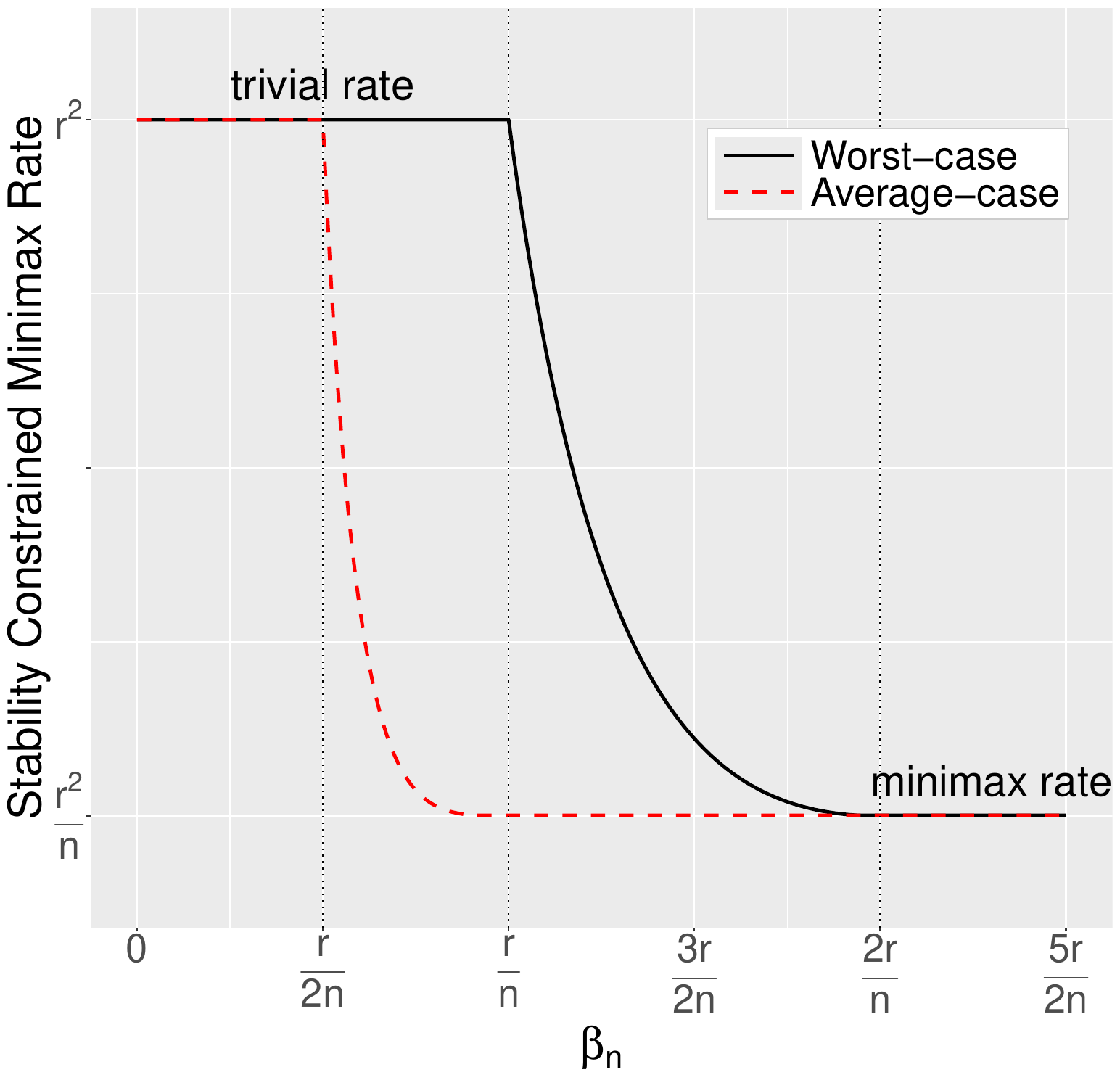}}
	\subfigure[Heavy-tailed Distribution Class]{\includegraphics[width=0.46\textwidth]{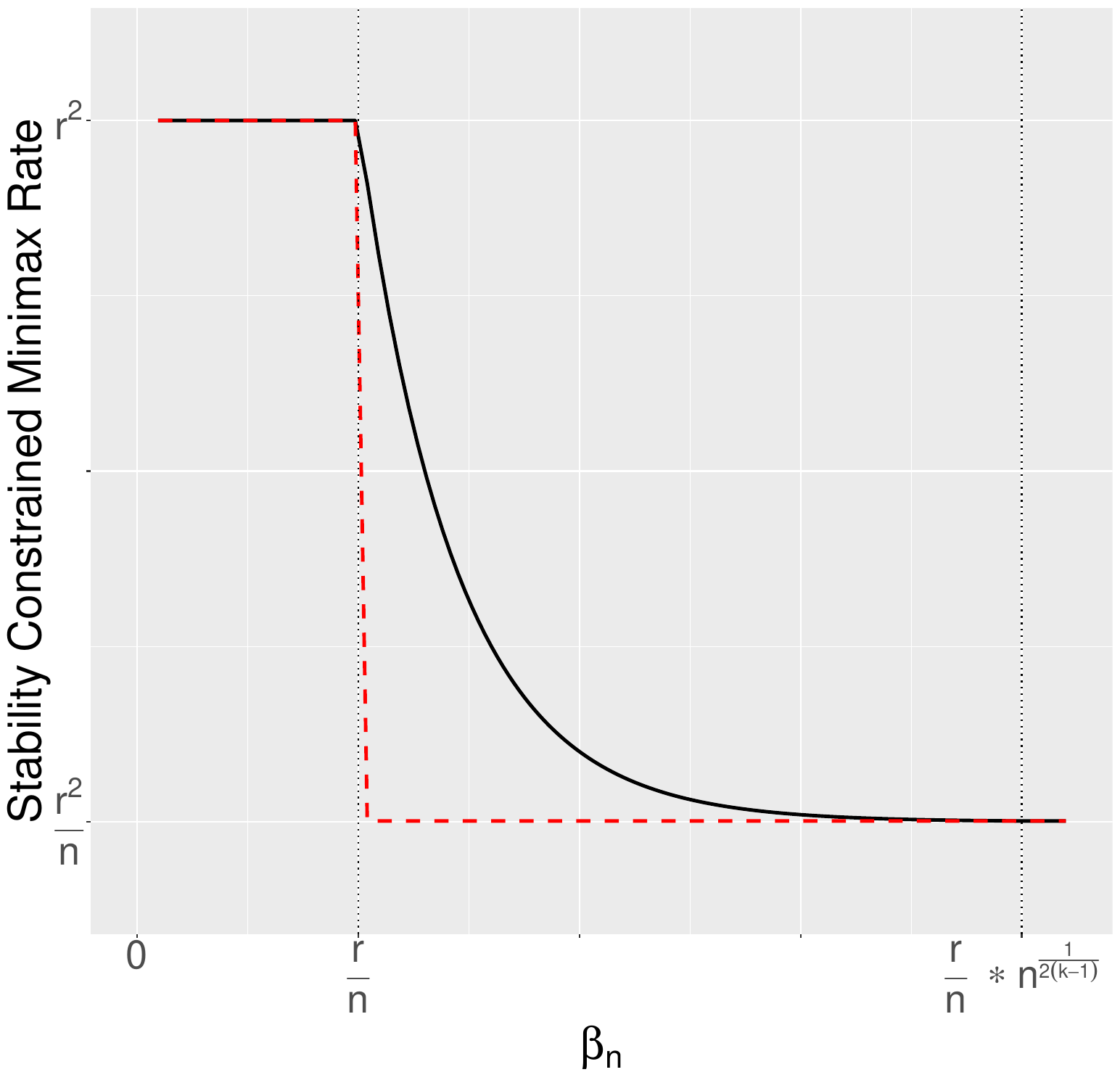}}
	\caption{Stability-vs.-accuracy trade-off curves for mean estimation under worst-case stability and average-case stability. Figure (a) is for the $1$-d mean estimation over a bounded distribution class; figure (b) is for the $1$-d mean estimation over a heavy-tailed distribution class with bounded $k$-th moment $(k \geq 2)$. Here, $r$ denotes the radius of the parameter space, and the vertical dotted lines mark the transition location of the curves.} \label{fig:phase-transition}
\end{figure}

From our results, we confirm the intuition that average-case stability is, in general, a qualitatively weaker constraint than the worst-case stability in statistical estimation. But we also find that the degree of the difference between these two stability notions can vary significantly across different problems. For example, in Figure \ref{fig:phase-transition}, we illustrate the stability and accuracy trade-off curves for mean estimation over the bounded class and heavy-tailed class. We can see that to achieve the same accuracy, the stability requirements under worst-case stability and average-case stability are roughly on the same order in bounded mean estimation, but this is not the case in heavy-tailed mean estimation. More comparison in other examples will be provided in Table \ref{tab:phase-transition} in Section \ref{sec:glimpse-result}.

\subsection{Paper Organization} \label{sec:organization}
We present the background and problem formulation for stability-constrained minimax estimation in Section~\ref{sec:formulation}. In
Section~\ref{sec:lower-bound}, we present our main results on the lower bounds for the stability-constrained minimax risk. { We then apply these lower bounds and construct stable estimators in four canonical problems: bounded mean estimation (Section \ref{sec:bounded-mean}), heavy-tailed mean estimation (Section \ref{sec:heavy-tailed-mean}), sparse mean estimation (Section \ref{sec:sparse-mean}), and nonparametric regression function estimation (Section \ref{sec:nonparametric}).} In Section~\ref{sec:discussion}, we discuss a conversion between stability and differential privacy and its implications, make a comparison to other variants of stability, and mention a few interesting open questions. All proofs are deferred to the Appendix.

\paragraph{Notation}
Given any two numbers $a, b \in \bbR$, let $a \wedge b := \min\{a,b\}$ and $a \vee b := \max\{a, b\}$. For any two nonnegative sequences $\{a_n\}$ and $\{b_n\}$, we write $a_n \asymp b_n$ if there exist constants $c, C>0$ such that $ca_n \leq b_n\leq Ca_n$ for all $n$; $a_n \lesssim b_n$ means that $a_n \leq C b_n$ holds for some constant $C > 0$ independent of $n$. We say $a_n \ll b_n$ if $\lim_{n \to \infty} a_n/b_n \to 0$. Given any vector $v \in \bbR^d$, we use $\|v\|_q$ to denote its $\ell_q$ norm. { For two datasets $\cD_n = \{X_1, \ldots, X_n\}, \cD_n' = \{X_1', \ldots, X_n'\}$ where $\{X_i, X_i'\}_{i=1}^n$ lie in the common space $\cX$, we let $d_{\Ham}(\cD_n, \cD_n') :=\#\{i=1,\cdots, n: X_i \neq X_i'\}$ denote the Hamming distance of two datasets.} The notation $P^{\otimes n}$ means the product distribution of $P$ with $n$ i.i.d. copies. Given any $x \in \bbR$, let $(x)_+ = \max\{x,0\}$. Finally, the notation $[x]_{\tau}$ denotes the value of $x$ clipped at level $\tau$, \[[x]_{\tau} := \left\{ \begin{array}{ll}
	 	\tau & \textnormal{ if } x > \tau,\\
	 	x & \textnormal{ if } -\tau \leq x \leq \tau,\\
	 	-\tau & \textnormal{ if } x < - \tau.
	 \end{array}  \right. \]

\section{Background and Problem Formulation} \label{sec:formulation}
In this section, we begin by reviewing the classical minimax framework, and then introduce the notions of stability studied in this paper and the minimax framework under stability constraints.

\subsection{Classical Minimax Framework}\label{sec:class-minimax}
Suppose we have a class of distributions $\cP$ with support $\cX$. Let $\theta(P) \in \Theta$ denote a functional defined on $\cP$. For simplicity, in this paper, we focus on $\Theta = \bbR^d$ for some $d \geq 1$. Given a fixed but unknown distribution $P \in \cP$ and a dataset $\cD_n:=(X_1, \ldots, X_n)$ consisting of $n$ i.i.d. data drawn from $P$, our goal is to estimate the target parameter $\theta(P)$. Formally, an estimator of $\theta(P)$ is a measurable function $\wh{\theta}: \cX^n\times[0,1]\to \bbR^d$, returning an estimate $\wh{\theta}(\cD_n,\xi)$, where $\xi\in[0,1]$ could be viewed as an (optional) random seed to allow randomness in the estimator. In order to assess the quality of any estimator, we let $\rho: \bbR^d \times \bbR^d \to [0, \infty)$ be some non-negative measure of the error, and consider the associated risk function:
\begin{equation*}
	R_n(\theta(P), \wh{\theta}) = \bbE[\rho(\wh{\theta}(\cD_n; \xi), \theta(P))],
\end{equation*}
where implicitly, the expected value is taken with respect to the distribution of $\cD_n\sim P^{\otimes n}$ and $\xi \sim \text{Unif}[0,1]$.
Throughout the paper, we will mainly focus on the mean-square risk, i.e., $\rho(\theta_1,\theta_2) = \|\theta_1 - \theta_2\|_2^2$. In the classical minimax setting, for each estimator $\wh{\theta}$, we can compute its worst-case risk $\sup_{P \in \cP} R_n(\theta(P), \wh{\theta})$, and rank estimators based on this quantity. An estimator $\wh{\theta}$ is called rate optimal (or simply optimal) if it can achieve the minimax rate up to a constant factor, $R_n(\theta(\cP),\wh{\theta}) \lesssim R_n(\theta(\cP))$, where the minimax rate is defined as
\begin{equation} \label{eq:class-minimax}
	R_n(\theta(\cP)) := \inf_{\wh{\theta}} \sup_{P \in \cP} R_n(\theta(P), \wh{\theta}) =\inf_{\wh{\theta}} \sup_{P \in \cP} \bbE \|\wh{\theta}(\cD_n; \xi)- \theta(P)\|_2^2.
\end{equation}

\subsection{Worst-case and Average-case Stability} \label{sec:stability-notions}
Generically speaking, the various notions of stability in the literature can be classified into two big categories: worst-case stability (or uniform stability) and average-case stability. Worst-case stability describes that the output of the procedure is insensitive to the change of any given data point in the dataset. So the influence of every data point should be low to satisfy worst-case stability. A few common worst-case-stable estimators have been identified in the literature \citep{bousquet2002stability,hardt2016train}, but it can be hard to achieve in many other scenarios. For example, \cite{xu2011sparse} showed that any algorithm that is constrained to obey sparsity cannot be worst-case-stable. This motivates researchers looking for relaxed stability notions that are still useful in the downstream statistical inference. Average-case stability is one such notion, which is based on the average perturbation on the output of the procedure if we replace one data point at random: this requires that, on average, the influence of each data point is low. Bagging \citep{breiman1996bagging} is a method that can achieve average-case stability in fairly general settings \citep{elisseeff2005stability,soloff2023bagging}.

Now, let us define these notions of stability precisely. First, we consider worst-case stability.
\begin{Definition} \label{def:worst-cases-stability}
	Given any estimator $\wh{\theta}: \cX^n\times[0,1] \to \bbR^d$, we say $\wh{\theta}$ satisfies {\it worst-case stability} with parameter $\beta_n$ (or, equivalently, we say $\wh{\theta}$ is $\beta_n$-worst-case-stable) if
\begin{equation} \label{def:worst-stability}
	\begin{split}
		\sup_{\cD_n\sim\cD_n'} \bbE_{\xi} \| \wh{\theta}(\cD_n; \xi) - \wh{\theta}(\cD_n'; \xi) \|_2 \leq \beta_n,
	\end{split}
\end{equation} where $\cD_n, \cD_n' \in \cX^n$, and the notation $\cD_n \sim \cD_n'$ means $\cD_n$ and $\cD'_n$ differ by only one data point (i.e., $d_{\Ham}(\cD_n, \cD_n')\leq1$).
\end{Definition}
Note that the expectation in Definition \ref{def:worst-cases-stability} is only taken with respect to the random seed $\xi$, while the data could be arbitrary. Worst-case stability is also known as `global sensitivity' in the differential privacy (DP) literature. In particular, the celebrated Laplace mechanism adds Laplace noise proportional to the global sensitivity to construct differentially private estimators \citep{dwork2008differential}. In Section \ref{sec:stability-dp}, we will provide a detailed exploration of the connection and difference of DP as a constraint compared to stability. 

Next, we turn to average-case stability.
\begin{Definition} \label{def:average-case-stability}
	Given any estimator $\wh{\theta}: \cX^n\times[0,1] \to \bbR^d$, we say $\wh{\theta}$ satisfies {\it average-case stability} with parameter $\beta_n$ (or, equivalently, we say $\wh{\theta}$ is $\beta_n$-average-case-stable) if
\begin{equation} \label{eq:average-case-stability}
	\frac{1}{(n+1)^2} \sup_{\cD_{n+1} \in \cX^{n+1}} \sum_{1 \leq i, j \leq n+1} \bbE_{\xi} \|\wh{\theta}(\cD_{n+1}^{\setminus i}; \xi )  - \wh{\theta}(\cD_{n+1}^{\setminus j}; \xi )\|_2 \leq \beta_n,
\end{equation} where $\cD_{n+1}$ denotes a dataset of size $n+1$, and the notation $\cD_{n+1}^{\setminus i}$ denotes the dataset obtained by deleting the $i$th data point from $\cD_{n+1}$.
\end{Definition}
 Note again that in \eqref{eq:average-case-stability}, the dataset could be arbitrary, but the influence of each data point on the estimator is averaged out. The existing notions of average-case stability are typically defined based on dropping one data point at random \citep{soloff2023bagging}. In practice, we expect the two definitions to be essentially equivalent. We refer readers to Appendix \ref{sec:two-average-case-stability} for more details. 
 
 In general, the worst-case and average-case stability can be unified under the following notion of $\ell_p$-stability.
 \begin{Definition} \label{def:l-gamma-stability}
	Given any $p \in [1,\infty]$ and an estimator $\wh{\theta}: \cX^n\times[0,1] \to \bbR^d$, we say $\wh{\theta}$ satisfies {\it $\ell_p$-stability} with parameter $\beta_n$ (or, equivalently, we say $\wh{\theta}$ is $\beta_n$-$\ell_p$-stable) if
\begin{equation} \label{eq:lp-stability}
	\left\{\frac{1}{(n+1)^2} \sup_{\cD_{n+1} \in \cX^{n+1}} \sum_{1 \leq i, j \leq n+1} \bbE_{\xi} \|\wh{\theta}(\cD_{n+1}^{\setminus i}; \xi )  - \wh{\theta}(\cD_{n+1}^{\setminus j}; \xi )\|^{p}_2 \right\}^{1/p}\leq \beta_n.
\end{equation}
(For the case $p=\infty$, we should interpret this definition as taking a limit as $p\to\infty$, i.e., we require
$\sup_{\cD_{n+1} \in \cX^{n+1}} \max_{1 \leq i, j \leq n+1} \bbE_{\xi} \|\wh{\theta}(\cD_{n+1}^{\setminus i}; \xi )  - \wh{\theta}(\cD_{n+1}^{\setminus j}; \xi )\|_2 \leq \beta_n$.)
\end{Definition}
Note that when $p=\infty$, $\ell_p$-stability reduces to worst-case stability (we can view $\cD_n$ and $\cD_n'$ in \eqref{def:worst-stability} as $\cD_{n+1}^{\setminus i}$ and $\cD_{n+1}^{\setminus j}$ in \eqref{eq:lp-stability}), while when $p=1$, this definition reduces to average-case stability. By definition, $\ell_{p}$-stability implies $\ell_q$-stability whenever $p \geq q$. In this paper, we will mainly focus on $p=\infty$ and $p=1$, i.e., worst-case stability and average-case stability, but in some examples we extend our results to general $p$.

{  Finally, we note that in the literature, there is another commonly used average-case stability notion where the average is over the randomness of the data. It has been used in showing the generalization of an algorithm and consistency of cross-validation risk. We will briefly discuss that notion and make a comparison in Section \ref{sec:compare-diff-stability}.}

\subsection{Stability-Constrained Minimax Risks}
Given the notions of stability of interest, we are now ready to define stability-constrained minimax risk:
\begin{equation*}
	\begin{split}
	R_{n,p}(\theta(\cP), \beta_n) := \inf_{\wh{\theta} \textnormal{ is $\beta_n$-$\ell_p$-stable} } \sup_{P \in \cP} \bbE_{\cD_n, \xi} \|  \wh{\theta}(\cD_n; \xi) - \theta(P) \|_2^2.
	\end{split}
\end{equation*} 
Notice that given any $p \in [1,\infty]$, $R_{n,p}(\theta(\cP), \beta_n)$ is a nonincreasing function with respect to $\beta_n$, and given any $\beta_n$, $R_{n,p}(\theta(\cP), \beta_n)$ is a nondecreasing function of $p$.
When $p=\infty$ or $1$, $R_{n,\infty}(\theta(\cP), \beta_n)$ and $R_{n,1}(\theta(\cP), \beta_n)$, respectively, denote the minimax risk under the worst-case and average-case stability constraint. In addition, when $\beta_n \to \infty$, the constraint on stability becomes vacuous, and $R_{n,p}(\theta(\cP), \beta_n)$ reduces to the classical minimax risk \eqref{eq:class-minimax}. { The main goal of this paper is to study $R_{n,p}(\theta(\cP), \beta_n)$ as a function of $\beta_n$ across different $p$'s and different problems, derive sharp or near-sharp bounds for these risks, construct stable estimators attaining the corresponding upper bounds, and compare how worst-case stability and average-case stability differ as constraints in statistical estimation tasks.}

We conclude this section by introducing some terminology to describe the phase transition patterns for the stability-vs.-accuracy trade-off curves, i.e., $R_{n,p}(\theta(\cP), \beta_n)$. We will use this terminology to discuss and compare our findings throughout the paper. 

\begin{Definition}[Sharp vs.\ gradual phase transition] \label{def:phase-transition-pattern}
	Given a stability-vs.-accuracy trade-off curve $R_n(\theta(\cP), \beta_n)$, we say it has a {\it sharp} phase transition at $\beta_n^*$ if $R_n(\theta(\cP), \beta_n) \asymp R_n(\theta(\cP), 0)$ when $\beta_n \ll \frac{\beta_n^*}{\textnormal{poly}\log n}$ and $R_n(\theta(\cP), \beta_n) \asymp R_n(\theta(\cP), \infty)$ when $\beta_n \gg \beta_n^* \cdot \textnormal{poly}\log n$. Otherwise, we say it has a {\it gradual} phase transition.
\end{Definition}
A sharp phase transition means that a certain stability budget is critical for this estimation task: if we require stronger stability below the critical threshold, then the best estimator we can hope for is a trivial constant estimator; while, as long as we only need stability above the critical threshold, then we can find an estimator that achieves the unconstrained minimax rate. In other words, we can get $\Omega(\beta_n^*)$ level stability for free while still maintaining the optimal performance, but when $\beta_n = o(\beta_n^*)$, nontrivial estimation is impossible, and there is no trade-off between stability and accuracy. When the stability-vs.-accuracy trade-off curve has a gradual phase transition, then we instead see a smoother trade-off between these two quantities.

\subsection{A Glimpse of the Result} \label{sec:glimpse-result}
With the notion of sharp and gradual phase transitions defined in Definition \ref{def:phase-transition-pattern}, the following table provides a brief summary of the characterization of $R_{n,p}(\theta(\cP), \beta_n)$ across the four canonical problems studied in this paper.

\begin{table}[h]
	\centering
	\begin{tabular}{c|c|c|c}
		\hline
		\multirow{3}{9em}{Estimation tasks} & \multicolumn{2}{c|}{Phase transition behavior} & Compare transition \\
		\cline{2-3}
		 &worst-case stability & average-case stability  &  thresholds for \\
         & ($\ell_p$ for $p=\infty$)& ($\ell_p$ for $p=1$) & $p=1$ vs $p=\infty$\\
		\hline
		Bounded mean (Sec. \ref{sec:bounded-mean}) & sharp & sharp & { same order} \\
		\hline
		Heavy-tailed mean (Sec. \ref{sec:heavy-tailed-mean}) & gradual & sharp & different \\
		\hline 
		Sparse mean (Sec. \ref{sec:sparse-mean})& sharp & sharp & { same order} \\
		\hline 
		Nonparametric reg. (Sec. \ref{sec:nonparametric}) & {gradual} & {sharp} & different\\
		\hline
	\end{tabular}
	\caption{Summary of stability-vs.-accuracy trade-off patterns across different problems under both worst-case stability and average-case stability (i.e., $\ell_\infty$ and $\ell_1$).} \label{tab:phase-transition}
\end{table} 
Perhaps surprisingly, in many examples, we find the stability-vs.-accuracy trade-off curve exhibits a sharp phase transition, and moreover, may sometimes exhibit the same phase transition at $p=1$ as at $p=\infty$, despite the strong differences between these two notions of stability (average-case or worst-case).

\section{General Tools for Lower Bounding $R_{n,p}(\theta(\cP), \beta_n)$} \label{sec:lower-bound}

In this section, we present two general results on lower-bounding the stability-constrained risk. Let us begin with a lower bound under worst-case stability.
\begin{Theorem}\label{th:worst-lower-bound-multi}
For the statistical estimation setting considered in Section \ref{sec:class-minimax} under the worst-case stability constraint, we have
	\begin{equation} \label{ineq:risk-lower-bound}
		R_{n,\infty}(\theta(\cP), \beta_n) \geq \sup_{P_0, P_1 \in \cP} \left[ \frac{ (\|\theta(P_0) - \theta(P_1)\|_2 - \bbE[ d_{\Ham}(\cD_n, \cD_n') ] \beta_n)_{+} }{2} \right]_+^2 \vee R_n(\theta(\cP)),
	\end{equation} where in the expected value, the joint distribution of $(\cD_n, \cD_n')$ could be any coupling between $P_0^{\otimes n}$ and $P_1^{\otimes n}$.
\end{Theorem} 
 Notice that the first term on the right-hand side of \eqref{ineq:risk-lower-bound} appears due to the stability constraint. In a typical setting, the two terms $\|\theta(P_0) - \theta(P_1)\|_2$ and $\bbE[ d_{\Ham}(\cD_n, \cD_n') ] \beta_n$ have a trade-off, as when $\|\theta(P_0) - \theta(P_1)\|_2$ increases, more likely, the Hamming distance between $\cD_n$ and $\cD_n'$ will also increase. The proof idea of Theorem \ref{th:worst-lower-bound-multi} is inspired by the classical Le Cam's two-point method for lower bounding the minimax risk \eqref{eq:class-minimax}, but it also differs in an important way. The classical Le Cam's method first reduces the unconstrained minimax risk to a two-point hypothesis testing problem and then lower bounds the testing risk. But in our context, it is not clear how to incorporate the stability constraint into the testing problem. To overcome that obstacle, here, instead of reducing to testing problems, we directly look at two datasets $\cD_n$ and $\cD_n'$ with i.i.d. data drawn from two distributions $P_0$ and $P_1$, respectively. { On the one hand, since $\wh{\theta}$ is $\beta_n$-worst-case-stable, we would expect $\wh{\theta}(\cD_n)$ and $ \wh{\theta}(\cD_n')$ to be close to each other. Concretely}
 \begin{equation} \label{ineq:stability-induced-upper-bound}
 	\begin{split}
 		\bbE\|\wh{\theta}(\cD_n) - \wh{\theta}(\cD_n')\|_2 \leq \beta_n \bbE  [ d_{\Ham}(\cD_n, \cD_n')],
 	\end{split}
 \end{equation} where the expectation is taken with respect to any given coupling between $\cD_n$ and $\cD_n'$. { On the other hand, if $\wh{\theta}(\cD_n)$ and $\wh{\theta}(\cD_n')$ are accurate estimators of $\theta(P_0)$ and $ \theta(P_1)$, respectively, their difference can also be lower bounded as follows}, 
 \begin{equation} \label{ineq:upperbound-risk}
 \begin{split}
 	\bbE\|\wh{\theta}(\cD_n) - \wh{\theta}(\cD_n')\|_2 &= \bbE\|\wh{\theta}(\cD_n) - \theta(P_0) + \theta(P_0) - \theta(P_1) + \theta(P_1) - \wh{\theta}(\cD_n')\|_2 \\ &\geq \|\theta(P_0) - \theta(P_1)\|_2 - 2 \sup_{P \in \cP} \bbE_P\|\wh{\theta} - \theta(P)\|_2.
 \end{split}
 \end{equation} Since \eqref{ineq:stability-induced-upper-bound} and \eqref{ineq:upperbound-risk} hold for any $\beta_n$-worst-case-stable estimator, combining these bounds yields a lower bound for $R_{n,\infty}(\theta(\cP), \beta_n)$. The detailed proof of Theorem \ref{th:worst-lower-bound-multi} will be provided in Appendix \ref{app:proof-lower-bound}.

Now, we move on to the statistical lower bound under the general $\ell_p$-stability constraint. When $p<\infty$, the bound in \eqref{ineq:stability-induced-upper-bound} no longer holds, but { the intuition is still correct: if the estimator is stable, then its evaluation on two different datasets can not be too far away. When $p<\infty$, we first find it is without loss of generality to assume $\wh{\theta}(\cD_n)$ is symmetric to its input and furthermore, we can modify the analysis in \eqref{ineq:stability-induced-upper-bound} to provide a different upper bound on $\|\wh{\theta}(\cD_n) - \wh{\theta}(\cD_n')\|_2$ using the $\ell_p$ metric by leveraging the symmetry property}; we defer details to the Appendix \ref{app:proof-lower-bound}. This next result provides a generic lower bound for $R_{n,p}(\theta(\cP), \beta_n)$. 

\begin{Theorem}\label{th:average-lower-bound-multi}
Given any $p \in [1, \infty]$,
	\begin{equation} \label{ineq:ave-risk-lower-bound}
		R_{n,p}(\theta(\cP), \beta_n) \geq \sup_{P_0, P_1 \in \cP} \left[ \frac{ (\|\theta(P_0) - \theta(P_1)\|_2 - (n+1) (\log n + 1)^{1/p} \beta_n)_{+} }{2} \right]^2 \vee R_n(\theta(\cP)).
	\end{equation} If, in addition, the functional of interest $\theta(P)$ is linear in $P\in\cP$, then for any $P_0, P_1 \in \cP$ such that $(1-t) P_0 + tP_1 \in \cP$ for all $t \in [0,1]$, 
\begin{equation}\label{ineq:average-lower-linear-fun}
		\begin{split}
			&R_{n,p}(\theta(\cP), \beta_n) \geq R_n(\theta(\cP))\vee \\
			&\sup_{\eta \in [0,1/2]}\left( \frac{ (1- 2 \eta) \|\theta(P_0)  - \theta(P_1) \|_2  - 2^{1/p} (n+1) \beta_n \left( \frac{1-2\eta + \frac{p}{2(n+1)} }{(2\eta)^{1/p}} \wedge (n+1)^{1/p} \right) }{2} \right)^2_{+}.
		\end{split}
	\end{equation}
\end{Theorem}
A corollary of Theorem \ref{th:average-lower-bound-multi} for linear functionals under the average-case stability, i.e., $p=1$, is given as follows.
\begin{Corollary} \label{coro:average-case-stability}
	Suppose $\theta(P)$ is linear in $P$ and $p=1$. Then for any $P_0, P_1 \in \cP$ such that $(1-t) P_0 + tP_1 \in \cP$ for all $t \in [0,1]$, $$R_{n,1}(\theta(\cP), \beta_n) \geq \left( \frac{1}{4} \|\theta(P_0)  - \theta(P_1) \|_2 - 2(n+1) \beta_n \right)_+^2 \vee R_n(\theta(\cP)).$$
\end{Corollary}
Since $R_{n,p}(\theta(\cP), \beta_n)$ is a nondecreasing function of $p$, Corollary \ref{coro:average-case-stability} immediately suggests that when the parameter space is unbounded and connected, then $R_{n,p}(\theta(\cP), \beta_n)$ is infinite for any $p \geq 1$ unless $\beta_n = \infty$. A similar phenomenon has also appeared in statistical estimation under the pure differential privacy constraint, where the task is only possible when the parameter space is bounded \citep{bun2021private}. Thus, for the rest of the article, we will be mainly interested in the setting where the parameter space is bounded. 

{ Next, we will show how Theorem \ref{th:worst-lower-bound-multi} and \ref{th:average-lower-bound-multi} can be used to derive lower bounds in a number of common statistical estimation tasks, including bounded mean estimation, heavy-tailed mean estimation, sparse mean estimation, and nonparametric regression. For bounded mean estimation and sparse mean estimation, we provide matching upper and lower bounds for $R_{n,p}(\theta(\cP), \beta_n)$ under general $\ell_p$-stability for all $p$. For heavy-tailed mean estimation and nonparametric regression, we focus on $p=\infty$ and $p=1$, i.e., worst-case stability and average-case stability; in the nonparametric setting, our upper and lower bounds are near-sharp up to logarithmic factors.} 

\section{Mean Estimation with Bounded Support} \label{sec:bounded-mean}
In this section, we study the first of our four examples: the problem of estimating the mean of a distribution, in a bounded setting.
Let $$\cP^d_{\textnormal{bound}}(r) = \{  \textnormal{distributions } P \textnormal{ supported on } \bbB_2(r) \subseteq \bbR^d \}$$ be the class of distributions on $\bbR^d$ with support on $\bbB_2(r)=\{x \in \bbR^d:\|x\|_2 \leq r\}$, a $\ell_2$-ball with radius $r$. Without any stability constraint, it is easy to check that the minimax risk for mean estimation in this family is of rate \[R_n(\theta(\cP)) \asymp \frac{r^2}{n},\] and it can be achieved by the sample mean. Under the stability constraint, we have the following result which establishes a sharp phase transition of $R_{n,p}\left(\theta(\cP^d_{\textnormal{bound}}(r)), \beta_n\right)$ at $\beta_n^* = \frac{r}{n}$ for any $p$. 
\begin{Theorem} \label{th:worst-bounded} Given any $p \in[1,\infty]$, \begin{equation*}
	\begin{split}
		R_{n,p}\left(\theta(\cP^d_{\textnormal{bound}}(r)), \beta_n\right) \asymp \left\{ \begin{array}{cc}
			\frac{r^2}{n} & \text{ if } \beta_n \geq \frac{2r}{n}\\
			r^2 & \text{ if } \beta_n \leq \frac{r}{10n}.
		\end{array} \right.
	\end{split}
\end{equation*} One estimator that achieves this rate is the shrinkage estimator $ \left( \frac{n \beta_n}{2r} \wedge 1 \right) \wb{X}$, where $\wb{X}$ is the sample mean.
\end{Theorem} In Theorem \ref{th:worst-bounded}, we can see that the phase transition of $R_{n,p}\left(\theta(\cP^d_{\textnormal{bound}}(r)), \beta_n\right)$ happens around $\beta_n \asymp \frac{r}{n}$ for all $p$. This suggests that stability can be obtained `for free' for mean estimation in this class at level $\beta_n\gtrsim r/n$ (in fact, in this case, we may simply use the sample mean as the estimator), while for smaller $\beta_n$, nontrivial estimation is impossible.

Next, we show that for this simple class, we can derive a more refined characterization of how the stability-vs.-accuracy trade-off curve behaves when $\beta_n \in [ \frac{r}{10n}, \frac{2r}{n}]$ in the special setting when $d = 1$.  

\subsection{Refined Characterization When $d = 1$}
Our first refined result is on characterizing the exact location where the phase transition occurs.
\begin{Theorem} \label{th:bounded-worst-refined}
	Given any $p \geq 1$ and any fixed constant $\delta > 0$, there exists a sufficiently large $C >0$ depending on $\delta$ only such that as long as $n \geq C$,
\begin{equation*}
	\begin{split}
		R_{n,p}\left(\theta(\cP^1_{\textnormal{bound}}(r)), \beta_n\right) \asymp \left\{ \begin{array}{cc}
			\frac{r^2}{n} & \text{ if } \beta_n \geq \frac{2^{1-1/p}\cdot r}{n}\\
			r^2 & \text{ if } \beta_n \leq \frac{2^{1-1/p}\cdot r}{n(1+ \delta)}.
		\end{array} \right.
	\end{split}
\end{equation*} One estimator that achieves this rate is $\left( \frac{n \beta_n}{2^{1-1/p}\cdot r} \wedge 1 \right) \wb{X}$.
\end{Theorem} Theorem \ref{th:bounded-worst-refined} shows that if we zoom in, we can see that different $p$'s in $\ell_p$-stability do have a different effect on the exact location where the phase transition happens. For example, the exact location where the phase transition happens shifts from $\beta_n^* = \frac{r}{n}$ in $\ell_{1}$-stability to $\beta_n^* = \frac{2r}{n}$ in $\ell_{\infty}$-stability.

Our second refined result is to further characterize how does $R_{n,p}\left(\theta(\cP^1_{\textnormal{bound}}(r)), \beta_n\right)$ transit from $r^2$ to $r^2/n$ in the narrow phase transition regime. We will focus on the problem when $p=\infty$ and $p=1$. In particular, when $p=\infty$, we will provide an exact formula for $R_{n,\infty}\left(\theta(\cP^1_{\textnormal{bound}}(r)), \beta_n\right)$ over the entire range of $\beta_n$; and when $p=1$, we are also able to capture the optimal rate of $R_{n,1}\left(\theta(\cP^1_{\textnormal{bound}}(r)), \beta_n\right)$ over the entire range of $\beta_n$. Perhaps surprisingly, we will see that if we look for optimal estimators in this narrow regime, then the naive shrinkage estimator is actually not enough.  
\begin{Theorem} \label{th:worst-bounded-refined} Suppose $d = 1$, then we have 
\begin{equation*}
	\begin{split}
		R_{n,\infty}\left(\theta(\cP^1_{\textnormal{bound}}(r)), \beta_n\right) = \left(\frac{\left( \frac{2r}{n \beta_n} - 1 \right)_{+}}{1 + \left( \frac{2r}{n \beta_n} - 1 \right)_{+} } \right)^2 r^2 \vee \frac{r^2}{\left(\sqrt{n} + 1 \right)^2},
	\end{split}
\end{equation*}and one optimal estimator, which achieves the above equality, is given by $\wh{\theta} = \left( \frac{n \beta_n}{2r} \wedge \frac{1}{1 + \frac{1}{\sqrt{n}} } \right)\bar{X}$.
\end{Theorem} { The exact risk computation of the estimator in Theorem \ref{th:worst-bounded-refined} involves a careful analysis on the trade-off of its bias and variance.}
When $\beta_n \leq \frac{r}{10n}$ and $\beta_n \geq \frac{2r}{n}$, we can see the results in Theorem \ref{th:worst-bounded} can be reduced from Theorem \ref{th:worst-bounded-refined}. But Theorem \ref{th:worst-bounded-refined} also tells us what happens when $\beta_n \in [\frac{r}{10n}, \frac{2r}{n}]$. In particular, for any $\delta_n \geq 0$, Theorem \ref{th:worst-bounded-refined} implies that
\begin{equation} \label{eq:worst-bounded-refined}
	R_{n,\infty}\left(\theta(\cP^1_{\textnormal{bound}}(r)), \frac{2r}{(1 + \delta_n) n}\right)  \asymp \left\{ \begin{array}{ll}
	   \frac{r^2}{n}  & \text{ if } \delta_n \leq n^{-1/2}  \\
	    \delta_n^2 r^2 & \text{ if }  n^{-1/2} <  \delta_n \leq 1 \\
        r^2 & \text{ if } \delta_n > 1.
	\end{array} \right.
\end{equation} 
 
At the same time, we can also characterize the optimal rate of $R_{n,1}\left(\theta(\cP^1_{\textnormal{bound}}(r)), \beta_n\right)$ for all ranges of $\beta_n$.
\begin{Theorem} \label{th:average-case-bounded-class} Suppose $d = 1$. There exists a large universal constant $C > 0$ such that for any $n \geq C$, we have \footnote{Throughout the paper, the notation $R_{n,p}(\cdot, \beta_n) \asymp f(n,\beta_n)$ means that there are universal constants $0< c < C$ independent of $\beta_n$ such that $c f(n,\beta_n) \leq R_{n,p}(\cdot, \beta_n) \leq C f(n,\beta_n)$. }
\begin{equation*}
		R_{n,1}\left(\theta(\cP^1_{\textnormal{bound}}(r)), \beta_n \right)  \asymp   \left(\frac{\left( \frac{r}{n \beta_n} - 1 \right)_{+} }{1 + \left( \frac{r}{n \beta_n} - 1 \right)_{+}}\right)^3 r^2 \vee \frac{r^2}{n}, 
	\end{equation*}
	and one optimal estimator is
	{ \begin{equation} \label{eq:ave-est}
		\widehat{\theta} = \left\{ \begin{array}{ll}
			\bar{X}, & \delta_n < 0\\
			\left( 1- \delta_n \left(1 \wedge \frac{2(\sqrt{\delta_n} + 1/\sqrt{n})r}{|\bar{X}|} \right) \right) \bar{X}, & \delta_n \in [0,1]\\
            0, & \delta_n > 1,
		\end{array} \right.
	\end{equation} where $\delta_n = \frac{r}{n \beta_n} -1 $.}
\end{Theorem}
As we have already seen in Theorem \ref{th:bounded-worst-refined}, the exact phase transition location of $R_{n,1}\left(\theta(\cP^1_{\textnormal{bound}}(r)), \beta_n \right)$ is lower than the one of $R_{n,\infty}\left(\theta(\cP_{\textnormal{bound}}), \beta_n \right)$. In addition, given any $\delta_n \geq 0$, Theorem \ref{th:average-case-bounded-class} shows
\begin{equation} \label{eq:average-case-delta-dependence}
	 R_{n,1}\left(\theta(\cP^1_{\textnormal{bound}}(r)), \frac{r}{(1 + \delta_n) n}\right)  \asymp  \left\{ \begin{array}{ll}
	   \frac{r^2}{n}  & \text{ if } \delta_n \leq n^{-1/3}  \\
	    \delta_n^3 r^2 & \text{ if } n^{-1/3} <  \delta_n \leq 1 \\
        r^2 & \text{ if } \delta_n > 1.
	\end{array} \right.
\end{equation}

Compared with \eqref{eq:worst-bounded-refined}, we can see that not only does the location of phase transition change in $R_{n,1}\left(\theta(\cP^1_{\textnormal{bound}}(r)), \beta_n \right)$ comparing with $R_{n,\infty}\left(\theta(\cP^1_{\textnormal{bound}}(r)), \beta_n \right)$, the estimation error rate also becomes strictly faster, i.e., shifting from $\delta_n^2 r^2$ to $\delta_n^3 r^2$, near the phase transition location under the average-case stability constraint. A pictorial illustration on the difference of $R_{n,\infty}\left(\theta(\cP^1_{\textnormal{bound}}(r)), \beta_n\right)$ and $R_{n,1}\left(\theta(\cP^1_{\textnormal{bound}}(r)), \beta_n \right)$ is provided in Figure \ref{fig:phase-transition}(a).

\begin{figure}[h]
	\centering
	\includegraphics[width=0.57\textwidth]{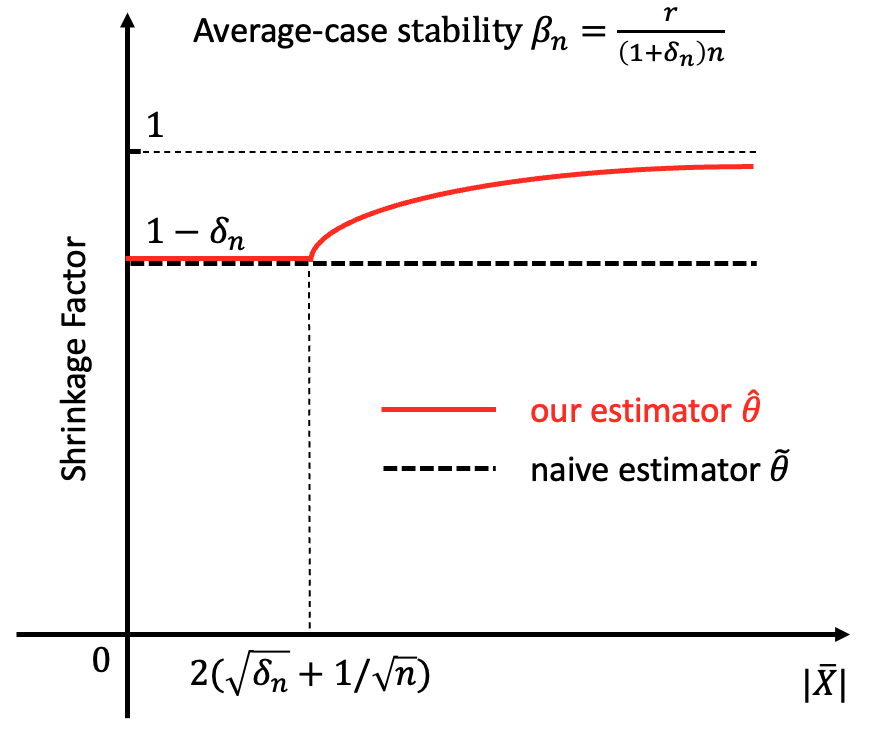}
	\caption{A pictorial illustration on the shrinkage factor of our estimator $\wh{\theta}$ \eqref{eq:ave-est}, and the naive estimator $\wt{\theta}$ \eqref{eq:naive-estimator} relative to the sample mean, given average-case stability $\beta_n = \frac{r}{(1+\delta_n)n}$ for some $\delta_n \in [0,1]$. } \label{fig:estimator}
\end{figure}

Finally, we would also like to comment on the novelty of the optimal estimators in Theorems \ref{th:worst-bounded-refined} and \ref{th:average-case-bounded-class}. Under worst-case stability, we first note that $\frac{\bar{X}}{1 + \frac{1}{\sqrt{n}} }$ is the estimator that can achieve the exact unconstrained minimax risk. To modify it to be stable, the basic idea is that if $\beta_n$ is greater than the worst-case stability level of $\frac{\bar{X}}{1 + \frac{1}{\sqrt{n}} }$, we just use $\frac{\bar{X}}{1 + \frac{1}{\sqrt{n}} }$, but if we require the estimator to be more stable than that, we can further shrink the sample mean until it achieves the desired stability level. This yields the estimator in Theorem \ref{th:worst-bounded-refined}. Under average-case stability, following a similar strategy, a naive extension of the estimator would be 
\begin{equation}\label{eq:naive-estimator}
    \wt{\theta} = \left\{ \begin{array}{ll}
			\bar{X}, & \delta_n < 0\\
			\left( 1- \delta_n\right) \bar{X}, & \delta_n \in [0,1]\\
            0, & \delta_n > 1
		\end{array} \right.
\end{equation}  where $\delta_n = \frac{r}{n \beta_n} -1 $. Then $\wt{\theta}$ can be shown to be $\beta_n$-average-case-stable. { We note $\wt{\theta}$ is the same as $\widehat{\theta}$ in \eqref{eq:ave-est} when $|\widebar{X}|$ is small.} However, in general, when $\beta_n = \frac{r}{(1+\delta_n)n}$ for some $\delta_n \in (n^{-1/3},1]$, the bias of $\wt{\theta}$ is of order $\delta_n^2 r^2$, which is bigger than the desired $\delta_n^3 r^2$ error rate. What we find, surprisingly, is that we can shrink the sample mean less aggressively when the magnitude of $|\bar{X}|$ is large while still maintaining the same level of average-case stability. As a result, the bias of the new estimator can be reduced to $\delta_n^3 r^2$. A pictorial illustration of the difference between $\wh{\theta}$ and $\wt{\theta}$ under average-case stability constraint is provided in Figure \ref{fig:estimator}. The main technical work here lies in checking that even for this data-driven shrinkage, the new estimator can still achieve the desired $\frac{r}{(1+\delta_n)n}$ average-case stability. { Readers are referred to the stability-guarantee part of the proof in Theorem \ref{th:average-case-bounded-class} and Lemma \ref{lm:averaged-case-lemma2} in the appendix for more details.}

\section{Heavy-tailed Mean Estimation}\label{sec:heavy-tailed-mean}
While the bounded setting considered in the previous section offers sharp theoretical results, in many practical settings the data does not lie in a bounded domain, and it is common to encounter heavy-tailed distributions \citep{catoni2012challenging}. In this section, we consider heavy-tailed mean estimation. Let $$\cP^{d}_{k}(r) = \{ \textnormal{distributions } P \textnormal{ supported on } \bbR^d \textnormal{ with } (\bbE_P[\|X\|_2^k])^{1/k} \leq r \}$$ be the class of distributions on $\bbR^d$ with $k$th moment bounded by $r^k$, where $k \geq 1$ is fixed. Without any stability constraint, when $k \geq 2$, the minimax risk for mean estimation in this family is of rate $\frac{r^2}{n}$ and the optimal estimator is the sample mean \citep{devroye2016sub}. When $k \in [1,2)$, the minimax risk for mean estimation in this family is of rate $ \frac{r^2}{n^{2(1-1/k)}}$ \citep{devroye2016sub} and this rate can be achieved by the truncated mean \citep{bubeck2013bandits}. At the same time, a trivial error rate is $r^2$ as $\|\bbE_P[X]\|_2 \leq r$ for any $P \in \cP^{d}_{k}(r)$. Under a stability constraint, we show that the stability-vs.-accuracy trade-off curve $R_{n,\infty}\left(\theta(\cP^{d}_{k}(r)), \beta_n\right)$ has a gradual phase transition, while $R_{n,1}\left(\theta(\cP^{d}_{k}(r)), \beta_n\right)$ still has a sharp phase transition. 
 \begin{Theorem} \label{th:worst-heavy-tail} 
 Given any $k \geq 1$,
 \begin{enumerate}
 	\item[(a)] Under worst-case stability (i.e., $p=\infty$), \begin{equation*}
		R_{n,\infty}\left(\theta(\cP^{d}_{k}(r)),  \beta_n \right) \asymp  \frac{r^{2k}}{(n\beta_n)^{2(k-1)}}\wedge r^2 + \frac{r^2}{n^{1 \wedge 2(1-1/k) }}. 
	\end{equation*} One optimal estimator is given by the truncated mean, $\wh{\theta} = \wb{Y}$, where $Y_i = X_i \indi(\|X_i\|_2 \leq \rho_k)$ for $i = 1, \ldots, n$ and where the truncation level is given by $\rho_k = \frac{n \beta_n}{2}\wedge rn^{1/k}$ for $k \in [1,2)$ and $\rho_k = \frac{n \beta_n}{2}$ for $k \geq 2$. 
	\item[(b)] Under average-case stability (i.e., $p=1$), \begin{equation*}
		R_{n,1}\left(\theta(\cP^{d}_{k}(r)),  \beta_n \right) \asymp \left\{ \begin{array}{ll}
			\frac{r^2}{n^{1 \wedge 2(1-1/k) }} & \textnormal{ if } \beta_n \geq \frac{24 r}{n}\\
			r^2 & \textnormal{ if } \beta_n \leq \frac{r}{10n}.
		\end{array} \right.
	\end{equation*}
One estimator that achieves the above rate is given by $\wh{\theta}= \left( \frac{\beta_n n}{24r} \wedge 1 \right) \wb{\theta}$ where 
	\begin{equation} \label{eq:heavy-tailed-aver-estimator}
		\wb{\theta} = \left(1 \wedge \frac{2nr}{\sum_{i=1}^n \|Y_i\|_2} \right) \bar{Y}.
	\end{equation}
	Again, $Y_i = X_i \indi(\|X_i\|_2 \leq \rho_k)$ for $i = 1, \ldots, n$, but now the truncation level is given by $\rho_k = rn^{1/k}$ for $k \in [1,2)$ and $\rho_k =\infty$ for $k \geq 2$.
 \end{enumerate}
	
\end{Theorem}
Comparing with Theorem \ref{th:worst-bounded}, we can see that the worst-case stability constraint has a more severe effect on heavy-tailed mean estimation. In particular, while $R_{n,\infty}\left(\theta(\cP^d_{\textnormal{bound}}(r)), \frac{2r}{n}\right) \asymp \frac{r^2}{n}$, now we instead have $R_{n,\infty}\left(\theta(\cP^{d}_{k}(r)), \frac{2r}{n}\right) \asymp r^2$. { The lower bound proof in Theorem \ref{th:worst-heavy-tail} Part (a) is based on Theorem \ref{th:worst-lower-bound-multi} with a carefully constructed coupling of $(\cD_n, \cD_n')$.} In terms of estimators, in contrast to the bounded class, here since the range of $X$ is unbounded, we have to truncate each data point to guarantee worst-case stability, and the optimal truncation level also depends on the stability requirement.

Comparing the Part (a) and Part (b) result in Theorem \ref{th:worst-heavy-tail}, we can see that $R_{n,1}\left(\theta(\cP^{d}_{k}(r)),  \beta_n \right)$ and $R_{n,\infty}\left(\theta(\cP^{d}_{k}(r)),  \beta_n \right)$ have a completely different phase transition behavior: the former one is sharp while the latter one is gradual (see a pictorial illustration in Figure \ref{fig:phase-transition}(b)). In particular, a $\frac{24r}{n}$-level average-case stability is enough to obtain the unconstrained minimax rate, while it is impossible to estimate the heavy-tailed mean under the same level worst-case stability constraint. This confirms the intuition that average-case stability is indeed a qualitatively much weaker constraint than the worst-case stability.

Finally, we would like to comment on the difference between the estimators under worst-case and average-case stability in Theorem \ref{th:worst-heavy-tail}. If we consider the sample mean, then 
\begin{equation*}
	\begin{split}
		\sup_{\cD_n \sim \cD_n'}\| \wh{\theta}(\cD_n) - \wh{\theta}(\cD_n') \|_2 \leq \sup \frac{\|X_i\|_2 + \|X_i'\|_2}{n}.
	\end{split}
\end{equation*} Thus to guarantee worst-case stability, we have to truncate each data point, with a truncation level that corresponds to the stability requirement. In contrast, under the average-case stability constraint, the sample mean has the following guarantee:
\begin{equation} \label{ineq:heavy-average-estimation-intuition}
	\begin{split}
		\frac{1}{(n+1)^2} \sup_{\cD_{n+1}} \sum_{1 \leq i,j\leq n+1} \| \wh{\theta}(\cD_{n+1}^{\setminus i}) - \wh{\theta}(\cD_{n+1}^{\setminus j}) \|_2 &\leq \frac{2}{(n+1)^2} \sup_{\cD_{n+1}} \sum_{1 \leq i<j\leq n+1} \frac{\|X_i\|_2 + \|X_j\|_2}{n} \\
		& = \frac{2}{(n+1)^2} \sup_{\cD_{n+1}} \left(\sum_{1 \leq i\leq n+1} \|X_i\|_2\right).
	\end{split}
\end{equation} Therefore, as long as $\sum_{1 \leq i\leq n+1} \|X_i\|_2$ is bounded, then we can control the average-case stability. That motivates the self-normalized shrinkage strategy in the estimator \eqref{eq:heavy-tailed-aver-estimator}, where we use the sample mean if $\sum_{i=1}^n \|Y_i\|_2$ is small, and normalize only if $\sum_{i=1}^n \|Y_i\|_2$ is large. { A rigorous proof on the average-case stability of estimator \eqref{eq:heavy-tailed-aver-estimator} is provided in Lemma \ref{lm:average-stability-bound} in the appendix.} Note that when $k \in [1,2)$, we still need to additionally truncate each data point, as that is needed for minimax optimality even when there is no stability constraint.

\section{Sparse Mean Estimation}\label{sec:sparse-mean}
In many modern applications, the dimension $d$ of the observed data may be larger than the sample size $n$, but a low-dimensional latent structure makes performing inference possible. In this section, we consider the task of sparse mean estimation under a stability constraint. Let $$\cP^d(r,s) = \{  \textnormal{distributions } P \textnormal{ supported on } \bbB_{\infty}(r) \subseteq \bbR^d \textnormal{ with } \|\bbE_P[X]\|_0 \leq s  \}$$ be the class of distributions with sparse mean vectors (here $\bbB_{\infty}(r)=\{x \in \bbR^d:\|x\|_\infty \leq r\}$  and $\|v\|_0 = \sum_{j=1}^d \indi\{v_j\neq 0\}$ denotes the sparsity, or $\ell_0$ `norm', of a vector $v\in\bbR^d$). Without the stability constraint, it is well known that estimators based on hard-thresholding or soft-thresholding the sample mean $\wb{X}$, i.e., 
\begin{equation} \label{eq:hard-soft-thres-est}
	\begin{split}
		\wh{\theta}^{\textnormal{hard}}_j = \wb{X}_j\indi(|\wb{X}|_j \geq \tau) \quad \textnormal{and}\quad \wh{\theta}^{\textnormal{soft}}_j = \textnormal{sign}(\wb{X}_j)(|\wb{X}_j| - \tau)_+ \quad \textnormal{ for } j = 1, \ldots, d,
	\end{split}
\end{equation} with $\tau \asymp \frac{r \sqrt{\log d}}{n} $ can yield the nearly optimal estimation error rate $ \frac{r^2 s \log d}{n}$ \citep{johnstone2017gaussian}. Here $\wb{X}_j$ denotes the $j$th coordinate of the sample mean. Under a stability constraint, we can show that the stability-vs.-accuracy trade-off curve $R_{n,p}\left(\theta(\cP^d(r,s)), \beta_n\right)$ has a sharp phase transition.
\begin{Theorem} \label{th:sparse-mean}
		Suppose $s \leq d^{1-c}$ for some $c \in (0,1]$ and $n \gtrsim \log d$. Given any $p \in[1,\infty]$,
		\begin{equation*}
		R_{n,p}\left(\theta(\cP^d(r,s)),  \beta_n \right) \asymp \left\{ \begin{array}{ll}
			\frac{r^2 s \log (4d)}{n} & \textnormal{ if } \beta_n \geq 4\sqrt 2\frac{ r \sqrt{s} }{n}\\
			r^2s & \textnormal{ if } \beta_n \leq \frac{1}{16}\frac{r \sqrt{s} }{n}.
		\end{array} \right.
	\end{equation*} One optimal estimator for achieving the above bound is given by a soft-thresholding estimator with a data-driven thresholding value:
	\begin{equation} \label{eq:data-driven-soft}
		\wh{\theta}_j = \left(1 \wedge\frac{n\beta_n}{4\sqrt{2s} r} \right)\textnormal{sign}(\wb{X}_j) (|\wb{X}_j| - |\wb{X}|_{(s+1)})_{+} \quad \textnormal{ for } j = 1, \ldots, d,
	\end{equation} where $|\wb{X}|_{(1)}\geq \dots \geq |\wb{X}|_{(d)}$ denote the order statistics of $\{|\wb{X}_i|\}_{i=1}^d$, so that $|\wb{X}|_{(s+1)}$ is the $(s+1)$-th largest (absolute) value in the sample mean.
\end{Theorem}
From Theorem \ref{th:sparse-mean}, we can see that for sparse mean estimation, the phase transition is sharp for $R_{n,p}\left(\theta(\cP^d(r,s)),  \beta_n \right)$, with a threshold $\beta_n^* = \frac{r\sqrt{s}}{n}$ that is independent of the ambient dimension $d$ and of the choice of $p$. 

In terms of the optimal estimator, our estimator is also novel and different from the estimators in \eqref{eq:hard-soft-thres-est} in an important way that we use a soft-thresholding estimator with a {\it data-driven} thresholding value, which is guaranteed to be $s$-sparse due to the particular thresholding we choose. Moreover, due to the Lipschitz property of the soft-thresholding operator, it turns out that such a sparsity-ensured soft-thresholding estimator can achieve the optimal stability we need. { The detailed analysis on the accuracy and stability of the proposed estimator are provided in Theorem \ref{thm:sparse-thresh-risk} and Lemma \ref{lem:sparse-thresh-stability}, respectively, in appendix.} In Appendix \ref{app:class-soft-hard}, we also illustrate that in general the estimators in \eqref{eq:hard-soft-thres-est} are not as stable as the one we propose. Interestingly, we note the same soft-thresholding operator as the one in \eqref{eq:data-driven-soft} is studied in Section 4.3 of \cite{liu2020between} in an optimization setting for the purpose of getting an optimal convergence rate. 

Finally, we comment that similar to bounded mean estimation setting, we also expect that the exact phase transition location in $R_{n,p}\left(\theta(\cP^d(r,s)),  \beta_n \right)$ would vary over $p$ and these
risks can also have different rates when $\beta_n$ is near the phase transition threshold. Narrowing down the exact phase transition location and identifying the quantitative difference between worst-case and average-case stability in the sparse mean estimation would likely require new lower bound and upper bound techniques, and we leave this question to future exploration.

\section{Nonparametric Regression Function Estimation} \label{sec:nonparametric}

While the previous three examples all study the problem of mean estimation in various settings, our last example departs from this question and considers the problem of function estimation in nonparametric regression. 

Suppose we have $n$ i.i.d.\ data points $(X_i, Y_i)_{i=1}^n$ generated from the following model:
\begin{equation}\label{eq:dgp-non-para}
	Y_i = f(X_i) + \xi_i,
\end{equation}
	where $f$ is some unknown function, $X_i \sim \textnormal{Unif}[0,1]$, and $\xi_i\sim N(0,\sigma^2)$. The goal is to estimate the function $f$ at a fixed point $x_0 \in (0,1)$, that is, to estimate $f(x_0)$, under a stability constraint. When there is no assumption on $f$, we cannot hope to estimate $f(x_0)$ well in general. Here, we assume $f \in \cB^s_{q_0,q_1}(1)$, where $\cB^s_{q_0,q_1}(1)$ denotes the unit Besov ball:
    \[\cB^s_{q_0,q_1}(1)  = \{ f \in L_2[0,1]: \|f\|_{\cB^s_{q_0,q_1}}\leq 1\}, \]
    where the Besov norm $\|\cdot \|_{\cB^s_{q_0,q_1}}$ will be defined below. In this section, we adopt the standard setting that $\nu := s - 1/q_0 > 0$, $1 \leq q_0 \leq \infty$ \citep{cai2003rates} and assume $1\leq q_1 \leq \infty$.

    Wavelet methods provide a natural framework for constructing both global and pointwise estimators in nonparametric problems, achieving optimal approximation rates through finite-resolution truncation (see, e.g.,~\cite{gine_mathematical_2016}).  
Formulating our analysis over Besov balls $\cB^s_{q_0,q_1}$ is convenient because they unify classical smoothness classes: Sobolev spaces correspond to $q_0=q_1=2$, while H\"older (Lipschitz) spaces correspond to $q_0=q_1=\infty$.  
Thus, our results simultaneously recover these familiar cases while also encompassing more refined smoothness regimes.

{ For concreteness, we take $\phi$ and $\psi$ to be, respectively, father (scaling) function and mother wavelet function from a compactly supported Daubechies dbN family \citet{daubechies1992ten} with smoothness order $N$ sufficiently large so that the wavelet system is $A$-regular with integer $A > s$. On $[0,1]$, we use the associated orthonormal Cohen--Daubechies--Vial boundary-corrected basis \citep{cohen1993wavelets}. For interior indices, its basis functions take the form
\[
\psi_{lk}(x) = 2^{l/2}\psi(2^l x - k),
\qquad
\phi_{l_0+1,m}(x) = 2^{(l_0+1)/2}\phi(2^{l_0+1}x - m).
\]
For indices whose supports intersect the endpoints of $[0,1]$, $\phi_{l_0+1,m}$ and $\psi_{lk}$ denote the corresponding Cohen--Daubechies--Vial (CDV) boundary-corrected functions.} { Here $l_0+1$ is the starting resolution level.} To use a single level index for the father and mother wavelets, define
\[
\mathcal K_l :=
\begin{cases}
    \{0,\ldots,2^{l_0+1}-1\}, & l=l_0,\\
    \{0,\ldots,2^l-1\}, & l\geq l_0+1.
\end{cases}
\]
With a mild abuse of notation, we write the father wavelets as $\psi_{l_0 k} := \phi_{l_0+1, k}$ for $k\in\mathcal K_{l_0}$, while $\psi_{lk}$ denotes the mother wavelets for $l\geq l_0+1$. We then represent any $f \in L_2[0,1]$ as
\[
f = \sum_{l=l_0}^{\infty}\sum_{k\in\mathcal K_l} f_{lk}\psi_{lk}, 
\qquad f_{lk} = \langle f,\psi_{lk}\rangle,
\]
so that $\|f\|_2^2 = \sum_{l,k} f_{lk}^2$ by orthonormality. With this representation of the function, we can define the Besov norm as
\[\|f \|_{\cB^s_{q_0,q_1}} = \begin{cases}
		\left( \underset{l=l_0}{\overset{\infty}{\sum}} \left(2^{l(s+1/2-1/q_0)} \left\| (f_{lk})_{k\in\mathcal K_l} \right\|_{q_0}\right)^{q_1}\right)^{1/q_1} &\text{ if } {1 \leq q_1 < \infty} \\
		\; \underset{l\geq l_0}{\sup} \; 2^{l(s+1/2-1/q_0)} \left\| (f_{lk})_{k\in\mathcal K_l} \right\|_{q_0}   &\text{ if } q_1 = \infty,
	\end{cases}\]
where $\|\cdot\|_{q_0}$ denotes the $\ell_{q_0}$ norm of a vector. Membership in $\cB^s_{q_0,q_1}(1)$ is thus characterized by the decay rate of the wavelet coefficients $\{f_{lk}\}$, with parameters $(s,q_0,q_1)$ controlling smoothness, integrability, and fine-scale behavior, respectively. 

{ In the conventional inhomogeneous Besov norm, the father-wavelet coefficient block is not multiplied by a resolution weight. Since the starting resolution $l_0$ is fixed, the factor $2^{l_0(s+1/2-1/q_0)}$ is a fixed positive constant, and hence the convention above yields a norm equivalent to the conventional Besov norm up to a scaling constant.}

A key property essential to our stability analysis is the \emph{local support} of wavelets: each $\psi_{lk}$ has compact support of size $O(2^{-l})$, ensuring that any single sample $(X_i, Y_i)$ influences only a bounded number of coefficients. This localization directly limits the estimator's sensitivity to perturbations in the data. Owing to these favorable structural properties, wavelet-based methods have also been widely employed in the design of differentially private estimators for nonparametric problems; see, for example, \cite{cai2024optimalfm,cai2024optimal}.

 In Appendix \ref{app:besov-intro}, we collect a few more basic properties of Besov spaces. 
    Interested readers are also referred to \cite{cohen1993multiresolution,johnstone2017gaussian} for a more detailed exploration of Besov space and wavelet basis.

When there is no stability constraint, it has been shown in \cite{cai2003rates} that the minimax pointwise $\ell_2$-risk for function estimation in Besov space is $n^{ \frac{-2\nu}{1 + 2\nu} }$. The rate-optimal estimator is based on the following truncated wavelet estimator:
\begin{equation} \label{eq:wavelet-est}
	\begin{split}
		 \wh{f}_{\cD_n}(x_0) = \frac{1}{n} \sum_{i=1}^n K_{L_{\text{opt}}}(X_i, x_0) Y_i
	\end{split}
\end{equation} where $ K_{L}(u, v) =\sum_{l=l_0}^{L} \sum_{k\in\mathcal K_l} \psi_{lk}(u) \psi_{lk}(v)$ is the projection kernel at resolution level $L$ and $L_{\text{opt}} = \lfloor \frac{\log_2 n}{2\nu+1} \rfloor$. { On the other hand, since $\nu=s-1/q_0>0$, the continuous Besov embedding $\cB^s_{q_0,q_1}([0,1])\hookrightarrow L_\infty([0,1])$ implies $\|f\|_\infty
    \leq M\|f\|_{\cB^s_{q_0,q_1}}
    \leq M$
uniformly over $f\in\cB^s_{q_0,q_1}(1)$; see, e.g., \citet[Proposition 4.3.9(i)]{gine_mathematical_2016}. Consequently, the constant estimator zero has pointwise squared risk of order at most $1$.} With the worst-case and average-case stability constraints, let us denote the $\ell_2$-pointwise risk for function estimation in Besov space as $R_{n,\infty}\left(\theta(\cB^s_{q_0,q_1}(1)), \beta_n\right)$ and $R_{n,1}\left(\theta(\cB^s_{q_0,q_1}(1)), \beta_n\right)$, respectively. The following theorem provides upper and lower bounds for $R_{n,\infty}\left(\theta(\cB^s_{q_0,q_1}(1)), \beta_n\right)$ and $R_{n,1}\left(\theta(\cB^s_{q_0,q_1}(1)), \beta_n\right)$. 

{
\begin{Theorem} \label{th:worst-nonparametric}
	\begin{enumerate}
		\item[(a)] Under worst-case stability ( i.e., $p=\infty$), \begin{equation*}
	 	\begin{split}
	 	R_{n,\infty}\left(\theta(\cB^s_{q_0,q_1}(1)), \beta_n\right)  \asymp (n \beta_n)^{-2\nu} \wedge 1 + n^{- \frac{2\nu}{2\nu + 1} }  .
	 	\end{split}
	 \end{equation*} One estimator achieving the upper bound is the transformed-clipped-truncated wavelet estimator:
	 \begin{equation} \label{eq:wavelet-est2}
	 \wh{f}_{\cD_n}(x_0) = h^{-1} \left( \pi_h \left( \frac{1}{n} \sum_{i=1}^n K_{L}(X_i, x_0) [Y_i]_{M}\right) \right) , 
	 \end{equation} where $h(t) = \bbE_{\xi \sim N(0, \sigma^2)}([\xi + t]_M)$ denotes a smooth and monotonically increasing transformation under Gaussian noise, $\pi_{h}(\cdot)$ denotes the Euclidean projection onto the interval $h([-M, M])$ and $L = \left\lfloor\log_2\left(c_{\psi}n \beta_n\right) \vee (\ell_0-1)\right \rfloor \wedge L_{\text{opt}}$. Here $M$ denotes the upper bound for $\|f\|_{\infty}$, $[\cdot]_M$ is the clipping operator at level $M$ (see its definition in Section \ref{sec:organization}), where $c_\psi > 0$ depending only on $(\sigma,s,q_0,q_1)$ and the wavelet basis. We take the convention that when $L = \ell_0-1$, then $ \wh{f}_{\cD_n}(x_0)  = 0$.
	 \item[(b)] Under average-case stability (i.e., $p=1$), {there exist constants $0<c_{\star}\leq C_{\star}<\infty$, depending only on $(\sigma,s,q_0,q_1)$ and the wavelet system, such that, for $n \geq 2,$} \begin{equation} \label{eq:wavelet-rate-ave}
	 	\begin{split}
	 		R_{n,1}\left(\theta(\cB^s_{q_0,q_1}(1)), \beta_n\right)  \asymp \left\{ \begin{array}{ll}
	 			n^{- \frac{2\nu}{2\nu + 1} } & \textnormal{ if } \beta_n \geq \frac{C_{\star}}{n} \\
	 		1 & \textnormal{ if } \beta_n \leq  \frac{c_{\star}}{n \log n}.
	 		\end{array} \right.
	 	\end{split}
	 \end{equation}
	 In addition, one estimator that achieves the above bound is
	 \begin{equation}\label{eq:wavelet-est-ave}
	 \wh{f}_{\cD_n}(x_0) =\left(\frac{n\beta_n }{24C_{\psi}} \wedge 1\right)  \left( \frac{2n C_{\psi}}{S} \wedge 1 \right) \bar{Z}  
	 \end{equation} where $Z_i = K_{L_{\text{opt}}}(X_i, x_0) Y_i$, $S = \sum_{i=1}^n |Z_i|$ and $C_{\psi} >0$ is some constant depending only on $(\sigma,s,q_0,q_1)$.
	\end{enumerate}
\end{Theorem}

 From Theorem \ref{th:worst-nonparametric}, we see that imposing worst-case stability has a dramatic impact on nonparametric function estimation. In particular, Part (a) shows that any estimator attaining the unconstrained minimax rate must have stability at least $\beta_n \asymp n^{- \frac{2\nu}{2\nu+1}}$. Without the transformation operation $h^{-1}(\pi_{h}(\cdot))$, the estimator in \eqref{eq:wavelet-est2} can be viewed as a clipped-truncation estimator. However, this estimator would suffer from large bias unless we allow the truncation level to grow as $\Theta(\log n)$ to accommodate the Gaussian tail. On the other hand, a growing truncation value leads to worse stability control. As a remedy, we consider $\sum_{i=1}^n K_{L}(X_i, x_0) [Y_i]_{M}$ as an estimator of the transformed function $h \circ f$ and then we apply an inverse transformation to transform it back. Furthermore, comparing with the classical wavelet estimator \eqref{eq:wavelet-est}, the resolution level of our estimator must depend on $\beta_n$ to ensure stability, which leads to a fundamentally different bias-variance trade-off. Under average-case stability, Part (b) characterizes two regimes: a sufficiently large stability budget for which the unconstrained minimax rate is attainable, and a sufficiently small stability budget for which only the trivial rate is possible. These two regimes demonstrate that average-case stability is a much weaker requirement: up to logarithmic factors, a stability budget of order $1/n$ suffices to achieve the optimal estimation rate, whereas the same stability level yields only trivial guarantees under worst-case stability. The corresponding estimator reflects this difference, requiring only selective shrinkage when the sum of the projection response $|Z_i|$ has a large magnitude, without the need for individual response clipping.}

\section{Discussion and Future Directions} \label{sec:discussion}
To conclude, we examine some connections to the literature, including the relationship between stability and privacy, and a comparison to alternative formulations of the stability constraint, and discuss some open questions raised by this work.

\subsection{Conversion between Stability and Differential Privacy} \label{sec:stability-dp}
In this section, we show that by leveraging the connections between worst-case stability and privacy, we obtain some interesting implications on the stability-vs.-accuracy trade-off curve and privacy-vs.-accuracy trade-off curve. 

To illustrate the idea clearly, we focus on the pure differential privacy (DP) in the setting of a $1$-dimensional parameter estimation problem and assume that the estimator $\wh{\theta}(\cD_n)$ is deterministic. Recall from \cite{dwork2006calibrating} that a randomized algorithm $M: \cX^n \to \bbR$ is $\epsilon$-differentially private if for every pair of adjacency datasets $\cD_n,\cD_n' \in \cX^n $ that differ by one individual data and every measurable $S \subseteq \bbR$, 
\begin{equation*}
	\bbP(M(\cD_n) \in S) \leq e^{\epsilon} \bbP(M(\cD_n') \in S)
\end{equation*} where the probability measure $\bbP$ is taken with respect to the randomness of $M$ only. The Laplace noise mechanism enables us to convert a worst-case-stable estimator to a differentially private estimator:
\begin{Lemma}[\cite{dwork2006calibrating}] \label{lm:stability-to-privacy}
	If $\wh{\theta}(\cD_n)$ is a $\beta_n$-worst-case-stable estimator, then $\wh{\theta}(\cD_n) + \xi$, where $\xi$ is independent of $\wh{\theta}(\cD_n)$ and follows $\textnormal{Laplace}(\beta_n/ \epsilon_n)$, satisfies $\epsilon_n$-differential privacy.
\end{Lemma}

In fact, we can also convert a differentially private estimator to a stable one if we have a mild assumption about the differentially private estimator.
\begin{Lemma}[Modified from Lemma 8 of \cite{wang2016learning}] \label{lm:privacy-to-stability}
	Suppose $M(\cD_n)$ satisfies $\epsilon_n$-differential privacy and $|M(\cD_n)| \leq r$ a.s. for every $\cD_n$, then $\wh{\theta}(\cD_n) = \bbE[M(\cD_n)|\cD_n]$ satisfies $r(e^{\epsilon_n}-1)$-worst-case stability. When $\epsilon_n < 1$, then $\wh{\theta}(\cD_n)$ is $2r\epsilon_n$-worst-case-stable.
\end{Lemma}
Similar to the development of the stability-constrained minimax risk in Section \ref{sec:formulation}, here we will also consider the privacy-constrained minimax risk \citep{barber2014privacy}:
 $$R_{\DP}( \theta(\cP) , \epsilon_n ) := \inf_{M \textnormal{ satisfies } \epsilon_n-\DP } \sup_{P \in \cP} \bbE [ M(\cD_n) - \theta(P) ]^2.$$
Then Lemma \ref{lm:stability-to-privacy} and \ref{lm:privacy-to-stability} together yield a few interesting conversions between $R_{\DP}( \theta(\cP) , \epsilon_n )$ and $R_{n,\infty}(\theta(\cP), \beta_n)$.

\begin{Proposition} \label{prop:convert-privacy-stability-curve}
	In the $1$-dimensional parameter estimation setting with $n$ i.i.d. samples, suppose the parameter of interest lies in $[-r,r]$. Then the stability-vs.-accuracy trade-off curve $R_{n,\infty}(\cdot)$ may be used to provide upper and lower bounds for privacy-vs.-accuracy trade-off curve $R_{\textnormal{DP}}(\cdot)$, as follows:
    \begin{equation}\label{eq:stability-dp}
        R_{n,\infty}(\theta(\cP), r(e^{\epsilon_n} - 1))\leq R_{\textnormal{DP}}( \theta(\cP) , \epsilon_n ) \leq \inf_{\beta_n} \left( R_{n,\infty}(\theta(\cP), \beta_n) + \frac{2 \beta_n^2 }{ \epsilon_n^2 } \right).
    \end{equation}
Equivalently, by inverting the above bounds, the privacy-vs.-accuracy trade-off curve $R_{\textnormal{DP}}(\cdot)$ may be used to provide upper and lower bounds on the stability-vs.-accuracy trade-off curve $R_{n,\infty}(\cdot)$, as follows:
\begin{equation}\label{eq:dp-stability}
    \sup_{\epsilon_n} \left(R_{\textnormal{DP}}( \theta(\cP) , \epsilon_n ) - \frac{2 \beta_n^2 }{ \epsilon_n^2 } \right) \leq R_{n,\infty}(\theta(\cP), \beta_n) \leq R_{\textnormal{DP}} \left(\theta(\cP),  \log\left( 1 + \frac{\beta_n}{r} \right)\right).
\end{equation}
\end{Proposition}

Proposition \ref{prop:convert-privacy-stability-curve} reveals that the stability-vs.-accuracy trade-off curve and the privacy-vs.-accuracy trade-off curve are closely related, and each may be used to derive bounds on the other. Moreover, we will illustrate in Appendix \ref{app:stability-privacy} using mean estimation examples to show a few implications of these inequalities. { Interestingly, the first inequality in \eqref{eq:dp-stability} provides an alternative approach to prove a lower bound under a worst-case stability constraint. On the other hand, we emphasize that our approach in Section \ref{sec:lower-bound} is still useful as it is more direct and can provide sharper bounds (in addition to generalizing to  $p\in[1,\infty]$, allowing for average-case in addition to worst-case stability).}

\subsection{A Comparison of Different Stability Notions} \label{sec:compare-diff-stability}
In the literature, there are many other notions of stability, which are developed for different purposes. One advantage of our framework is that it provides a unified way to compare different notions of stability from a statistical decision theory point of view. Other than the worst-case stability, average-case stability, and in general the $\ell_p$-stability we have discussed before, here we mention a few other common stability notions. For example, in many statistical applications, a type of { distributional-}average-case stability, where the average is over the randomness in the data, has also been considered:  we might say $\wh{\theta}$ is $\beta_n(P)$-distributional-stable with respect to $P$ if 
\begin{equation} \label{eq:average-randomness-def}
	\begin{split}
		\sup_{i,j \in [n+1]}\bbE_{\cD_{n+1}\overset{i.i.d.}\sim P; \xi}\|\wh{\theta}(\cD_{n+1}^{\setminus i}; \xi) - \wh{\theta}(\cD_{n+1}^{\setminus j}; \xi)\|_2 \leq \beta_n(P)
	\end{split}
\end{equation} for data $\cD_{n+1}$ drawn i.i.d.\ from $P$ (rather than for an arbitrary dataset, as in the definitions used throughout this paper). One application of such a stability condition is the problem of inference on cross-validation estimate of model risk \citep{austern2020asymptotics,bayle2020cross,kissel2022high}.  However, it is known that testing whether $\beta_n(P)$ is small or not is impossible \citep{kim2021black,luo2024algorithmic}. { For any symmetric estimator, it is easy to check that if $\wh{\theta}$ is $\beta_n$-$\ell_1$-stable, then it is also $\beta_n$-distributional-stable with respect to any $P$. Conceptually, distributional-average-case stability is distribution dependent, so to study it in a decision-theoretic framework, we can require $\wh{\theta}$ to satisfy the distributional-average-case stability in \eqref{eq:average-randomness-def} uniformly over all distributions $P \in \cP$ for some distribution class $\cP$ and ask what the best possible risk we can achieve over this distribution class, i.e.,
\begin{equation} \label{eq:distributional-risk}
	\begin{split}
	R_{n,\textnormal{dis}}(\theta(\cP), \beta_n) := \inf_{\wh{\theta} \textnormal{ is $\beta_n$-distributional-stable} } \sup_{P \in \cP} \bbE_{\cD_n, \xi} \|  \wh{\theta}(\cD_n; \xi) - \theta(P) \|_2^2.
	\end{split}
\end{equation} 
 A general minimax lower bound under the distributional-average-case stability constraint is provided as follows. 
	\begin{Theorem} \label{th:dis-risk}
		For any $P_0, P_1 \in \cP$ such that $P_t:=(1-t) P_0 + tP_1 \in \cP$ for all $t \in [0,1]$,
		\begin{equation*}
				R_{n,\textnormal{dis}}\left(\theta(\cP),  \beta_n \right) \geq R_n(\theta(\cP))\vee  \sup_{0<\eta<1/2} \frac{1}{4}\left( \|\theta(P_\eta) - \theta(P_{1-{\eta}})\|_2 - n \beta_n \log \left(\frac{1-\eta}{\eta}\right)  \right)^2_{+}.
		\end{equation*} In particular, if $\theta(P)$ is linear in $P$, then 
		\begin{equation*}
			R_{n,\textnormal{dis}}\left(\theta(\cP),  \beta_n \right) \geq R_n(\theta(\cP))\vee \frac{1}{16}\left( \|\theta(P_0) - \theta(P_1)\|_2 - \log (9)  n \beta_n\right)^2_{+}.
		\end{equation*}
	\end{Theorem}
    When $\theta(P)$ is linear in $P$, we can see that the above Theorem \ref{th:dis-risk} and Corollary \ref{coro:average-case-stability} provide almost the same lower bound under distributional average-case stability and $\ell_1$-stability, which demonstrates the similarity of these two stability notions.
    }

Finally, we note that for all the stability notions we have considered so far, the measure of the closeness of the outputs given a perturbation in the inputs is taken with respect to a metric on the parameter or function space. Such a notion has a long history as a tool for ensuring generalization \citep{rogers1978finite,devroye1979distribution,kearns1997algorithmic,bousquet2002stability}. More recently, lots of progress has been made on another class of distributional algorithmic stability, where the output object is a distribution and the metric of closeness is based on distributional closeness. Broadly speaking, this includes differential privacy \citep{dwork2006calibrating} (as discussed in Section \ref{sec:stability-dp}), perfect generalization \citep{cummings2016adaptive}, replicability \citep{impagliazzo2022reproducibility},  KL-stability and TV-stability \citep{bassily2016algorithmic}. We refer readers to \cite{bun2023stability,kalavasis2023statistical,moran2023bayesian} for a recent summary of these notions and their connections.

\subsection{Additional Directions and Open Questions}
The stability property on its own demonstrates some robustness of the statistical procedure to data perturbation. In the robust statistics literature, one related notion is the influence function \citep{huber2011robust, hampel2011robust}: given any functional $T$ of a distribution $F$ defined on $\bbR$, the influence curve $\textnormal{IF}(\cdot; T, F): \cX \to \bbR$ at $x$ is defined by 
\begin{equation*}
	\begin{split}
		\textnormal{IF}(x; T, F) = \lim_{t \to 0}\frac{T((1-t) F + t \delta_x ) - T(F)}{t},
	\end{split}
\end{equation*}  where $\delta_x$ denotes the point mass at $x$. The maximum absolute value of the influence curve, i.e., $\gamma^*(T, F) = \sup_x |\textnormal{IF}(x; T, F)|$, is called the gross error sensitivity \citep{hampel2011robust}.
If the gross error sensitivity is unbounded, then outliers may be problematic. To estimate the gross error sensitivity, we can use the following empirical version:
\begin{equation*}
	\begin{split}
		\gamma^*(T; \{x_1, \ldots, x_n\}) = \sup_{x} \left| (n+1) \left(T\left((1 - \frac{1}{n+1}) F_n + \frac{1}{n+1} \delta_x \right)- T(F_n) \right) \right|,
	\end{split}
\end{equation*} where $F_n:= \frac{1}{n} \sum_{i=1}^n \delta_{x_i} $ denotes the empirical distribution of the data. We note that this quantity can be viewed as a ``local'' version of the worst-case stability, as it fixes the dataset $\{x_1, \ldots, x_n\}$ in computing the worst-case perturbation of replacing one data point. In robust statistics, a number of robust estimators have been developed to offer control over the sensitivity. It is interesting to consider whether some of these ideas may be useful for deriving a more principled approach to develop stable estimators. 

{ Finally, in the two-regime sense of Definition \ref{def:phase-transition-pattern}, all four estimation problems explored in this paper exhibit a sharp phase transition under average-case stability, while under worst-case stability the phase transition may be sharp or gradual, depending on the setting (see Table \ref{tab:phase-transition}). This statement summarizes the behavior in the sufficiently small and sufficiently large stability-budget regimes and does not assert a complete characterization of the trade-off curve at every stability level; in particular, the intermediate regime for nonparametric regression under average-case stability remains open. This suggests an interesting question of whether the two-regime pattern is a common phenomenon or whether there is an example exhibiting a gradual phase transition even under average-case stability.}

\section*{Acknowledgments} We thank the Associate Editor and the two anonymous referees for their helpful suggestions, which helped improve the presentation and quality of this paper. R.F.B.\ was supported by the Office of Naval Research via grant N00014-24-1-2544
and by the National Science Foundation via grant DMS-2023109. Y. L. would like to thank Ankit Pensia for helpful discussions. The authors thank the American Institute of Mathematics (AIM) for hosting the workshop ``Algorithmic stability: mathematical foundations for the modern era'' in 2025, which helped us to develop this work.

\bibliographystyle{apalike}
\bibliography{reference}
%
%
		\appendix

\section{Connection of Two Average-case Stability Notions} \label{sec:two-average-case-stability}
If $\wh{\theta}$ is also well defined on a dataset of size $n+1$, then another average-case stability notion based on dropping one data point at random would read as
\begin{equation}\label{eq:average-case-stability2}
	\frac{1}{n+1} \sup_{\cD_{n+1} \in \cX^{n+1}} \sum_{1 \leq i \leq n+1} \bbE_{\xi} \|\wh{\theta}(\cD_{n+1}; \xi )  - \wh{\theta}(\cD_{n+1}^{\setminus i}; \xi )\|_2 \leq \beta_n. 
\end{equation} We first show that the version in Definition \ref{def:average-case-stability} can be implied by \eqref{eq:average-case-stability2}.
\begin{equation*}
	\begin{split}
		& \frac{1}{(n+1)^2} \sup_{\cD_{n+1} \in \cX^{n+1}} \sum_{1 \leq i, j \leq n+1} \bbE_{\xi} \|\wh{\theta}(\cD_{n+1}^{\setminus i}; \xi )  - \wh{\theta}(\cD_{n+1}^{\setminus j}; \xi )\|_2 \\
		= & \frac{1}{(n+1)^2} \sup_{\cD_{n+1} \in \cX^{n+1}} \sum_{1 \leq i, j \leq n+1} \bbE_{\xi} \|\wh{\theta}(\cD_{n+1}^{\setminus i}; \xi )  - \wh{\theta}(\cD_{n+1}; \xi ) + \wh{\theta}(\cD_{n+1}; \xi )  - \wh{\theta}(\cD_{n+1}^{\setminus j}; \xi )\|_2 \\
		\leq & \frac{1}{(n+1)^2} \sup_{\cD_{n+1} \in \cX^{n+1}} \sum_{1 \leq i, j \leq n+1} \left( \bbE_{\xi} \|\wh{\theta}(\cD_{n+1}^{\setminus i}; \xi )  - \wh{\theta}(\cD_{n+1}; \xi )\|_2 + \bbE_{\xi} \| \wh{\theta}(\cD_{n+1}; \xi )  - \wh{\theta}(\cD_{n+1}^{\setminus j}; \xi )\|_2 \right) \\
		= & \frac{2}{(n+1)^2} \sup_{\cD_{n+1} \in \cX^{n+1}} \sum_{1 \leq i, j \leq n+1}  \bbE_{\xi} \| \wh{\theta}(\cD_{n+1}; \xi )  - \wh{\theta}(\cD_{n+1}^{\setminus j}; \xi )\|_2 \\
		\overset{ \eqref{eq:average-case-stability2} }\leq & 2 \beta_n.
	\end{split}
\end{equation*} Next, we show that the average-case stability defined in Definition \ref{def:average-case-stability} can also imply the one in \eqref{eq:average-case-stability2} if we define $\wh{\theta}$ properly on a dataset of size $n+1$. Suppose $\wh{\theta}$ is already defined on a dataset of size $n$, then one natural estimator on a dataset of size $n+1$ would be $\wh{\theta}(\cD_{n+1}; \xi) = \frac{1}{n+1} \sum_{j=1}^{n+1} \wh{\theta} (\cD_{n+1}^{\setminus j}; \xi) $. With this definition of $\wh{\theta}(\cD_{n+1};\xi)$, we have
\begin{equation*}
	\begin{split}
		&\frac{1}{n+1} \sup_{\cD_{n+1} \in \cX^{n+1}} \sum_{1 \leq i \leq n+1} \bbE_{\xi} \|\wh{\theta}(\cD_{n+1}; \xi )  - \wh{\theta}(\cD_{n+1}^{\setminus i}; \xi )\|_2 \\
		=& \frac{1}{n+1} \sup_{\cD_{n+1} \in \cX^{n+1}} \sum_{1 \leq i \leq n+1} \bbE_{\xi} \left\| \frac{1}{n+1} \sum_{j=1}^{n+1} \left(\wh{\theta} (\cD_{n+1}^{\setminus j};\xi) - \wh{\theta}(\cD_{n+1}^{\setminus i}; \xi ) \right) \right\|_2 \\
		\leq & \frac{1}{(n+1)^2} \sup_{\cD_{n+1} \in \cX^{n+1}} \sum_{1 \leq i, j \leq n+1} \bbE_{\xi} \|\wh{\theta}(\cD_{n+1}^{\setminus i}; \xi )  - \wh{\theta}(\cD_{n+1}^{\setminus j}; \xi )\|_2\\
		&  \overset{\textnormal{Definition } \ref{def:average-case-stability} }\leq  \beta_n. 
	\end{split}
\end{equation*}

\section{Stability of Standard Soft- and Hard-Thresholding Estimators}
\label{app:class-soft-hard}

We now compute the worst-case stability for the classical soft- and hard-thresholding estimators introduced in~\eqref{eq:hard-soft-thres-est}.  
Throughout, we assume $X_i \in \bbB_\infty(r) \subseteq \bbR^d$ for some $r>0$, and denote the sample mean by
\[
\wb{X} = \frac{1}{n}\sum_{i=1}^n X_i,
\qquad
\wb{X}_j = \frac{1}{n}\sum_{i=1}^n X_{ij}.
\]

\subsection{Soft-thresholding}
The soft-thresholding estimator is defined coordinatewise as
\[
\wh{\theta}^{\textnormal{soft}}_j
\;=\;
\textnormal{sign}(\wb{X}_j)\,\bigl(|\wb{X}_j| - \tau\bigr)_+,
\qquad j=1,\dots,d,
\]
for some threshold $\tau>0$.  

Consider replacing one sample $X_i$ with $X_i'$.  
For any coordinate $j$,
\[
\bigl|\wb{X}_j - \wb{X}_j'\bigr|
= \frac{1}{n}\,\bigl|X_{ij}-X_{ij}'\bigr|
\le \frac{2r}{n}.
\]
The soft-thresholding operator $x\mapsto \textnormal{sign}(x)(|x|-\tau)_+$ is $1$-Lipschitz, so this perturbation translates directly to the estimate:
\[
\bigl|\wh{\theta}^{\textnormal{soft}}_{j}(\cD_n) - \wh{\theta}^{\textnormal{soft}}_{j}(\cD_n')\bigr|
\le \frac{2r}{n}.
\]
In the worst case, the equality can be achieved. Aggregating across all coordinates,
\[
\bigl\|\wh{\theta}^{\textnormal{soft}}(\cD_n) - \wh{\theta}^{\textnormal{soft}}(\cD_n')\bigr\|_2
\le \sqrt{d}\,\frac{2r}{n}.
\]
Hence, the worst-case stability of $\wh{\theta}^{\textnormal{soft}}$ is bounded by
\[
\beta_n^{\textnormal{soft}}
\;\le\;
\frac{2r\sqrt{d}}{n},
\] and again the equality can be achieved in the worst case.

\subsection{Hard-thresholding}
The hard-thresholding estimator is defined as
\[
\wh{\theta}^{\textnormal{hard}}_j
\;=\;
\wb{X}_j\,\mathbf{1}\{|\wb{X}_j|\ge\tau\},
\qquad j=1,\dots,d.
\]
Replacing one sample again perturbs $\wb{X}_j$ by at most $2r/n$.  
When $|\wb{X}_j|$ lies far from the threshold $\tau$, this change only induces a perturbation of order $O(r/n)$.  
However, if $|\wb{X}_j|$ lies within $O(r/n)$ of $\tau$, a single replacement may flip the indicator $\mathbf{1}\{|\wb{X}_j|\ge\tau\}$, causing a discontinuous jump as large as $O(\tau)$ in that coordinate.  
Consequently,
\[
\beta_n^{\textnormal{hard}}
\;\le\;
\left(\tau + \frac{2r}{n}\right)\sqrt{d}.
\]
When $\tau \asymp \frac{r \sqrt{\log d}}{n} $, the hard-thresholding estimator is $C \frac{r\sqrt{d \log d}}{n}$-worst-case-stable for some constant $C>0$.

\subsection{Comparing to our proposed estimator}
From the above two subsections, we can see that the soft and hard-thresholding estimators in \eqref{eq:hard-soft-thres-est} have worst-case stability that scales at least as large as $\frac{r\sqrt{d}}{n}$ and it depends on the ambient dimension $d$. While as we illustrated in the proof of Theorem \ref{th:sparse-mean}, specifically Lemma \ref{lem:sparse-thresh-stability}, our estimator is $\frac{4 \sqrt{2} r \sqrt{s} }{n}$-worst-case-stable. The stability of our estimator only scales with the intrinsic sparsity level, not with the ambient dimension.

\section{Conversion Between Stability and Differential Privacy} \label{app:stability-privacy}
In this section, we apply Proposition \ref{prop:convert-privacy-stability-curve} on two concrete examples: one dimensional mean estimation in a bounded class and heavy-tailed mean estimation. We will use these two examples to show that the reduction from stability to privacy, i.e., the upper bound in \eqref{eq:stability-dp} and the lower bound in \eqref{eq:dp-stability}, are often tight, but the reduction from privacy to stability, i.e., the lower bound in \eqref{eq:stability-dp} and the upper bound in \eqref{eq:dp-stability}, are often loose.

First, the minimax rates for mean estimation under DP constraint in bounded class and heavy-tailed class with $k \geq 2$ have been established in \cite{barber2014privacy}:
\begin{equation} \label{eq:privacy-bounds}
\begin{split}
	R_{\DP}( \theta(\cP^1_{\textnormal{bound}}(r)) , \epsilon_n ) & \asymp \frac{r^2}{n} + \left(\frac{r^2}{n^2 \epsilon_n^2} \wedge r^2\right),\\
	R_{\DP}( \theta(\cP^1_k(r)) , \epsilon_n ) & \asymp \frac{r^2}{n} + \left(\frac{r^{2}}{(n \epsilon_n)^{2-2/k}}\wedge r^2\right).
\end{split}
\end{equation} The minimax rates for mean estimation in bounded and heavy-tailed classes under the worst-case-stability constraint have been established in Theorem \ref{th:worst-bounded-refined} and Theorem \ref{th:worst-heavy-tail}, respectively. Then a corollary of Proposition \ref{prop:convert-privacy-stability-curve} is given as follows and its proof is provided in the subsequent subsections.
\begin{Corollary} \label{coro:bounded}
	Suppose $d = 1$. Then a combination of Proposition \ref{prop:convert-privacy-stability-curve} and Theorem \ref{th:worst-bounded-refined} yields
	\begin{subequations}
		\begin{align}
			& R_{\textnormal{DP}}( \theta(\cP^1_{\textnormal{bound}}(r)) , \epsilon_n ) \lesssim \frac{r^2}{n} + \left(\frac{r^2}{n^2 \epsilon_n^2} \wedge r^2\right), \label{eq:dp-upper-bounded}\\
			& R_{\textnormal{DP}}( \theta(\cP^1_{\textnormal{bound}}(r)) , \epsilon_n ) \gtrsim  \frac{r^2}{n} + \left(\frac{\left( \frac{2}{n(e^{\epsilon_n}-1)} - 1 \right)_{+}}{1 + \left( \frac{2}{n (e^{\epsilon_n}-1)} - 1 \right)_{+} } \right)^2 r^2. \label{eq:dp-lower-bounded}
		\end{align}
	\end{subequations}
	
	At the same time, a combination of Proposition \ref{prop:convert-privacy-stability-curve} and \eqref{eq:privacy-bounds} yields
	\begin{subequations}
		\begin{align}
		R_{n,\infty}(\theta(\cP^1_{\textnormal{bound}}(r)), \beta_n) &\gtrsim \frac{r^2}{n} + \left( c_0 r^2 - 2 \beta_n^2 n^2 \right)_{+}  \label{eq:stability-lower-bounded}\\
			R_{n,\infty}(\theta(\cP^1_{\textnormal{bound}}(r)), \beta_n) &\lesssim \frac{r^2}{n} + \left(\frac{r^2}{n^2  \log^2\left( 1 + \frac{\beta_n}{r} \right)} \wedge r^2\right). \label{eq:stability-upper-bounded}
		\end{align}
	\end{subequations}
for some universal constant $c_0 > 0$.  
\end{Corollary}

We can see that the upper bound in \eqref{eq:dp-upper-bounded} and the lower bound in \eqref{eq:stability-lower-bounded} are tight. These two implications are deduced from the reduction from stability to privacy. On the other hand, the results coming from the reduction of privacy to stability, i.e., \eqref{eq:dp-lower-bounded} and \eqref{eq:stability-upper-bounded}, are loose. For example, in $\eqref{eq:dp-lower-bounded}$, if we take $\epsilon_n = \frac{2}{n}$, then we can only get $r^2/n$ lower bound, but the tight lower bound would be $r^2$ based on \eqref{eq:privacy-bounds}. We have a similar result for the heavy-tailed mean estimation.

\begin{Corollary} \label{coro:heavy-tailed}
	Suppose $d = 1$ and $k \geq 2$. Then a combination of Proposition \ref{prop:convert-privacy-stability-curve} and Theorem \ref{th:worst-heavy-tail} Part (a) yields
	\begin{subequations}
		\begin{align}
				& R_{\textnormal{DP}}( \theta(\cP^1_k(r)) , \epsilon_n ) \lesssim \frac{r^2}{n} + \left(\frac{r^2}{(n \epsilon_n)^{2-2/k}}\wedge r^2\right), \label{eq:dp-upper-heavy}\\
			& R_{\textnormal{DP}}( \theta(\cP^1_k(r)) , \epsilon_n ) \gtrsim  	\frac{r^2}{n} + \left( \frac{r^2}{(n(e^{\epsilon_n}-1))^{2(k-1)}} \right) \wedge r^2.  \label{eq:dp-lower-heavy}
		\end{align}
	\end{subequations}
	
	At the same time, a combination of Proposition \ref{prop:convert-privacy-stability-curve} and \eqref{eq:privacy-bounds} yields
	\begin{subequations}
		\begin{align}
				R_{n,\infty}(\theta(\cP^1_k(r)), \beta_n) &\gtrsim \frac{r^2}{n} + \left( \frac{r^{2k}}{(n\beta_n)^{2(k-1)}} \wedge r^2\right) , \label{eq:stability-lower-heavy}\\
				R_{n,\infty}(\theta(\cP^1_k(r)), \beta_n) &\lesssim \frac{r^2}{n} + \left(\frac{r^2}{ \left(n  \log\left( 1 + \frac{\beta_n}{r} \right) \right)^{2-2/k}} \wedge r^2\right). \label{eq:stability-upper-heavy}
		\end{align}
	\end{subequations}	
\end{Corollary} 
Again, we can see that \eqref{eq:dp-upper-heavy} and \eqref{eq:stability-lower-heavy} are tight, but \eqref{eq:dp-lower-heavy} and \eqref{eq:stability-upper-heavy} are loose.

\subsection{Proof of Corollary \ref{coro:bounded}}
First, we note that \eqref{eq:dp-lower-bounded} and \eqref{eq:stability-upper-bounded} are straightforward. Next, we show \eqref{eq:dp-upper-bounded}:	\begin{equation*}
		\begin{split}
			R_{\textnormal{DP}}( \theta(\cP^1_{\textnormal{bound}}(r)) , \epsilon_n ) &\overset{ \eqref{eq:stability-dp} }\leq \inf_{\beta_n} \left( R_{n,\infty}(\theta(\cP^1_{\textnormal{bound}}(r)), \beta_n) + \frac{2 \beta_n^2 }{ \epsilon_n^2 } \right) \\
			& \overset{ \textnormal{Theorem } \ref{th:worst-bounded-refined} }= \inf_{\beta_n} \left( \left(\frac{\left( \frac{2r}{n \beta_n} - 1 \right)_{+}}{1 + \left( \frac{2r}{n \beta_n} - 1 \right)_{+} } \right)^2 r^2 \vee \frac{r^2}{\left(\sqrt{n} + 1 \right)^2} +   \frac{2 \beta_n^2 }{ \epsilon_n^2 } \right)\\
			& \lesssim  \frac{r^2}{n} + \frac{r^2}{n^2 \epsilon_n^2}\wedge r^2. \quad\left( \textnormal{take } \beta_n = \left\{\begin{array}{ll}
				2r/n & \epsilon_n \geq 2/n\\
				r\epsilon_n & \epsilon_n < 2/n
			\end{array} \right.\right)
		\end{split}
	\end{equation*} 

Finally, we show \eqref{eq:stability-lower-bounded}. By \eqref{eq:privacy-bounds}, we know there exists some universal constant $c_0$ such that $R_{\DP}( \theta(\cP^1_{\textnormal{bound}}(r)) , \epsilon_n ) \geq c_0 \left( \frac{r^2}{n} + \left(\frac{r^2}{n^2 \epsilon_n^2} \wedge r^2\right) \right)$. Then by the lower bound in \eqref{eq:dp-stability}, we have
\begin{equation}\label{eq:trivial-lower-bound}
\begin{split}
     R_{n,\infty}(\theta(\cP^1_{\textnormal{bound}}(r)), \beta_n) &\geq \sup_{\epsilon_n} \left( c_0 \left( \frac{r^2}{n} + \left(\frac{r^2}{n^2 \epsilon_n^2} \wedge r^2\right) \right) - \frac{2\beta_n^2}{\epsilon_n^2} \right) \\
     & = \max\left\{ \frac{c_0 r^2}{n} + c_0 r^2 - 2 \beta_n^2 n^2, \sup_{\epsilon_n \geq 1/n} \left(\frac{c_0 r^2}{n} + \frac{1}{\epsilon_n^2 n^2} \left( c_0 r^2 - 2 \beta_n^2 n^2 \right)  \right) \right\}\\
     & \geq \frac{c_0 r^2}{n} + \left( c_0 r^2 - 2 \beta_n^2 n^2 \right)_{+}.
\end{split}   
\end{equation} 
\subsection{Proof of Corollary \ref{coro:heavy-tailed} }
First, we note that \eqref{eq:dp-lower-heavy} and \eqref{eq:stability-upper-heavy} are straightforward. Next, we show \eqref{eq:dp-upper-heavy}:
	\begin{equation*}
		\begin{split}
			R_{\textnormal{DP}}( \theta(\cP^1_k(r)) , \epsilon_n )  &\overset{ \eqref{eq:stability-dp} }\leq \inf_{\beta_n} \left( R_{n,\infty}(\theta(\cP^1_k(r)), \beta_n) + \frac{2 \beta_n^2 }{ \epsilon_n^2 } \right) \\
			& \overset{ \textnormal{Theorem } \ref{th:worst-heavy-tail} \textnormal{ Part (a)} }= \inf_{\beta_n} \left( \left( \frac{r^{2k}}{(n\beta_n)^{2(k-1)}} \wedge r^2 + \frac{r^2}{n} \right) +   \frac{2 \beta_n^2 }{ \epsilon_n^2 } \right)\\
			& \lesssim  \frac{r^2}{n} + \frac{r^{2}}{(n \epsilon_n)^{2-2/k}}\wedge r^2. \quad \left( \textnormal{take } \beta_n = \left\{\begin{array}{ll}
				rn^{-1}(n \epsilon_n)^{1/k} & \epsilon_n \geq 1/n\\
				r\epsilon_n & \epsilon_n < 1/n
			\end{array} \right.\right)
		\end{split}
	\end{equation*}

Finally, we show \eqref{eq:stability-lower-heavy}. By \eqref{eq:privacy-bounds}, we know there exists some universal constant $c_0$ such that $R_{\DP}( \theta(\cP^1_k(r)) , \epsilon_n ) \geq c_0 \left( \frac{r^2}{n} + \left(\frac{r^{2}}{(n \epsilon_n)^{2-2/k}}\wedge r^2\right)  \right)$.  Then when $\epsilon_n \leq \frac{1}{n}$, we have $R_{\DP}( \theta(\cP^1_k(r)) , \epsilon_n ) \geq c_0 r^2$. So when $\beta_n \leq \frac{\sqrt{c_0}r}{2n}$, then by the lower bound in \eqref{eq:dp-stability},
\begin{equation*}
    R_{n,\infty}(\theta(\cP^1_k(r)), \beta_n) \geq R_{\textnormal{DP}}( \theta(\cP^1_k(r)) , 1/n ) - c_0r^2/2  \geq c_0 r^2/2 \gtrsim r^2. 
\end{equation*}
When $\beta_n > \frac{\sqrt{c_0}r}{2n}$, we take $\epsilon_n = \sqrt{\frac{2\beta_n^2}{R_{n,\infty}(\theta(\cP^1_k(r)), \beta_n)}}$, then by the lower bound in \eqref{eq:dp-stability}, 
\begin{equation*}
	\begin{split}
		R_{n,\infty}(\theta(\cP^1_k(r)), \beta_n) & \gtrsim  R_{\textnormal{DP}}\left( \theta(\cP^1_k(r)) , \sqrt{ \frac{2\beta_n^2}{R_{n,\infty}(\theta(\cP^1_k(r)), \beta_n)} } \right) \\
		& \overset{ \eqref{eq:privacy-bounds} } \gtrsim  \frac{r^2}{n} + \left(\frac{R_{n,\infty}(\theta(\cP^1_k(r)), \beta_n)}{2n^2 \beta_n^2} \right)^{1-1/k}r^2 \wedge r^2\\
		\Longrightarrow  & R_{n,\infty}(\theta(\cP^1_k(r)), \beta_n)  \gtrsim \left( \frac{r^{2k}}{(n\beta_n)^{2(k-1)}} \wedge r^2\right) + \frac{r^2}{n} 
	\end{split}
\end{equation*}
By combining two cases, we have $R_{n,\infty}(\theta(\cP^1_k(r)), \beta_n)  \gtrsim \left( \frac{r^{2k}}{(n\beta_n)^{2(k-1)}} \wedge r^2\right) + \frac{r^2}{n} $.

\section{Proofs for Section \ref{sec:lower-bound} } \label{app:proof-lower-bound}

\subsection{Proof of Theorem \ref{th:worst-lower-bound-multi}}
 First, $R_{n,\infty}(\theta(\cP), \beta_n) \geq  R_n(\theta(\cP))$ is trivial. So we just need to show that $R_{n,\infty}(\theta(\cP), \beta_n)$ is greater than the first term on the right-hand side of \eqref{ineq:risk-lower-bound}.
	For any estimator $\wh{\theta}$ that is $\beta_n$-worst-case-stable, let 
	\begin{equation} \label{eq:risk-upper-bound}
		\sigma_n^2 := \sup_{P \in \cP} \bbE_{\cD_n, \xi} \|  \wh{\theta}(\cD_n; \xi) - \theta(P) \|_2^2.
	\end{equation}
	For any $P_0, P_1 \in \cP$ and any coupling $(\cD_n, \cD_n')$ between $P_0^{\otimes n}$ and $P_1^{\otimes n}$, by the calculation in \eqref{ineq:upperbound-risk}, we have
	\begin{equation} \label{ineq:upper-bound}
		\begin{split}
			& \bbE\|\wh{\theta}(\cD_n) - \wh{\theta}(\cD_n')\|_2 \geq \|\theta(P_0) - \theta(P_1)\|_2 - 2\sigma_n.
		\end{split}
	\end{equation} By combining it with \eqref{ineq:stability-induced-upper-bound}, we get 
	\begin{equation} \label{ineq:lower-bound-sigma}
		\begin{split}
			& \beta_n \bbE  [ d_{\Ham}(\cD_n, \cD_n')] \geq \|\theta(P_0) - \theta(P_1)\|_2 - 2\sigma_n \\
			\overset{\textnormal{as } \sigma_n \geq 0}\Longleftrightarrow & \sigma_n \geq \frac{ (\|\theta(P_0) - \theta(P_1)\|_2 - \bbE  [ d_{\Ham}(\cD_n, \cD_n')] \beta_n)_{+} }{2}
		\end{split}
	\end{equation}
	Since the lower bound in \eqref{ineq:lower-bound-sigma} holds for any estimator, any $P_0, P_1 \in \cP$, and any coupling $(\cD_n, \cD_n')$ between $P_0^{\otimes n}$ and $P_1^{\otimes n}$, we have shown that $R_{n,\infty}(\theta(\cP), \beta_n)$ is greater than the first term on the right-hand-side of \eqref{ineq:risk-lower-bound}.

\subsection{Proof of Theorem \ref{th:average-lower-bound-multi} }
First, we note again we just need to show the first term on the right-hand side of \eqref{ineq:ave-risk-lower-bound}. In addition, without loss of generality, we are going to assume $\wh{\theta}(\cD_n, \xi)$ is symmetric with respect to $\cD_n$, i.e., for any permutation $\sigma$ on $[n]$, $\wh{\theta}(\cD_n, \xi) = \wh{\theta}(\sigma(\cD_n), \xi)$. This is because for any estimator, we can always symmetrize it so that the symmetrized estimator has no worse stability guarantee and has risk no bigger than the original one. In particular, given any $\wh{\theta}(\cD_n, \xi)$, define $\tilde{\theta}(\cD_n, \xi) = \frac{1}{|\mathcal{S}_n|} \sum_{\sigma \in \mathcal{S}_n} \wh{\theta}(\sigma(\cD_n), \xi)$, where $\mathcal{S}_n$ denotes the set of all permutations supported on $[n]$, we always have
\begin{equation} \label{ineq:symmetrization}
	\begin{split}
		&\frac{1}{(n+1)^2} \sup_{\cD_{n+1}} \sum_{1 \leq i, j \leq n+1} \bbE_{\xi} \|\wt{\theta}(\cD_{n+1}^{\setminus i}; \xi )  - \wt{\theta}(\cD_{n+1}^{\setminus j}; \xi )\|^{p}_2 \\
		& = \frac{1}{(n+1)^2} \sup_{\cD_{n+1}} \sum_{1 \leq i, j \leq n+1} \bbE_{\xi} \left\| \frac{1}{|\mathcal{S}_n|} \sum_{\sigma \in \mathcal{S}_n} \left( \wh{\theta}(\sigma(\cD_{n+1}^{\setminus i}); \xi )  - \wh{\theta}(\sigma(\cD_{n+1}^{\setminus j}); \xi ) \right) \right\|_2^{p} \\
		& \overset{(a)}\leq \frac{1}{(n+1)^2} \sup_{\cD_{n+1}} \sum_{1 \leq i, j \leq n+1} \frac{1}{|\mathcal{S}_n|} \sum_{\sigma \in \mathcal{S}_n} \bbE_{\xi} \left\| \wh{\theta}(\sigma(\cD_{n+1}^{\setminus i}); \xi )  - \wh{\theta}(\sigma(\cD_{n+1}^{\setminus j}); \xi ) \right\|_2^{p}\\
		& \leq \frac{1}{|\mathcal{S}_n|} \sum_{\sigma \in \mathcal{S}_n}  \frac{1}{(n+1)^2} \sup_{\cD_{n+1}}\sum_{1 \leq i, j \leq n+1} \bbE_{\xi} \left\| \wh{\theta}(\sigma(\cD_{n+1}^{\setminus i}); \xi )  - \wh{\theta}(\sigma(\cD_{n+1}^{\setminus j}); \xi ) \right\|_2^{p}  \leq \beta_n
	\end{split}
\end{equation} where (a) is by Jensen's inequality as $\|\cdot\|_2^{p}$ is a convex function when $p \geq 1$ and
\begin{equation*}
	\begin{split}
		\bbE_{\cD_n, \xi} \|  \tilde{\theta}(\cD_n; \xi) - \theta(P) \|_2^2 \leq \frac{1}{|\mathcal{S}_n|^2} |\mathcal{S}_n| \sum_{\sigma \in \mathcal{S}_n} \bbE_{\cD_n, \xi} \|  \wh{\theta}(\sigma(\cD_n); \xi) - \theta(P) \|_2^2 = \bbE_{\cD_n, \xi} \|  \wh{\theta}(\cD_n; \xi) - \theta(P) \|_2^2.
	\end{split}
\end{equation*}

We divide the rest of the proof into two parts: Part I is for estimating a functional that is linear in $P$ and Part II is for the general setting.

\vskip.2cm
{\bf (Part I: Linear functional setting)}
We will prove the lower bound based on a two-point mixture distribution construction. Take any $P_0, P_1 \in \cP$ such that $(1-t) P_0 + t P_1 \in \cP$ for all $t \in [0,1]$. Given any $\eta \in [0,1/2]$, let $P = \eta P_0 + (1 - \eta) P_1 $ and $P' = \eta P_1 + (1 - \eta) P_0 $. Then $\theta(P) = \eta \theta(P_0) + (1 - \eta)  \theta(P_1)$ and $\theta(P') = \eta \theta(P_1) + (1 - \eta) \theta(P_0)$. For any symmetric estimator that is $\beta_n$-average-case-stable, let 
	\begin{equation*} 
		\sigma_n^2 := \sup_{P \in \cP} \bbE_{\cD_n \overset{i.i.d.}\sim P, \xi} \|  \wh{\theta}(\cD_n;\xi) - \theta(P) \|_2^2.
	\end{equation*} 
	Let us generate $\cD_n = (X_1, \ldots, X_n) \sim P^{\otimes n}$ and $\cD_n' = (X_1', \ldots, X_n') \sim P^{'\otimes n}$ in a coupled way as follows: 
	\begin{itemize}
		\item first, generate $B_i \overset{i.i.d.}\sim \textnormal{Bernoulli}(2\eta)$ for $i = 1, \ldots, n$;
		\item for each $i = 1, \ldots, n$, if $B_i$ = 1, then independently generate $A_i \sim \textnormal{Bernoulli}(1/2)$, if $A_i =1$, generate $X_i = X_i' \sim P_0$ and if $A_i = 0$, generate $X_i = X_i' \sim P_1$; if $B_i = 0$, then set $A_i =-1$ generate $X_i \sim P_1$ and independently generate $X_i' \sim P_0$.
	\end{itemize}
	
	First, we have
	\begin{equation} \label{ineq:lower-bound}
		\begin{split}
			\|\bbE[ \wh{\theta}(\cD_n;\xi) - \wh{\theta}(\cD'_n; \xi) ] \|_2 &= \|\bbE[ \wh{\theta}(\cD_n;\xi) - \theta(P) + \theta(P) - \theta(P') + \theta(P') - \wh{\theta}(\cD'_n; \xi) ] \|_2 \\
			& \geq \|\theta(P) - \theta(P') \|_2 - 2 \sigma_n = (1 - 2 \eta) \|\theta(P_0)  - \theta(P_1)\|_2 - 2\sigma_n.
		\end{split}
	\end{equation}
	
	Second, given $T_1 = \sum_{i=1}^n B_i$ and $T_2 = \sum_{i=1}^n B_i A_i$, by construction we have that $\cD_n$ and $\cD_n'$ differ on $(n-T_1)$ data points and among the $T_1$ data points $\cD_n$ and $\cD_n'$ share, $T_2$ of them come from $P_0$ and $T_1 - T_2$ of them come from $P_1$. Since $\wh{\theta}$ is symmetric in its input, without loss of generality, we can assume that given $T_1, T_2$, $\cD_n$ and $\cD'_n$ can be written as
	\begin{equation*}
		\begin{split}
			\cD_n  &= (\tilde{X}_1, \ldots, \tilde{X}_{T_2}, X_{T_2+1}, \ldots, X_{T_1}, X_{T_1+1}, \ldots, X_n), \\
			\cD'_n  &= (\tilde{X}_1, \ldots, \tilde{X}_{T_2}, X_{T_2+1}, \ldots, X_{T_1}, \tilde{X}_{T_1+1}, \ldots, \tilde{X}_n),
		\end{split}
	\end{equation*} where $(\tilde{X}_1, \ldots, \tilde{X}_{T_2}, \tilde{X}_{T_1+1}, \ldots, \tilde{X}_n) \overset{i.i.d.}\sim P_0$ and $(X_{T_2+1}, \ldots, X_n) \overset{i.i.d.}\sim P_1$. So we can construct a sequence of datasets $\cD^{(0)}, \cdots, \cD^{(n-T_1)}$ with $\cD^{(0)} = \cD_n$, $\cD^{(n-T_1)} = \cD_n'$, and 
	\begin{equation*}
		\begin{split}
			\cD^{(i)} = (\tilde{X}_1, \ldots, \tilde{X}_{T_2}, X_{T_2 + 1} \ldots, X_{T_1}, \tilde{X}_{T_1 + 1} , \ldots, \tilde{X}_{T_1+i}, X_{T_1+i+1}, \ldots, X_n)
		\end{split}
	\end{equation*} such that
	\begin{equation} \label{eq:chain-decompo}
		 \wh{\theta}(\cD_n;\xi) - \wh{\theta}(\cD'_n; \xi) = \sum_{i=1}^{n-T_1} (\wh{\theta}(\cD^{(i-1)};\xi) - \wh{\theta}(\cD^{(i)};\xi)).
	\end{equation} Now let us consider upper bound $\|\bbE[ \wh{\theta}(\cD^{(i-1)};\xi) - \wh{\theta}(\cD^{(i)};\xi) | T_1, T_2 ] \|_2^{p}$ for each $i = 1, \ldots, n - T_1$. Let us denote $$\cD_{n+1} = \cD^{(i-1)} \bigcup \cD^{(i)} = (\tilde{X}_1, \ldots, \tilde{X}_{T_2}, X_{T_2 + 1} \ldots, X_{T_1}, \tilde{X}_{T_1 + 1} , \ldots, \tilde{X}_{T_1+i}, X_{T_1 + i}, X_{T_1+i+1}, \ldots, X_n).$$ Notice that $\cD_{n+1}$ above also depends on $i$, but we suppress the notion for simplicity. By construction, we have
	\begin{equation} \label{eq:special-leave-one-out}
		\begin{split}
			\cD^{(i-1)} = \cD_{n+1}^{\setminus T_1+i} \quad \textnormal{ and }\quad \cD^{(i)} = \cD_{n+1}^{\setminus T_1+i +1}.
		\end{split}
	\end{equation} In addition 
	\begin{equation} \label{eq:expectation-equal}
		\begin{split}
			\bbE[\wh{\theta}(\cD_{n+1}^{\setminus T_1+i}; \xi)| T_1, T_2 ] &=  \bbE[\wh{\theta}(\cD_{n+1}^{\setminus j};\xi)| T_1, T_2 ], \quad \forall j = 1, \ldots, T_2, T_1 +1, \ldots, T_1 + i,\\
			\bbE[\wh{\theta}(\cD_{n+1}^{\setminus T_1+1+i};\xi)| T_1, T_2 ] &=  \bbE[\wh{\theta}(\cD_{n+1}^{\setminus j};\xi)| T_1, T_2 ], \quad \forall j = T_2+1, \ldots,  T_1 , T_1 + i + 1, \ldots, n+1.
		\end{split}
	\end{equation} Thus,
	\begin{equation*}
		\begin{split}
			 & (T_2 + i) (T_1 - T_2 + n + 1 - T_1 - i) \|\bbE[\wh{\theta}(\cD^{(i-1)};\xi) - \wh{\theta}(\cD^{(i)};\xi)| T_1, T_2 ]\|_2^{p} \\
			= & (T_2 + i) (n - T_2 + 1 - i) \|\bbE[\wh{\theta}(\cD^{(i-1)};\xi) - \wh{\theta}(\cD^{(i)};\xi)| T_1, T_2 ]\|_2^{p} \\
			\overset{ \eqref{eq:special-leave-one-out}, \eqref{eq:expectation-equal} }= & \left\| \bbE_{\cD_n,\xi} \left[ \sum_{j \in [T_2] \bigcup ([T_1 + i] \setminus [T_1]) } \sum_{k \in ([T_1]\setminus [T_2]) \bigcup ([n+1] \setminus [T_1 + i])}  \left(  \wh{\theta}(\cD_{n+1}^{\setminus j};\xi) - \wh{\theta}(\cD_{n+1}^{\setminus k})\right) | T_1, T_2 \right] \right\|_2^{p} \\
			\leq & \frac{1}{2} \bbE_{\cD_n} \left[ \bbE_{\xi} \left[ \sum_{j,k \in [n+1], j \neq k } \left\|  \wh{\theta}(\cD_{n+1}^{\setminus j};\xi) - \wh{\theta}(\cD_{n+1}^{\setminus k})\right\|_2^{p}\right] \Big| T_1, T_2  \right]\\
			\overset{(a)}\leq & (n+1)^2 \beta^{p}_n/2,
		\end{split}
	\end{equation*} where (a) is by the assumption on average-case stability. Thus, we obtain
	\begin{equation} \label{ineq:individual-term-bound}
		 \|\bbE[\wh{\theta}(\cD^{(i-1)};\xi) - \wh{\theta}(\cD^{(i)};\xi)| T_1, T_2 ]\|_2^{p} \leq \frac{(n+1)^2 \beta^{p}_n}{2(T_2 + i) (n - T_2 + 1 - i)}.
	\end{equation}  As a result
	\begin{equation} \label{ineq:bound-condi-T}
		\begin{split}
			&\|\bbE[ \wh{\theta}(\cD_n;\xi) - \wh{\theta}(\cD'_n; \xi) | T_1, T_2 ] \|_2^{p} \\
			& \overset{\eqref{eq:chain-decompo}}\leq (n-T_1)^{p-1} \sum_{i=1}^{n-T_1} \|\bbE (\wh{\theta}(\cD^{(i-1)};\xi) - \wh{\theta}(\cD^{(i)};\xi)| T_1, T_2 )\|_2^{p} \\
			& \overset{ \eqref{ineq:individual-term-bound} }\leq (n-T_1)^{p-1}\sum_{i=1}^{n-T_1} \frac{(n+1)^2 \beta^{p}_n}{2(T_2 + i) (n - T_2 + 1 - i)} \\
			& =(n-T_1)^{p-1}\frac{(n+1)\beta^{p}_n}{2} \sum_{i=1}^{n-T_1} \left( \frac{1}{T_2 + i} + \frac{1}{n - T_2 + 1 - i} \right) \\
			& \overset{(a)}\leq  (n-T_1)^{p-1}\frac{(n+1)\beta^{p}_n}{2}  \left( \log \left( \frac{n + T_2 - T_1}{T_2 +1} \right) + \frac{1}{T_2 +1}   +  \log \left( \frac{n - T_2}{T_1 - T_2 +1} \right) + \frac{1}{T_1 - T_2 +1} \right),
		\end{split}
	\end{equation} where (a) is because 
	\begin{equation*}
		\begin{split}
			\sum_{i=1}^{n-T_1} \frac{1}{T_2 + i} \leq \int_1^{n - T_1} \frac{1}{T_2 + x} dx + \frac{1}{T_2 + 1} = \log \left( \frac{n + T_2 - T_1}{T_2 +1} \right) + \frac{1}{T_2 +1}
		\end{split}
	\end{equation*} and similarly, we also have $\sum_{i=1}^{n-T_1} \frac{1}{n - T_2 + 1 - i} \leq   \log \left( \frac{n - T_2}{T_1 - T_2 +1} \right) + \frac{1}{T_1 - T_2 +1}$. We note that the upper bound in \eqref{ineq:bound-condi-T} holds when $T_1 \leq n-1$, while $T_1 = n$, the difference is zero and the rest of the proof still follows with minor modification.

	By construction, we have $T_1 \sim \textnormal{Binomial}(n, 2\eta)$ and $T_2 | T_1 \sim \textnormal{Binomial}(T_1, 1/2)$. Since, $T_1 - T_2$ and $T_2$ follows the same distribution given $T_1$, we have 
	\begin{equation} \label{eq:equal-expectation}
		\bbE\left(  \log \left( \frac{n + T_2 - T_1}{T_2 +1} \right) + \frac{1}{T_2 +1}  | T_1 \right) = \bbE \left(  \log \left( \frac{n - T_2}{T_1 - T_2 +1} \right) + \frac{1}{T_1 - T_2 +1} | T_1 \right).
	\end{equation}Thus,
	\begin{equation} \label{ineq:exp-diff-upper}
		\begin{split}
			\|\bbE[ \wh{\theta}(\cD_n;\xi) - \wh{\theta}(\cD'_n; \xi) ] \|_2^{p} &= \left\|\bbE\left[ \bbE\left[  \wh{\theta}(\cD_n;\xi) - \wh{\theta}(\cD'_n; \xi) | T_1 , T_2\right]  \right] \right\|_2^{p} \\
			& \overset{ \eqref{ineq:bound-condi-T}, \eqref{eq:equal-expectation} }\leq (n+1) \beta^{p}_n \bbE  \left(  (n-T_1)^{p-1} \left(\log \left( \frac{n - T_1 + T_2}{T_2 +1} \right) + \frac{1}{T_2 +1} \right) \right).
		\end{split}
	\end{equation} The following lemma provides an upper bound for $\bbE  \left(   \log \left( \frac{n - T_1 + T_2}{T_2 +1} \right) + \frac{1}{T_2 +1} \right)$ and its conditional version. The proof is provided in the subsequent subsections.
\begin{Lemma} \label{lm:binomial-expectation}
	Suppose $T_1 \sim \textnormal{Binomial}(n, q)$ and given $T_1$, $T_2 \sim \textnormal{Binomial}(T_1, 1/2)$, then 
	\begin{equation} \label{eq:binomial-exp-1}
		\bbE  \left(  \log \left( \frac{n - T_1 + T_2}{T_2 +1} \right) + \frac{1}{T_2 +1} \Big| T_1 \right) \leq   \log \left( \frac{2n}{T_1 + 1} -1 \right) + \frac{2}{T_1 + 1}. 
	\end{equation} Moreover,
	\begin{equation}\label{eq:binomial-exp-2} 
		\bbE  \left(  \log \left( \frac{n - T_1 + T_2}{T_2 +1} \right) + \frac{1}{T_2 +1} \right) \leq \log\left( \frac{2n}{q(n+1)} - 1 \right) + \frac{2}{q (n+1)} .
	\end{equation}
\end{Lemma} Then, we can further have
\begin{equation} \label{ineq:upper-bound-diff}
	\begin{split}
		&\|\bbE[ \wh{\theta}(\cD_n;\xi) - \wh{\theta}(\cD'_n; \xi) ] \|_2^{p} \\
		\overset{ \eqref{ineq:exp-diff-upper},\eqref{eq:binomial-exp-1} }\leq & (n+1) \beta^{p}_n \bbE  \left(  (n-T_1)^{p-1} \left(\log \left( \frac{2n}{T_1 +1} -1 \right) + \frac{2}{T_1 +1} \right) \right) \\
		\leq & (n+1) \beta^{p}_n \bbE  \left(  (n-T_1)^{p-1}\left(\left( \frac{2n}{T_1 +1} -2 \right) + \frac{2}{T_1 +1} \right) \right) \quad(\textnormal{as } \log(1+x) \leq x, \forall x > -1)\\
		= &2(n+1) \beta^{p}_n \bbE \left( \frac{(n-T_1)^{p}}{T_1 + 1} \right).
	\end{split}
\end{equation}
The next lemma provides an upper bound for $\bbE \left( \frac{(n-T_1)^{p}}{T_1 + 1} \right)$ and the proof of this lemma is provided in the subsequent subsections.

\begin{Lemma} \label{lm:bino-exp-lm2}
	Suppose $T_1 \sim \textnormal{Binomial}(n, q)$, then for any $p > 0$,
	\begin{equation*}
		\bbE \left( \frac{(n-T_1)^{p}}{T_1 + 1} \right) \leq \frac{((n+1)(1-q))^{p} \left( 1 + \frac{p}{2 (n+1)(1-q)} \right)^{p}}{(n+1)q} \wedge n^{p}.
	\end{equation*}
\end{Lemma}

By replacing $q$ by $2\eta$ in Lemma \ref{lm:bino-exp-lm2}, we obtain from \eqref{ineq:upper-bound-diff} that 
\begin{equation}\label{ineq:upper-bound-diff2}
	\begin{split}
		\|\bbE[ \wh{\theta}(\cD_n;\xi) - \wh{\theta}(\cD'_n; \xi) ] \|_2^{p} &\leq 2 \beta^{p}_n \left( \frac{((n+1)(1-2\eta))^{p} \left( 1 + \frac{p}{2 (n+1)(1-2\eta)} \right)^{p}}{2\eta} \wedge n^{p}(n+1) \right)\\
		& \leq 2 \beta^{p}_n \left( \frac{((n+1)(1-2\eta))^{p} \left( 1 + \frac{p}{2 (n+1)(1-2\eta)} \right)^{p}}{2\eta} \wedge (n+1)^{p+1} \right)
	\end{split}
\end{equation} Combining \eqref{ineq:upper-bound-diff2} and \eqref{ineq:lower-bound}, we get 
\begin{equation*}
	\begin{split}
		2^{1/p} \beta_n (n+1)  \left( \frac{(1-2\eta) \left(1 + \frac{p}{2(n+1)(1 - 2 \eta)} \right)}{(2\eta)^{1/p}} \wedge (n+1)^{1/p} \right) \geq (1 - 2 \eta) \|\theta(P_0)  - \theta(P_1)\|_2 - 2\sigma_n.
	\end{split}
\end{equation*} Thus, we get 
\begin{equation*}
	\sigma_n \geq \left( \frac{ (1 - 2 \eta) \|\theta(P_0)  - \theta(P_1)\|_2 - 2^{1/p} \beta_n (n+1)  \left( \frac{(1-2\eta) \left(1 + \frac{p}{2(n+1)(1 - 2 \eta)} \right)}{(2\eta)^{1/p}} \wedge (n+1)^{1/p} \right) }{2} \right)_+.
\end{equation*} Since this holds for any $\eta \in [0,1/2]$ and any estimator satisfies the stability constraint, we have finished the proof.

\vskip.2cm
{\bf (Part II: General setting)}  In this setting, we can set $\eta = 0$ in the proof of Part I, and the proof still goes through in the general setting. When $\eta = 0$, $T_1 = T_2 = 0$. From \eqref{ineq:lower-bound}, we know
	\begin{equation}
		\begin{split}
			\|\bbE[ \wh{\theta}(\cD_n;\xi) - \wh{\theta}(\cD'_n; \xi) ] \|_2 \geq  \|\theta(P_0)  - \theta(P_1)\|_2 - 2\sigma_n.
		\end{split}
	\end{equation} At the same time, following the same analysis as in \eqref{ineq:exp-diff-upper} with $T_1 = T_2 = 0$, we have 
	\begin{equation*}
		\begin{split}
			\|\bbE[ \wh{\theta}(\cD_n;\xi) - \wh{\theta}(\cD'_n; \xi) ] \|^{p}_2 \leq  (n+1)^{p} \beta_n^{p} \left( \log n + 1 \right). 
		\end{split}
	\end{equation*} A combination of these two yields
	\begin{equation*}
		\begin{split}
			\sigma_n \geq \left( \|\theta(P_0)  - \theta(P_1)\|_2 - (n+1) \beta_n (\log n + 1)^{1/p}\right)_+/2. 
		\end{split}
	\end{equation*}  Since this holds for any $P_0, P_1 \in \cP$ and any estimator satisfies the stability constraint, we have finished the proof.

\subsubsection{Proof of Lemma \ref{lm:binomial-expectation}}
First, it is easy to check that if $X \sim \textnormal{Binomial}(n, q)$, then
\begin{equation} \label{eq:binomial-inverse-expectation}
	\begin{split}
		\bbE \left( \frac{1}{X + 1} \right) = \sum_{k=0}^{n} \frac{1}{1 + k} {n \choose k} q^k (1-q)^{n-k} =  \sum_{k=0}^{n} \frac{1}{(1 + n) q} {n +1 \choose k +1} q^{k+1} (1-q)^{n-k} = \frac{1 - (1-q)^{n + 1} }{(n + 1)q}.
	\end{split}
\end{equation}
Since $T_2 \sim \textnormal{Binomial}(T_1, 1/2)$, then $\bbE \left( \frac{1}{T_2 + 1} | T_1 \geq 1 \right) = \frac{1 - (1/2)^{T_1 + 1} }{(T_1 + 1)/2}$. Notice that even when $T_1 = 0$, we still have $\bbE \left( \frac{1}{T_2 + 1} | T_1 = 0 \right) = \frac{1 - (1/2)^{T_1 + 1} }{(T_1 + 1)/2}$. Thus, $\bbE \left( \frac{1}{T_2 + 1}|T_1\right) = \frac{1 - (1/2)^{T_1 + 1} }{(T_1 + 1)/2} \leq \frac{2}{T_1 + 1} $. Next,
\begin{equation*}
	\begin{split}
		&\bbE \left(   \log \left( \frac{n - T_1 + T_2}{T_2 +1} \right) | T_1 \leq n-1 \right) \\
		&= \bbE \left(   \log \left( 1+\frac{n - T_1 -1}{T_2 +1} \right) | T_1 \leq n-1 \right) \\
		 & \leq \log \left( 1 + (n- T_1 - 1) \bbE \left( \frac{1}{T_2 + 1} | T_1 \leq n-1 \right) \right) \quad( \text{Jensen's inequality)} \\
		 & \overset{\eqref{eq:binomial-inverse-expectation}} \leq \log \left( 1 +   (n- T_1 - 1) \cdot \frac{2}{T_1 + 1}  \right) = \log \left( \frac{2n}{T_1 + 1} -1 \right).
	\end{split}
\end{equation*} In addition, when $T_1 = n$,
\begin{equation*}
	\bbE \left(   \log \left( \frac{n - T_1 + T_2}{T_2 +1} \right) | T_1 = n \right) = \bbE \left(   \log \left( \frac{ T_2}{T_2 +1} \right)\right) = -\infty \leq \log \left( \frac{2n}{n+ 1} -1 \right) = \log \left( \frac{2n}{T_1 + 1} -1 \right).
\end{equation*} Thus, we have shown $\bbE \left(   \log \left( \frac{n - T_1 + T_2}{T_2 +1} \right)| T_1 \right) \leq \log \left( \frac{2n}{T_1 + 1} -1 \right)$ holds almost surely. In summary,
\begin{equation*}
	\begin{split}
		\bbE  \left(  \log \left( \frac{n - T_1 + T_2}{T_2 +1} \right) + \frac{1}{T_2 +1} \right)& \leq \bbE \left( \log \left( \frac{2n}{T_1 + 1} -1 \right) + \frac{2}{T_1 + 1}  \right) \\
		& \leq   \log \left( 2n \bbE \left(\frac{1}{T_1 + 1}\right) -1 \right) + 2  \bbE \left(\frac{1}{T_1 + 1} \right) \quad( \text{Jensen's inequality)} \\
		& \overset{ \eqref{eq:binomial-inverse-expectation} }= \log\left( 2n \cdot \frac{1 - (1 - q)^{n+1} }{q (n+1)} - 1 \right) + 2 \frac{1 - (1 - q)^{n+1}}{q(n+1)} \\
		& \leq \log\left( \frac{2n}{q(n+1)} - 1 \right) + \frac{2}{q (n+1)}.
	\end{split}
\end{equation*}

\subsubsection{Proof of Lemma \ref{lm:bino-exp-lm2}} 
First, we note that the $n^{p}$ upper bound is straightforward. Now we prove the other part of the upper bound. Note that 
\begin{equation*}
	\begin{split}
		\bbE \left( \frac{(n-T_1)^{p}}{T_1 + 1}  \right) &= \sum_{k=0}^{n} \frac{(n-k)^{p}}{k+1}  {n \choose k} q^k (1-q)^{n-k} \\
		& = \frac{1}{(n+1)q} \sum_{k=1}^{n+1} \frac{(n+1 - k)^{p} (n+1)!}{k!(n+1-k)!} q^k (1 - q)^{n+1-k}\\
		& \leq \frac{1}{(n+1)q} \bbE_{Y \sim \textnormal{Binom}(n+1,1-q) }[Y^{p}] \\
		& \overset{(a)}\leq   \frac{((n+1)(1-q))^{p} \left( 1 + \frac{p}{2 (n+1)(1-q)} \right)^{p}}{(n+1)q}.
	\end{split}
\end{equation*} where (a) is by the bound on the raw moment of the Binomial provided in the following lemma from Corollary 1 of \cite{ahle2022sharp}.
\begin{Lemma}\label{lm:binomial-moment}
	Suppose $X \sim \textnormal{Binom}(n,p)$. Then $\bbE[X^k] \leq (np)^k \left(1 + \frac{k}{2np} \right)^k$ for any $k \geq 1$.
\end{Lemma} This finishes the proof.

\subsection{Proof of Corollary \ref{coro:average-case-stability} } The result follows if we plug in $p=1$ and $\eta = 1/4$ in Theorem \ref{th:average-lower-bound-multi} \eqref{ineq:average-lower-linear-fun}.

\section{Proofs for Section \ref{sec:bounded-mean}}

\subsection{Proof of Theorem \ref{th:worst-bounded} }
Notice that for this class, $\|X_i\|_2 \leq r$, thus the sample mean is $\frac{2r}{n}$-$\ell_\infty$-stable. By monotonicity, it is also $\frac{2r}{n}$-$\ell_p$-stable for all $p \geq 1$. So when $\beta_n \geq \frac{2r}{n}$, the sample mean can achieve the upper bound, and the lower bound is the unconstrained minimax risk and it can be proved using a standard Le Cam's two-point argument \citep{yu1997assouad}. For simplicity, we omit it here. When $\beta_n \leq \frac{r}{10n} $, 
\begin{equation*}
	R_{n,p}(\theta(\cP), \beta_n) \geq R_{n,1}(\theta(\cP), \beta_n) \gtrsim r^2,
\end{equation*}
where the last inequality is by Corollary \ref{coro:average-case-stability} where we choose $P_0$ and $P_1$ such that $\theta(P_0) = r e_1$ and $\theta(P_1) = -re_1$, and $e_1$ is a unit vector with first coordinate to be $1$ and others are zero. The upper bound is also straightforward as $\|\wb{X}\|_2 \leq r$. 

\subsection{Proof of Theorem \ref{th:bounded-worst-refined} }
We first show that $R_{n,p}\left(\theta(\cP^1_{\textnormal{bound}}(r)), \beta_n\right) \asymp \frac{r^2}{n}$ when $\beta_n \geq \frac{2^{1-1/p}r}{n}$. In particular, let us compute the $\ell_p$ stability guarantee for the sample mean:
\begin{equation} \label{eq:lp-stab-maximum}
	\begin{split}
		&\frac{1}{(n+1)^2} \sup_{\cD_{n+1} \in \cX^{n+1}} \sum_{1 \leq i, j \leq n+1} \bbE_{\xi} |\wh{\theta}(\cD_{n+1}^{\setminus i}; \xi )  - \wh{\theta}(\cD_{n+1}^{\setminus j}; \xi )|^{p} \\
		= & \frac{1}{(n+1)^2} \sup_{\cD_{n+1} \in \cX^{n+1}} \sum_{1 \leq i, j \leq n+1} \frac{|X_i - X_j|^{p}}{n^{p}} = \frac{2}{(n+1)^2 n^{p}} \sup_{\cD_{n+1} \in \cX^{n+1}} \sum_{1 \leq i< j \leq n+1} |X_i - X_j|^{p}  .
	\end{split}
\end{equation} Given any $(X_1, \ldots, X_{n+1})$ with each of them lying in $[-r,r]$, let $f(X_1, \ldots, X_n) := \sum_{1 \leq i < j \leq n+1} |X_i - X_j|^{p}$. The first claim is that the maximum of $f(X_1, \ldots, X_{n+1})$ over different $\cD_{n+1} \in \cX^{n+1}$ is obtained when $X_i \in \{-r,r\}$ for all $i \in [n+1]$. To see that, let us consider fixing $(X_2, \ldots, X_{n+1})$ and optimizing $f(X_1, \ldots, X_n)$ over $X_1$, i.e., we want to solve $\max_{X_1 \in [-r,r]} g(X_1) = \sum_{2\leq j \leq n+1} |X_j - X_1|^{p} $. Notice that $g(X_1)$ is convex in $X_1$ as $p \geq 1$ and it is optimized over a compact set, then a standard result in convex optimization theory says that the maximum of a convex function on a compact set is achieved at the boundary point \citep[Theorem 32.2]{rockafellar1997convex}, i.e., $X_1 \in \{r,-r\}$ here. By a similar argument for each coordinate, we get that the maximum of $f(X_1, \ldots, X_{n+1})$ is obtained when $X_i \in \{-r,r\}$ for all $i \in [n+1]$. Let $m$ denote the number of $r$'s in $\cD_{n+1}$, then 
\begin{equation*}
	\begin{split}
		\sup_{\cD_{n+1} \in \cX^{n+1}} f(X_1, \ldots, X_n) = \sup_{m \in [n+1]} m (n+1-m) (2r)^{p}\leq \frac{(n+1)^2}{4} (2r)^{p},
	\end{split}
\end{equation*} and the last inequality can be achieved if $n$ is odd and $m = \frac{n+1}{2}$. By plugging it into \eqref{eq:lp-stab-maximum}, we get that for the sample mean
\begin{equation*}
	\frac{1}{(n+1)^2} \sup_{\cD_{n+1} \in \cX^{n+1}} \sum_{1 \leq i, j \leq n+1} \bbE_{\xi} |\wh{\theta}(\cD_{n+1}^{\setminus i}; \xi )  - \wh{\theta}(\cD_{n+1}^{\setminus j}; \xi )|^{p}  \leq \frac{1}{2} \left( \frac{2r}{n} \right)^{p}.
\end{equation*} Thus, when $\beta_n \geq \frac{2^{1-1/p}r}{n}$, sample mean is $\beta_n$-$\ell_p$-stable and it can achieve the $\frac{r^2}{n}$ error bound. The lower bound in this regime is the minimax lower bound without any stability constraint.

Next, we show $R_{n,p}\left(\theta(\cP^1_{\textnormal{bound}}(r)), \beta_n\right) \asymp r^2$ when $\beta_n \leq \frac{2^{1-1/p}r}{n(1 + \delta)}$ for any small constant $\delta > 0$ given $n$ is greater than a sufficiently large constant depending on $\delta$ only. First, we note that in terms of the upper bound, it is trivial as $|\wb{X}| \leq r$ and $|\theta| \leq r$. So we just need to focus on the lower bound. We note that the result in Corollary \ref{coro:average-case-stability} is not sharp enough to yield such a result and we will use \eqref{ineq:average-lower-linear-fun}. In addition, we note that $\theta(\cP^1_{\textnormal{bound}}(r))$ is a convex class, i.e., given any $P_0, P_1 \in \theta(\cP^1_{\textnormal{bound}}(r))$, their convex combination lies in $\theta(\cP^1_{\textnormal{bound}}(r))$ as well. We divide the proof into two cases based on the value of $p$.
\begin{itemize}[leftmargin=*]
	\item Let us first consider $p \geq C_1 \log n$, where $C_1 > 0$ is some constant we will choose later. Then by \eqref{ineq:average-lower-linear-fun} with $\eta = 0$, we get
\begin{equation}
\begin{split}
	R_{n,p}\left(\theta(\cP^1_{\textnormal{bound}}(r)), \beta_n\right) &\geq \sup_{P_0, P_1 \in \cP^1_{\textnormal{bound}}(r)} \left( \frac{|\theta(P_0)  - \theta(P_1) | - 2^{1/p} (n+1)^{1+1/p} \beta_n }{2} \right)_+^2  \\
	& = \left( \frac{2r- 2^{1/p} (n+1)^{1+1/p} \beta_n }{2} \right)_+^2 \\
	& \geq  \left( \frac{2r- 2r \frac{(n+1)^{1+1/p}}{n(1 + \delta)} }{2} \right)_+^2 \quad(\textnormal{as } \beta_n \leq \frac{2^{1-1/p}r}{n(1 + \delta)})\\
	& = r^2 \left(1 - \frac{n+1}{n} \frac{(n+1)^{1/p}}{1 + \delta} \right)_+^2 \geq r^2 \left(1 - \frac{n+1}{n} \frac{e^{2/C_1}}{1 + \delta} \right)_+^2. \quad(\textnormal{as } p \geq C_1 \log n)
\end{split}
\end{equation} So as long as $n \geq C$ and $C_1 \geq C'$ for some large constants $C,C' > 0$ depending on $\delta$ only, then we have $R_{n,p}\left(\theta(\cP^1_{\textnormal{bound}}(r)), \beta_n\right) \geq c r^2$ for some $c > 0$ depending on $\delta$ only when $\beta_n \leq \frac{2^{1-1/p}r}{n(1 + \delta)}$.
\item Next, let us consider $1\leq p < C_1 \log n$. Again, by using \eqref{ineq:average-lower-linear-fun} with the other term in the minimum, we get
\begin{equation} \label{ineq:lp-stability-lower-2}
	\begin{split}
		&R_{n,p}\left(\theta(\cP^1_{\textnormal{bound}}(r)), \beta_n\right) \\
		&\geq \sup_{P_0, P_1 \in \cP} \sup_{\eta \in [0,1/2]}\left( \frac{ (1- 2 \eta) |\theta(P_0)  - \theta(P_1) |  - 2^{1/p} (n+1) \beta_n \left( \frac{1-2\eta + \frac{p}{2(n+1)} }{(2\eta)^{1/p}} \right) }{2} \right)^2_{+} \\
		& = \sup_{\eta \in [0,1/2]}\left( \frac{ 2r(1- 2 \eta)  - 2^{1/p} (n+1) \beta_n \left( \frac{1-2\eta + \frac{p}{2(n+1)} }{(2\eta)^{1/p}} \right) }{2} \right)^2_{+} \\
		& \geq \sup_{\eta \in [0,1/2]}\left( \frac{ 2r(1- 2 \eta)  - 2r \frac{n+1}{n(1 + \delta)} \left( \frac{1-2\eta + \frac{p}{2(n+1)} }{(2\eta)^{1/p}} \right) }{2} \right)^2_{+} \quad(\textnormal{as } \beta_n \leq \frac{2^{1-1/p}r}{n(1 + \delta)}) \\
		& = r^2 \sup_{\eta \in [0,1/2]} \left( (1 - 2 \eta) - \frac{n+1}{n(1 + \delta)} \left( \frac{1-2\eta + \frac{p}{2(n+1)} }{(2\eta)^{1/p}} \right)  \right)^2_+ \\
		& =  r^2 \sup_{\eta \in [0,1/2]} \left( (1 - 2 \eta) - \frac{1}{(1 + \delta)}  \frac{1-2\eta  }{(2\eta)^{1/p}} - \frac{1-2\eta}{n(1+\delta)(2\eta)^{1/p}} - \frac{p}{2n(1 + \delta)(2\eta)^{1/p}}  \right)^2_+ \\
		& \geq  r^2 \sup_{\eta \in [1/4,1/2]} \left( (1 - 2 \eta)  - \frac{1}{(1 + \delta)}  \frac{1-2\eta  }{2\eta} - \frac{2 + 2C_1 \log n}{n} \right)_+^2 (\textnormal{as } 1\leq p < C_1 \log n, \eta \in [ 1/4,1/2]). \\
	\end{split}
\end{equation} Let $g(\eta) = (1 - 2 \eta)  - \frac{1}{(1 + \delta)}  \frac{1-2\eta  }{2\eta}$, then $g'(\eta) = -2 + \frac{1}{2(1 + \delta) \eta^2} $. This means that $g(\eta)$ is increasing on $\eta \in [0, \frac{1}{2\sqrt{1 + \delta}} ]$ and decreasing on $\eta \in [ \frac{1}{2\sqrt{1 + \delta}},1/2 ]$. Since $g(1/2) = 0$, it means that $\sup_{\eta \in [1/4,1/2]} g(\eta) $ is some positive constant depending on $\delta$ only. Combining this with \eqref{ineq:lp-stability-lower-2} suggests that as long as $n$ is greater than a sufficiently large constant depending on $\delta$ only, then the right-hand side of \eqref{ineq:lp-stability-lower-2} is lower bounded by $c r^2$ with some $c > 0$ depending on $\delta$ only. 
\end{itemize} This finishes the proof of this theorem.

\subsection{Proof of Theorem \ref{th:worst-bounded-refined} }
Let us begin with the upper bound. We begin by deriving a generic upper bound for an estimator of the form $\wh{\theta} = \frac{\widebar{X}}{1 + \delta}$ for some $\delta \geq 0$.
	\begin{equation} \label{ineq:bounded-class-worst-bound}
	\begin{split}
		\bbE(\wh{\theta} - \theta)^2  &= \bbE\left( \frac{ \bar{X}}{1+\delta} - \theta \right)^2  = \bbE\left( \frac{\bar{X}^2}{(1 + \delta)^2} + \theta^2 - 2 \frac{\theta \wb{X}}{1+ \delta} \right)  = \frac{1}{(1 + \delta)^2} \left(\frac{\var(X)}{n} + \theta^2 \delta^2 \right),
	\end{split}
	\end{equation} where $\theta = \bbE[X]$. Since $|X| \leq r$, $\var(X) + \theta^2 = \bbE|X|^2 \leq r^2$. Thus,
\begin{equation} \label{ineq:generic-bound-1}
		\begin{split}
			\sup_{P \in \cP^1_{\textnormal{bound}}(r)} \bbE(\wh{\theta} - \theta)^2 = \sup_{P \in \cP^1_{\textnormal{bound}}(r)}   \frac{1}{(1 + \delta)^2} \left(\frac{\var(X)}{n} + \theta^2 \delta^2 \right) &\leq \sup_{P \in \cP^1_{\textnormal{bound}}(r)}   \frac{1}{(1 + \delta)^2} \left(\frac{r^2 - \theta^2}{n} + \theta^2 \delta^2 \right) \\
			& = \left\{ \begin{array}{ll}
				\frac{r^2}{n(1 + \delta)^2} & \delta \leq 1/\sqrt{n}\\
				\frac{\delta^2r^2}{(1 + \delta)^2} & \delta > 1/\sqrt{n}.
			\end{array} \right.
		\end{split}
	\end{equation}

For our estimator, we can easily verify that it is $\beta_n$-worst-case-stable as the sample mean is $\frac{2r}{n}$-worst-case-stable. In terms of the error analysis, we divide it into two cases. 
\begin{itemize}[leftmargin=*]
	\item Case 1: $\beta_n \geq  \frac{2r}{n(1 + 1/\sqrt{n})}$. In this case, $\wh{\theta} = \frac{\widebar{X}}{1 + \frac{1}{\sqrt{n}}}$ and by \eqref{ineq:generic-bound-1}, $$\sup_{P \in \cP^1_{\textnormal{bound}}(r)} \bbE(\wh{\theta} - \theta)^2 = \frac{r^2}{(\sqrt{n} + 1)^2} =  \left(\frac{\left( \frac{2r}{n \beta_n} - 1 \right)_{+}}{1 + \left( \frac{2r}{n \beta_n} - 1 \right)_{+} } \right)^2 r^2 \vee \frac{r^2}{\left(\sqrt{n} + 1 \right)^2},$$ where in the second equality, we leverage the observation that $\left(\frac{\left( \frac{2r}{n \beta_n} - 1 \right)_{+}}{1 + \left( \frac{2r}{n \beta_n} - 1 \right)_{+} } \right)^2$ is a decreasing function of $\beta_n$ and when $\beta_n = \frac{2r}{n(1 + 1/\sqrt{n})}$, $\left(\frac{\left( \frac{2r}{n \beta_n} - 1 \right)_{+}}{1 + \left( \frac{2r}{n \beta_n} - 1 \right)_{+} } \right)^2 = \frac{1}{\left(\sqrt{n} + 1 \right)^2}$.
	\item Case 2: $\beta_n < \frac{2r}{n(1 + 1/\sqrt{n})}$. In this case, $\wh{\theta} = \frac{\widebar{X}}{ \frac{2r}{n \beta_n} } = \frac{\widebar{X}}{1 + \frac{2r}{n \beta_n} - 1 }$ and by \eqref{ineq:generic-bound-1}, $$\sup_{P \in \cP^1_{\textnormal{bound}}(r)} \bbE(\wh{\theta} - \theta)^2 = \frac{\left(\frac{2r}{n \beta_n} - 1 \right)^2}{(1 + \frac{2r}{n \beta_n} - 1)^2}r^2 = \left(\frac{\left( \frac{2r}{n \beta_n} - 1 \right)_{+}}{1 + \left( \frac{2r}{n \beta_n} - 1 \right)_{+} } \right)^2 r^2 \vee \frac{r^2}{\left(\sqrt{n} + 1 \right)^2}.$$ 
\end{itemize} 

Now, we move on to the lower bound. Again, we divide the proof into two cases.
\begin{itemize}[leftmargin=*]
	\item Case 1: $\beta_n \geq \frac{2r}{n(1 + 1/\sqrt{n})}$. Notice that in this regime, $ \left(\frac{\left( \frac{2r}{n \beta_n} - 1 \right)_{+}}{1 + \left( \frac{2r}{n \beta_n} - 1 \right)_{+} } \right)^2 \vee \frac{1}{\left(\sqrt{n} + 1 \right)^2} =  \frac{1}{\left(\sqrt{n} + 1 \right)^2}$. Let us consider a subclass $\cP_{\textnormal{binary}} \subseteq \cP^1_{\textnormal{bound}}(r)$, where $\cP_{\textnormal{binary}}$ contains the class of distribution supported on $\{-r,r\}$ and $P(X = r) = p$ for some $p \in [0,1]$. Then it is known that 
	\begin{equation*}
		\inf_{\wh{\theta}}\sup_{P \in  \cP_{\textnormal{binary}}} \bbE(\wh{\theta} - \theta)^2 \geq \frac{r^2}{(\sqrt{n} + 1)^2}
	\end{equation*} using the fact that a Bayes estimator with constant risk is minimax (see Example 4.18 in \cite{shao2006mathematical}). As a result,
	\begin{equation*}
		\begin{split}
			R_{n,\infty}(\theta(\cP^1_{\textnormal{bound}}(r)), \beta_n) \geq \inf_{\wh{\theta}}\sup_{P \in  \cP_{\textnormal{binary}}} \bbE(\wh{\theta} - \theta)^2 \geq \frac{r^2}{(\sqrt{n} + 1)^2}.
		\end{split}
	\end{equation*}
	\item Case 2: $\beta_n < \frac{2r}{n(1 + 1/\sqrt{n})}$. In this regime, $\left(\frac{\left( \frac{2r}{n \beta_n} - 1 \right)_{+}}{1 + \left( \frac{2r}{n \beta_n} - 1 \right)_{+} } \right)^2 r^2\vee \frac{r^2}{\left(\sqrt{n} + 1 \right)^2} =  \left(\frac{ \frac{2r}{n \beta_n} - 1}{1 + \frac{2r}{n \beta_n} - 1 } \right)^2 r^2= \left( \frac{2r - n \beta_n}{2} \right)^2$. By Theorem \ref{th:worst-lower-bound-multi} with $\bbE[ d_{\Ham}(\cD_n, \cD_n') ]$ being replaced by a trivial upper bound $n$, we have 
	\begin{equation*}
		\begin{split}
			R_{n,\infty}(\theta(\cP^1_{\textnormal{bound}}(r)), \beta_n) &\geq \left( \frac{ 2r - n\beta_n  }{2} \right)^2.
		\end{split}
	\end{equation*}
\end{itemize}
In summary, we have shown 
\begin{equation*}
	R_{n,\infty}(\theta(\cP^1_{\textnormal{bound}}(r)), \beta_n) \geq \left(\frac{\left( \frac{2r}{n \beta_n} - 1 \right)_{+}}{1 + \left( \frac{2r}{n \beta_n} - 1 \right)_{+} } \right)^2 r^2 \vee \frac{r^2}{\left(\sqrt{n} + 1 \right)^2}.
\end{equation*}


\subsection{Proof of Theorem \ref{th:average-case-bounded-class}}

We divide the proof into two parts. In Part I, we show the upper bound and in Part II, we show the lower bound.
\vskip.2cm
{\noindent \bf Part I: Upper Bound.} We begin by showing the estimation error guarantee of the proposed estimator and then check it also satisfies the desired average-case stability.
\vskip.2cm
{\bf \noindent (Error Guarantee)} We consider the estimation error based on different values of $\beta_n$.
\begin{itemize}[leftmargin=*]
	\item When $\beta_{n} \leq r/(2n)$, it is easy to check $\wh{\theta} = 0$, thus $\bbE(\wh{\theta} - \theta)^2 \lesssim r^2 \asymp \frac{\left( \frac{r}{n \beta_n} - 1 \right)^3_{+}}{\left(1 + \left( \frac{r}{n \beta_n} - 1 \right)_{+} \right)^3} r^2 \vee \frac{r^2}{n}.$
	\item When $\beta_n \geq r/n$, $\wh{\theta} = \bar{X}
	$. Then $\bbE(\wh{\theta} - \theta)^2 \lesssim  \frac{r^2}{n} \asymp \frac{\left( \frac{r}{n \beta_n} - 1 \right)^3_{+}}{\left(1 + \left( \frac{r}{n \beta_n} - 1 \right)_{+} \right)^3} r^2 \vee \frac{r^2}{n}$.
	\item  When $\beta_n \in [ \frac{r}{2n}, \frac{r}{n} ].$ We can write $ \beta_n = \frac{r}{n(1 + \delta)}$, where $\delta = \frac{r}{n \beta_n} -1 \in [0,1]$. Then our proposed estimator can be written as 
	\begin{equation} \label{eq:average-case-estimator}
		\wh{\theta} = \left( 1- \delta \left(1 \wedge \frac{2(\sqrt{\delta} + 1/\sqrt{n})r}{|\bar{X}|} \right) \right) \bar{X}.
	\end{equation} Then
\begin{equation*}
	\begin{split}
		& \bbE(\wh{\theta} - \theta)^2 \\
		&= \bbE[(\wh{\theta} - \theta)^2 \indi(|\widebar{X}| \leq 2r(\sqrt{\delta} + 1/\sqrt{n}) ) ] + \bbE[(\wh{\theta} - \theta)^2 \indi(|\widebar{X}| > 2r(\sqrt{\delta} + 1/\sqrt{n}) ) ] \\
		& \leq  \bbE\left[\left(\widebar{X}(1 - \delta) - \theta\right)^2 \indi(|\widebar{X}| \leq 2r(\sqrt{\delta} + 1/\sqrt{n}) ) \right] + \bbE\left[\left(\widebar{X} \left(1 - \delta \frac{2(\sqrt{\delta} + 1/\sqrt{n})r}{|\bar{X}|} \right) - \theta\right)^2 \right] \\
		& \lesssim  \bbE[(\widebar{X} - \theta)^2 ]  +  \bbE\left[\delta^2 \widebar{X}^2 \indi(|\widebar{X}| \leq 2r(\sqrt{\delta} + 1/\sqrt{n}) ) \right] + \delta^2 (\sqrt{\delta} + 1/\sqrt{n})^2 r^2 \\
		& \lesssim \frac{r^2}{n} + \delta^3 r^2 \asymp \frac{\left( \frac{r}{n \beta_n} - 1 \right)^3_{+}}{\left(1 + \left( \frac{r}{n \beta_n} - 1 \right)_{+} \right)^3} r^2 \vee \frac{r^2}{n}.
	\end{split}
\end{equation*}
\end{itemize} We note the constant in $\bbE(\wh{\theta} - \theta)^2 \lesssim \frac{\left( \frac{r}{n \beta_n} - 1 \right)^3_{+}}{\left(1 + \left( \frac{r}{n \beta_n} - 1 \right)_{+} \right)^3} r^2 \vee \frac{r^2}{n}$ can be chosen to be independent of $\beta_n$.

\vskip.2cm
{\bf \noindent (Stability Guarantee)} Next, we show that the estimator satisfies the desired average-case stability. We divide the proof based on different cases of $\beta_n$.

\begin{itemize}[leftmargin=*]
	\item Case 1: $\beta_n \geq r/n$. In this case, our estimator is the sample mean and the following lemma shows that the sample mean is $ \frac{r}{n}$-average-case-stable, so it is also $\beta_n$-average-case-stable. 
	\begin{Lemma} \label{lm:average-case-bounded-max}
	Suppose $\wh{\theta}(\cD_n; \xi) = \widebar{X}$, then 
	\begin{equation*}
	\frac{1}{(n+1)^2} \sup_{\cD_{n+1} \in \cX^{n+1}} \sum_{1 \leq i, j \leq n+1} \bbE_{\xi} |\wh{\theta}(\cD_{n+1}^{\setminus i}; \xi )  - \wh{\theta}(\cD_{n+1}^{\setminus j}; \xi )| \leq \frac{r}{n}.
\end{equation*} Moreover, the maximum can be achieved when $\cD_{n+1}$ consists of half $X_i$'s are $r$ and half $X_i$'s are $-r$. 
\end{Lemma}
	
	\item Case 2: $\beta_n < r/(2n)$. In this case, $\widehat{\theta} = 0$ and the stability guarantee clearly holds. 
	
	\item Case 3: $\beta_n \in [\frac{r}{2n}, \frac{r}{n}]$. In this case, let us denote $\delta = \frac{r}{n \beta_n} -1$. Lemma \ref{lm:average-case-bounded-max} suggests the natural shrinage estimator $\wh{\theta} = \frac{n \beta_n}{r} \bar{X} =  \frac{\widebar{X}}{1+ \delta}$ can be $\beta_n$-average-case-stable. In particular, the most unstable scenario appears when half of the data is $-r$ and half data is $r$, or in another way of speaking, when the sample mean is zero. Unfortunately, the bias of this estimator is too large. Our proposed estimator \eqref{eq:average-case-estimator} shrinks the sample mean in a different way based on its magnitude. The rest of the proof is devoted to showing that we can still guarantee the estimator in \eqref{eq:average-case-estimator} is $\beta_n$-average-case-stable. In particular, the following lemma is the key: it suggests that if the magnitude of $\wb{X}$ is large, then actually using the original sample mean can also guarantee the desired average-case stability.
\begin{Lemma} \label{lm:averaged-case-lemma2}
	Given any $\delta \in (0,1)$ and $n \geq 3$, let $T = \argmax\{i\in \bbZ: i \leq  2(n\sqrt{\delta} +  \sqrt{n} )-1, i \textnormal{ is even} \}$. Suppose $\wh{\theta}(\cD_n; \xi) = \widebar{X}$, then 
	\begin{equation}
	\frac{1}{(n+1)^2} \max_{\cD_{n+1} \in \cX^{n+1}: |\sum_{i=1}^{n+1}X_i| \geq Tr } \sum_{1 \leq i, j \leq n+1} \bbE_{\xi} |\wh{\theta}(\cD_{n+1}^{\setminus i}; \xi )  - \wh{\theta}(\cD_{n+1}^{\setminus j}; \xi )| \leq \frac{r}{n} - \frac{2\delta r}{n}.
\end{equation}
\end{Lemma}

Next, we consider two scenarios. Given any dataset $\cD_{n+1}$, without loss of generality, we will assume $-r \leq X_1 \leq \cdots \leq X_{n+1} \leq r$ and we also denote $\wb{\cD}_{n+1}^{\setminus i} = \frac{\sum_{j \in [n+1], j \neq i} X_j }{n}$ for notation convenience. 
	
	{\noindent \bf Scenario 1: $\wb{\cD}_{n+1}^{\setminus n+1} \geq - 2r(\sqrt{\delta} + \frac{1}{\sqrt{n}})$ and $\wb{\cD}_{n+1}^{\setminus 1} \leq 2r(\sqrt{\delta} + \frac{1}{\sqrt{n}})$.} In this scenario, we have
	\begin{equation} \label{eq:scenario1}
		|\wb{\cD}_{n+1}^{\setminus i}| \leq 2r(\sqrt{\delta} + \frac{1}{\sqrt{n}})\quad \textnormal{for all  } i \in [n+1].
	\end{equation}
Thus,
		\begin{equation*}
		\begin{split}
			& \frac{1}{(n+1)^2} \sup_{\cD_{n+1} \in \cX^{n+1} \textnormal{ satisfies } \eqref{eq:scenario1} } \sum_{1 \leq i, j \leq n+1} \bbE_{\xi} |\wh{\theta}(\cD_{n+1}^{\setminus i}; \xi )  - \wh{\theta}(\cD_{n+1}^{\setminus j}; \xi )| \\
	& \overset{ \eqref{eq:average-case-estimator} }= \frac{1-\delta}{(n+1)^2} \sup_{\cD_{n+1} \in \cX^{n+1} \textnormal{ satisfies } \eqref{eq:scenario1}} \sum_{1 \leq i, j \leq n+1} |\wb{\cD}_{n+1}^{\setminus i} - \wb{\cD}_{n+1}^{\setminus j}|\\
	& \leq \frac{1- \delta}{(n+1)^2} \sup_{\cD_{n+1} \in \cX^{n+1}} \sum_{1 \leq i, j \leq n+1} |\wb{\cD}_{n+1}^{\setminus i} - \wb{\cD}_{n+1}^{\setminus j}|\\
	& \overset{\textnormal{Lemma } \ref{lm:average-case-bounded-max} } \leq \frac{(1 - \delta)r}{n} \leq  \frac{r}{n(1 + \delta)} = \beta_n 
		\end{split}
\end{equation*} So our estimator is $\beta_n$-average-case-stable.

{\noindent \bf Scenario 2: $\wb{\cD}_{n+1}^{\setminus n+1} < - 2r(\sqrt{\delta} + \frac{1}{\sqrt{n}})$ or $\wb{\cD}_{n+1}^{\setminus 1} > 2r(\sqrt{\delta} + \frac{1}{\sqrt{n}})$.} We focus on the proof for the scenario 
\begin{equation} \label{eq:scenario-2}
	\wb{\cD}_{n+1}^{\setminus 1} > 2r(\sqrt{\delta} + \frac{1}{\sqrt{n}}),
\end{equation} while the proof for the other case is similar. We first notice that when $\wb{\cD}_{n+1}^{\setminus 1} > 2r(\sqrt{\delta} + \frac{1}{\sqrt{n}})$, then it implies 
\begin{equation} \label{eq:2-implied-condition}
\begin{split}
	&(i): \quad \sum_{i=1}^{n+1} X_i \geq n \wb{\cD}_{n+1}^{\setminus 1} - r > 2r(n \sqrt{\delta} + \sqrt{n} ) - r \\
	&(ii): \quad \wb{\cD}_{n+1}^{\setminus n+1} = \frac{n \wb{\cD}_{n+1}^{\setminus 1} - X_{n+1} + X_1 }{n} \geq 2r(\sqrt{\delta} + \frac{1}{\sqrt{n}}) - \frac{2r}{n} > 0.
\end{split}
\end{equation} This implies that there exists a $k \in [n+1]$ such that 
\begin{equation} \label{eq:smoothness}
\begin{split}
(i):\,&\wb{\cD}_{n+1}^{\setminus i} > 2r(\sqrt{\delta} + \frac{1}{\sqrt{n}}), \text{for all } i =1,\ldots, k\\
	(ii):\,&0 < 2r(\sqrt{\delta} + \frac{1}{\sqrt{n}}) - \frac{2r}{n} \leq \wb{\cD}_{n+1}^{\setminus i} \leq 2r(\sqrt{\delta} + \frac{1}{\sqrt{n}}), \text{ for all } i =k+1,\ldots, n+1
\end{split}
\end{equation}
 for all $i =k+1,\ldots, n+1$. Thus,
\begin{equation*}
	\begin{split}
			& \frac{1}{(n+1)^2} \max_{\cD_{n+1} \in \cX^{n+1} \textnormal{ satisfies } \eqref{eq:scenario-2} } \sum_{1 \leq i, j \leq n+1} \bbE_{\xi} |\wh{\theta}(\cD_{n+1}^{\setminus i}; \xi )  - \wh{\theta}(\cD_{n+1}^{\setminus j}; \xi )| \\
			& = \frac{2}{(n+1)^2} \max_{\cD_{n+1} \in \cX^{n+1} \textnormal{ satisfies } \eqref{eq:scenario-2}} \Bigg( \sum_{1 \leq i < j \leq k} |\wh{\theta}(\cD_{n+1}^{\setminus i}; \xi )  - \wh{\theta}(\cD_{n+1}^{\setminus j}; \xi )| \\
			& \quad + \sum_{1 \leq i \leq k < j \leq n+1} |\wh{\theta}(\cD_{n+1}^{\setminus i}; \xi )  - \wh{\theta}(\cD_{n+1}^{\setminus j}; \xi )|  + \sum_{k+1< i < j \leq n+1} |\wh{\theta}(\cD_{n+1}^{\setminus i}; \xi )  - \wh{\theta}(\cD_{n+1}^{\setminus j}; \xi )| \Bigg)  \\
			& \overset{ \eqref{eq:average-case-estimator}, \eqref{eq:smoothness} }= \frac{2}{(n+1)^2} \max_{\cD_{n+1} \in \cX^{n+1}  \textnormal{ satisfies } \eqref{eq:scenario-2}} \\
			&\Bigg( \sum_{1 \leq i < j \leq k} \left| \bar{\cD}_{n+1}^{\setminus i} \left(1 - \delta \frac{2(\sqrt{\delta} + 1/\sqrt{n})r}{\bar{\cD}_{n+1}^{\setminus i}} \right) - \bar{\cD}_{n+1}^{\setminus j} \left(1 - \delta \frac{2(\sqrt{\delta} + 1/\sqrt{n})r}{\bar{\cD}_{n+1}^{\setminus j}} \right) \right| \\
			& \quad  +  \sum_{1 \leq i \leq k < j \leq n+1}\left|\bar{\cD}_{n+1}^{\setminus i} \left(1 - \delta \frac{2(\sqrt{\delta} + 1/\sqrt{n})r}{\bar{\cD}_{n+1}^{\setminus i}} \right) - \bar{\cD}_{n+1}^{\setminus j}(1 - \delta) \right| + \sum_{k+1< i < j \leq n+1} \frac{|X_j - X_i|}{n}(1 - \delta) \Bigg) \\
			& \overset{\eqref{eq:smoothness}(ii)}\leq   \frac{2}{(n+1)^2} \max_{\cD_{n+1} \in \cX^{n+1}  \textnormal{ satisfies } \eqref{eq:scenario-2}} \left(  \sum_{1 \leq i < j \leq n+1}  \frac{|X_j - X_i|}{n} +  \sum_{1 \leq i \leq k < j \leq n+1} \left| \delta \left( 2(\sqrt{\delta} + 1/\sqrt{n}) -  \bar{\cD}_{n+1}^{\setminus j}\right) \right| \right) \\
			& \overset{\eqref{eq:smoothness}(ii)}\leq  \frac{2}{(n+1)^2} \max_{\cD_{n+1} \in \cX^{n+1}  \textnormal{ satisfies } \eqref{eq:scenario-2}} \left( \sum_{1 \leq i < j \leq n+1} \frac{|X_j - X_i|}{n} \right) + \frac{2\delta}{(n+1)^2} k(n+1-k) \frac{2r}{n}  \\
			& \overset{ \eqref{eq:2-implied-condition} (i) }\leq  \frac{2}{(n+1)^2}  \max_{\cD_{n+1} \in \cX^{n+1}: \sum_{i=1}^{n+1} X_i \geq 2(n\sqrt{\delta} + \sqrt{n} )r - r } \left( \sum_{1 \leq i < j \leq n+1} \frac{|X_j - X_i|}{n} \right) + \frac{\delta r}{n} \\
			& \overset{(a)}\leq \frac{1}{(n+1)^2}  \max_{\cD_{n+1} \in \cX^{n+1}: \sum_{i=1}^{n+1} X_i \geq Tr } \left( \sum_{1 \leq i , j \leq n+1} \frac{|X_j - X_i|}{n} \right) + \frac{\delta r}{n} \\
			& \overset{\textnormal{Lemma \ref{lm:averaged-case-lemma2}} }\leq \frac{r}{n} - \frac{\delta r}{n} \leq \frac{r}{(1 + \delta)n}.
	\end{split}
\end{equation*} 
where (a) is because $T = \argmax\{i\in \bbZ: i \leq  2(n\sqrt{\delta} +  \sqrt{n} )-1, i \textnormal{ is even} \}$.
\end{itemize}

\vskip.2cm
{\noindent \bf Part II: Lower Bound.} To show the lower bound, we need a slightly sharper lower bound than the one in \eqref{ineq:average-lower-linear-fun}. Its proof is provided in the subsequent subsections.

\begin{Theorem}\label{th:sharper-average-lower-bound}
Consider the same setting as in Theorem \ref{th:average-lower-bound-multi}. Suppose $p=1$ and the functional of interest $\theta(P)$ is linear in $P$, i.e., $\theta((1-t) P_0 + tP_1 ) = (1-t)\theta(P_0) + t\theta(P_1)$, then for any $P_0, P_1 \in \cP$ such that $(1-t) P_0 + t P_1 \in \cP$ for all $t \in [0,1]$, we have
\begin{equation*}
		\begin{split}
			R_{n,1}(\theta(\cP), \beta_n) \geq \frac{1}{4} \sup_{\eta \in [1/4,1/2]}\left( (1- 2 \eta) \|\theta(P_0)  - \theta(P_1) \|_2  - n \beta_n \log\left( \frac{1}{\eta} - 1 \right) - 6 \beta_n \right)^2_{+} \vee R_n(\theta(\cP)).
		\end{split}
	\end{equation*}
\end{Theorem}

 We divide the proof into different cases based on the value of $\beta_n$.
\begin{itemize}[leftmargin=*]
	\item When $\beta_n \geq r/n$, then by a standard two-point Le Cam's method \citep{yu1997assouad}, we can show $R_{n,1}\left(\theta(\cP^1_{\textnormal{bound}}(r)), \beta_n \right) \geq R_{n,1}\left(\theta(\cP^1_{\textnormal{bound}}(r)), \infty \right) \gtrsim \frac{r^2}{n} \asymp  \frac{\left( \frac{r}{n \beta_n} - 1 \right)^3_{+}}{\left(1 + \left( \frac{r}{n \beta_n} - 1 \right)_{+} \right)^3} r^2 \vee \frac{r^2}{n}$. For simplicity, we omit the details here.
	\item When $\beta_n \leq r/(2n)$, then by taking $\eta = 1/4$, $\theta(P_0) = r$ and $\theta(P_1) = -r$ in Theorem \ref{th:sharper-average-lower-bound}, we get
	\begin{equation*}
		\begin{split}
			R_{n,1}(\theta(\cP^1_{\textnormal{bound}}(r)), \beta_n) \geq r^2(1 - \log(3)/2 - 3/n)^2_+ \gtrsim r^2 \asymp \frac{\left( \frac{r}{n \beta_n} - 1 \right)^3_{+}}{\left(1 + \left( \frac{r}{n \beta_n} - 1 \right)_{+} \right)^3} r^2 \vee \frac{r^2}{n} .
		\end{split}
	\end{equation*}
	\item When $\beta_n \in [r/(2n),r/n]$, let $\delta = \frac{r}{n \beta_n} - 1 \in [0,1]$. Thus $\beta_n$ can be rewritten as $\frac{r}{n(1+\delta)}$. First, we notice that if $\delta \leq n^{-1/3}$, then $\frac{\left( \frac{r}{n \beta_n} - 1 \right)^3_{+}}{\left(1 + \left( \frac{r}{n \beta_n} - 1 \right)_{+} \right)^3} r^2 \vee \frac{r^2}{n} \asymp \frac{r^2}{n} $ and the lower bound for $R_{n,1}(\theta(\cP^1_{\textnormal{bound}}(r)), \beta_n)$ can be proved in the same way as $\beta_n \leq r/n$. So we just need to focus on the setting $\delta \geq n^{-1/3}$. By Theorem \ref{th:sharper-average-lower-bound} with $\theta(P_0) = r$ and $\theta(P_1) = -r$, we get
	\begin{equation} \label{ineq:average-lower-bound-ineq}
		\begin{split}
			R_{n,1}(\theta(\cP^1_{\textnormal{bound}}(r)), \beta_n) \geq \frac{1}{4} \sup_{\eta \in [1/4,1/2]} r^2\left( 2(1-2\eta) - \frac{\log(1/\eta - 1)}{1+\delta}  - \frac{6}{n} \right)_+^2.
		\end{split}
	\end{equation} Let $f(\eta) =  2(1-2\eta) - \frac{\log(1/\eta - 1)}{1+\delta} $. Then $f'(\eta) = -4 + \frac{1}{(1+\delta)\eta(1-\eta)}$. It is easy to check $f(\eta)$ increases on $\eta \in [0, \eta^*]$ and then decreases on $\eta \in [\eta^*, 1/2]$, where $\eta^* = \frac{1 - \sqrt{ \frac{\delta}{1+ \delta} } }{2} $. Moreover, $f(1/2) = 0$. If $\delta \in [1/3,1]$, i.e., $\eta^* \leq 1/4$, then we plug in $\eta = 1/4$ in \eqref{ineq:average-lower-bound-ineq} and it implies that 
	\begin{equation*}
		\begin{split}
			R_{n,1}(\theta(\cP^1_{\textnormal{bound}}(r)), \beta_n) \geq\frac{1}{4} r^2(1 - \log(3)/(1 + \delta) - 6/n )_+^2 \gtrsim r^2 \asymp \frac{\left( \frac{1}{n \beta_n} - 1 \right)^3_{+}}{\left(1 + \left( \frac{1}{n \beta_n} - 1 \right)_{+} \right)^3} r^2 \vee r^2,
		\end{split}
	\end{equation*} as $n$ is large and $\delta \geq 1/3$. Let $c > 0$ to a small enough constant to be specified later. If $\delta \in [n^{-1/3},c)$, then $\eta^* \in [1/4,1/2]$ and by plugging in $\eta = \eta^*$ in  \eqref{ineq:average-lower-bound-ineq}, we get
	\begin{equation} \label{ineq:case3-lower-bound}
		\begin{split}
			R_{n,1}(\theta(\cP^1_{\textnormal{bound}}(r)), \beta_n) &\geq\frac{1}{4} r^2\left( 2(1-2\eta^*) - \frac{\log(1/\eta^* - 1)}{1+\delta}  - \frac{6}{n} \right)_+^2 \\
			& =\frac{1}{4}r^2\left( 2\sqrt{ \frac{\delta}{1+ \delta} } - \frac{\log\left(1 + \frac{2 \sqrt{ \frac{\delta}{1+ \delta} } }{1 - \sqrt{ \frac{\delta}{1+ \delta} } } \right)}{1+\delta}  - \frac{6}{n} \right)_+^2
		\end{split}
	\end{equation} 
	Notice that \begin{equation*}
		\begin{split}
			\log\left(1 + \frac{2 \sqrt{ \frac{\delta}{1+ \delta} } }{1 - \sqrt{ \frac{\delta}{1+ \delta} } } \right) \leq \frac{2 \sqrt{ \frac{\delta}{1+ \delta} } }{1 - \sqrt{ \frac{\delta}{1+ \delta} } } - \frac{1}{2} \left(\frac{2 \sqrt{ \frac{\delta}{1+ \delta} } }{1 - \sqrt{ \frac{\delta}{1+ \delta} } } \right)^2 + \frac{1}{3} \left(\frac{2 \sqrt{ \frac{\delta}{1+ \delta} } }{1 - \sqrt{ \frac{\delta}{1+ \delta} } } \right)^3. 
		\end{split}
	\end{equation*} due to the fact $\log(1 + x) \leq x - x^2/2 + x^3/3$ holds for all $x >-1$. Plug the above inequality into \eqref{ineq:case3-lower-bound}, we have 
	\begin{equation*}
		\begin{split}
			&2\sqrt{ \frac{\delta}{1+ \delta} } - \frac{\log\left(1 + \frac{2 \sqrt{ \frac{\delta}{1+ \delta} } }{1 - \sqrt{ \frac{\delta}{1+ \delta} } } \right)}{1+\delta} \\
			& \geq 2\sqrt{ \frac{\delta}{1+ \delta} }- \frac{1}{1+\delta} \left( \frac{2 \sqrt{ \frac{\delta}{1+ \delta} } }{1 - \sqrt{ \frac{\delta}{1+ \delta} } } - \frac{1}{2} \left(\frac{2 \sqrt{ \frac{\delta}{1+ \delta} } }{1 - \sqrt{ \frac{\delta}{1+ \delta} } } \right)^2 + \frac{1}{3} \left(\frac{2 \sqrt{ \frac{\delta}{1+ \delta} } }{1 - \sqrt{ \frac{\delta}{1+ \delta} } } \right)^3 \right) \\
			& = 2 \sqrt{ \frac{\delta}{1+ \delta} } \frac{\delta - \sqrt{ \frac{\delta}{1+ \delta} } - \delta \sqrt{ \frac{\delta}{1+ \delta} }}{(1 + \delta)(1 - \sqrt{ \frac{\delta}{1+ \delta} })} + 2 \frac{\delta}{(1 + \delta)^2 \left( 1- \sqrt{ \frac{\delta}{1+ \delta} } \right)^2} - \frac{8 \delta^{3/2} }{3 (1 + \delta)^{5/2} \left( 1- \sqrt{ \frac{\delta}{1+ \delta} } \right)^3 } \\
			& = \frac{2 \delta^{3/2}}{(1 + \delta)^{1.5} \left(1-\sqrt{ \frac{\delta}{1+ \delta} }\right) } \left( 1+\frac{1}{(1 + \delta) \left( 1- \sqrt{ \frac{\delta}{1+ \delta} } \right) } - \frac{4}{3(1 + \delta)\left( 1- \sqrt{ \frac{\delta}{1+ \delta} } \right)^2}  - \frac{\sqrt{\delta}}{(1 + \delta)^{0.5} } \right)  \\
			& \gtrsim \delta^{3/2}(1 - C \sqrt{\delta})
		\end{split}
	\end{equation*} for some universal constant $C > 0$. Thus, we can take $c$ to be small enough so that $C \sqrt{\delta} \leq 1/2$, then we have $2\sqrt{ \frac{\delta}{1+ \delta} } - \frac{\log\left(1 + \frac{2 \sqrt{ \frac{\delta}{1+ \delta} } }{1 - \sqrt{ \frac{\delta}{1+ \delta} } } \right)}{1+\delta} \gtrsim \delta^{3/2}$. By plugging it back to \eqref{ineq:case3-lower-bound}, when $n$ is sufficiently large than, then $R_{n,1}(\theta(\cP^1_{\textnormal{bound}}(r)), \beta_n) \gtrsim \delta^3 r^2 \asymp \frac{\left( \frac{r}{n \beta_n} - 1 \right)^3_{+}}{\left(1 + \left( \frac{r}{n \beta_n} - 1 \right)_{+} \right)^3} r^2 \vee r^2 $ as $\delta \geq n^{-1/3}$.

	Finally, if $\delta \in [c,1/3)$, we still have $\eta^* \in [1/4,1/2]$, thus plugging in $\eta = \eta^*$ in  \eqref{ineq:average-lower-bound-ineq}, we get
	\begin{equation} \label{ineq:case3-lower-bound2}
		\begin{split}
			R_{n,1}(\theta(\cP^1_{\textnormal{bound}}(r)), \beta_n) &\geq\frac{1}{4} r^2\left( 2(1-2\eta^*) - \frac{\log(1/\eta^* - 1)}{1+\delta}  - \frac{6}{n} \right)_+^2 \\
			& =\frac{1}{4}r^2\left( 2\sqrt{ \frac{\delta}{1+ \delta} } - \frac{\log\left(1 + \frac{2 \sqrt{ \frac{\delta}{1+ \delta} } }{1 - \sqrt{ \frac{\delta}{1+ \delta} } } \right)}{1+\delta}  - \frac{6}{n} \right)_+^2 \\
			& \overset{(a)}\gtrsim r^2 \asymp \frac{\left( \frac{r}{n \beta_n} - 1 \right)^3_{+}}{\left(1 + \left( \frac{r}{n \beta_n} - 1 \right)_{+} \right)^3} r^2 \vee r^2 ,
		\end{split}
	\end{equation} where (a) is because $f(1/2) = 0$, so $f(\eta^*)$ is some positive quantity depending on $\delta$ only, which is of constant level when $\delta$ is at a constant level and in addition, we can take $n$ to be large enough.
	\end{itemize}
Finally, we note the constant in $\bbE(\wh{\theta} - \theta)^2 \gtrsim \frac{\left( \frac{r}{n \beta_n} - 1 \right)^3_{+}}{\left(1 + \left( \frac{r}{n \beta_n} - 1 \right)_{+} \right)^3} r^2 \vee \frac{r^2}{n}$ can be chosen to be independent of $\beta_n$.

\subsubsection{Proofs of Lemma \ref{lm:average-case-bounded-max} }
 Without loss of generality, we can assume $-r\leq X_1 \leq X_2 \leq \cdots \leq X_{n+1} \leq r$. Then given any $\cD_{n+1}$,
\begin{equation} \label{eq:average-case-exact-calculation}
	\begin{split}
			& \frac{1}{(n+1)^2} \sum_{1 \leq i, j \leq n+1} \bbE_{\xi} |\wh{\theta}(\cD_{n+1}^{\setminus i}; \xi )  - \wh{\theta}(\cD_{n+1}^{\setminus j}; \xi )| 
			=  \frac{2}{(n+1)^2} \sum_{1 \leq i < j \leq n+1} \frac{X_j - X_i}{n} \\
			= & \frac{2}{n(n+1)^2} \sum_{i=0}^n (-n + 2i ) X_{i+1}.  
	\end{split}
\end{equation} When $n$ is odd, the above quantity is maximized when $X_1 = \cdots = X_{(n-1)/2+1} = -r$ and $X_{(n-1)/2+2} = \ldots =X_{n+1} = r$ and when $n$ is even, the above quantity is maximized when $X_1 = \ldots = X_{n/2} = -r$ and $X_{n/2+1} = \ldots = X_{n+1} = r$. The maximum values are $\frac{r}{n}$ (when $n$ is odd) and $ \frac{r}{n+1}\left( 1 + \frac{1}{n+1} \right)$ (when $n$ is even) and both of them are bounded by $\frac{r}{n}$.

\subsubsection{Proof of Lemma \ref{lm:averaged-case-lemma2} }
	Without loss of generality, let us assume $r = 1$ and given any $\cD_{n+1}$, assume $-1 \leq X_1 \leq \cdots \leq X_{n+1} \leq 1$. In addition, we only need to consider the case $T \leq n+1$; otherwise, there is no valid dataset that satisfies the constraint. Moreover, we will just consider the case that $\sum_{i=1}^{n+1}X_i \geq T$, the proof for the case $\sum_{i=1}^{n+1}X_i \leq -T$ is similar. Following the same computation as in \eqref{eq:average-case-exact-calculation}, when $\wh{\theta}(\cD_n; \xi) = \widebar{X}$,
	\begin{equation} \label{ineq:constrianed-opt}
		\begin{split}
			&\frac{1}{(n+1)^2} \max_{\cD_{n+1} \in \cX^{n+1}: \sum_{i=1}^{n+1}X_i \geq T } \sum_{1 \leq i, j \leq n+1} \bbE_{\xi} |\wh{\theta}(\cD_{n+1}^{\setminus i}; \xi )  - \wh{\theta}(\cD_{n+1}^{\setminus j}; \xi )| \\
			&=  \frac{2}{n(n+1)^2} \max_{\cD_{n+1} \in \cX^{n+1}: \sum_{i=1}^{n+1}X_i \geq T } \sum_{i=0}^n (-n + 2i ) X_{i+1}.
		\end{split}
	\end{equation} Next, we focus on the case $n$ is odd while the proof for $n$ is even is similar. Let $\eta = T/2$, $k_1 = \frac{n+3}{2}  - \eta$ and $k_2 = \frac{n+1}{2}  + \eta$. Since $T$ is even, $\eta, k_1$ and $k_2$ are all integers. We divide the set $\{X_i\}_{i=1}^{n+1}$ into three groups:
	\begin{equation*}
		\begin{split}
			\cD_{n+1} = \cX_{\textnormal{bottom}} \cup  \cX_{\textnormal{bulk}} \cup \cX_{\textnormal{top}}
		\end{split}
	\end{equation*} where $\cX_{\textnormal{bottom}} = \{X_i: i = 1, \ldots, k_1-1\}$, $\cX_{\textnormal{bulk}} = \{X_i: i = k_1, \ldots, k_2\}$ and $\cX_{\textnormal{top}} = \{X_i: i = k_2+1, \ldots, n+1\}$. By construction $|\cX_{\textnormal{bottom}}| = |\cX_{\textnormal{top}} | = \frac{n+1}{2} - \eta $ and $|\cX_{\textnormal{bulk}}| = T$. Next, we are going to show that for the constrained optimization problem in \eqref{ineq:constrianed-opt}, i.e., 
	\begin{equation*}
	\begin{split}
		&\max_{  \substack{X_1, \ldots, X_{n+1} \in [-1,1]\\ \sum_{i=1}^{n+1}X_i \geq T } } f(\{X_i\}_{i=1}^{n+1})\\
		& \quad \quad \quad := \frac{2}{n(n+1)^2} \left( -n X_1 - (n-2)X_2 - \cdots - X_{\frac{n-1}{2}} + X_{\frac{n+1}{2} } + \cdots + n X_{n+1} \right),
	\end{split}
	\end{equation*} the maximum value is achieved when all $X_i$ values in $\cX_{\textnormal{bottom}}$ are $-1$ and all $X_i$ values in $\cX_{\textnormal{bulk}} \cup \cX_{\textnormal{top}}$ are $1$. We denote this particular configuration as $\cX_{\textnormal{bottom}}^*$, $\cX_{\textnormal{bulk}}^* $ and $\cX_{\textnormal{top}}^*$. It is easy to check that this configuration satisfies the constraints as $\sum_{i=1}^{n+1}X_i  = T$, and with this configuration,
	\begin{equation*}
		\begin{split}
			f(\cX_{\textnormal{bottom}}^* \cup \cX_{\textnormal{bulk}}^* \cup \cX_{\textnormal{top}}^*) & \overset{\textnormal{Lemma \ref{lm:average-case-bounded-max}}} \leq \frac{1}{n} -  \frac{2}{n(n+1)^2} \times 2 \times \sum_{i=1}^{\eta} (2i-1) = \frac{1}{n} - \frac{4 \eta^2}{n (n+1)^2}  \\
			& \leq \frac{1 }{n} - \frac{4 \delta n}{(n+1)^2} \quad(\textnormal{as } \eta = T/2 \geq n \sqrt{\delta} \textnormal{ when } n \geq 3) \\
			& \leq \frac{1 }{n} \left( 1- 2 \delta \right).
		\end{split}
	\end{equation*}
	
	So if we can prove the claim regarding the constrained optimization problem in \eqref{ineq:constrianed-opt}, then we are done. We will prove the claim by showing that for any other configuration, the objective value will not be greater than the one achieved by $\cX_{\textnormal{bottom}}^* \cup \cX_{\textnormal{bulk}}^* \cup \cX_{\textnormal{top}}^*$. We divide the proof into two steps.
	\begin{itemize}[leftmargin=*]
		\item Step 1: first we observe that for any configuration such that values in $\cX_{\textnormal{bottom}}^*$ are not all $-1$ and values in $\cX_{\textnormal{top}}^*$ are not all $1$, we can always keep increasing the objective by decreasing values in $\cX_{\textnormal{bottom}}^*$ and increasing values in $\cX_{\textnormal{top}}^*$ by the same amount, while we still satisfies the constraint. This means that the optimal configuration satisfies at least one of the conditions: (i) all values in $\cX_{\textnormal{bottom}}^*$ are $-1$ or (ii) all values in $\cX_{\textnormal{top}}^*$ are $1$. 
		\item Step 2 (i): if one configuration satisfies that all values in $\cX_{\textnormal{bottom}}^*$ are $-1$ and also satisfies the constraint, then we argue that the rest of the values in $\cX_{\textnormal{bulk}}^* \cup \cX_{\textnormal{top}}^*$ must be $1$. If this is not the case, i.e, not all values in $\cX_{\textnormal{bulk}}^* \cup \cX_{\textnormal{top}}^*$ are $1$, then $\sum_{i=1}^{n+1}X_i <  2 \eta = T$, which does not satisfy the constraint. 
		\item Step 2 (ii): if one configuration satisfies that all values in $\cX_{\textnormal{top}}^*$ are $1$, then we argue that the objective is no bigger than the optimal configuration. To see this, first we observe that since in the objective \eqref{ineq:constrianed-opt}, the absolute value of the coefficient in front of $X_i$s in $\cX_{\textnormal{bottom}}$ is greater than the ones in $\cX_{\textnormal{bulk}} $, so we can always increases the values in $\cX_{\textnormal{bulk}}$, while decreasing the values in $\cX_{\textnormal{bottom}}$ by the same amount so that the constraint is still satisfied. This process can continue until either (i) all values in $\cX_{\textnormal{bulk}}$ are $1$ or (ii) all values in $\cX_{\textnormal{bottom}}$ are $-1$. If it is the case (ii), then this configuration can satisfy the constraint only when all values in $\cX_{\textnormal{bulk}}$ are all $1$, which is the same as the optimal configuration. If it is the case (i), then we can keep increasing the objective value by decreasing the values in $\cX_{\textnormal{bottom}}$ until they are all $-1$, and in that case the constraint can still be satisfied. Moreover, this configuration is the same as the optimal configuration.  
	\end{itemize} 
	In summary, we have proved the claim and this finishes the proof of this lemma.

\subsection{Proof of Theorem \ref{th:sharper-average-lower-bound} }
The proof of this theorem is very similar to the Part I analysis of Theorem \ref{th:average-lower-bound-multi}. For any $\eta \in [0,1/2]$, we still consider the mixture distribution $P = \eta P_0 + (1 - \eta) P_1 $ and $P' = \eta P_1 + (1 - \eta) P_0 $. Following the same proof as in Part I proof of Theorem \ref{th:average-lower-bound-multi} until equation \eqref{ineq:exp-diff-upper}, by plugging $p=1$ in \eqref{ineq:exp-diff-upper}, we get
\begin{equation} \label{ineq:average-sharper1} 
	\|\bbE[ \wh{\theta}(\cD_n;\xi) - \wh{\theta}(\cD'_n; \xi) ] \|_2 \leq (n+1) \beta_n \bbE  \left( \left(\log \left( \frac{n - T_1 + T_2}{T_2 +1} \right) + \frac{1}{T_2 +1} \right) \right).
\end{equation} Next, we will change the analysis in the equation \eqref{ineq:upper-bound-diff}. In particular, instead of applying \eqref{eq:binomial-exp-1}, we apply \eqref{eq:binomial-exp-2} to upper bound the above quantity (replacing $q$ by $2\eta$):
\begin{equation}
\begin{split}
	\|\bbE[ \wh{\theta}(\cD_n;\xi) - \wh{\theta}(\cD'_n; \xi) ] \|_2 &\overset{ \eqref{ineq:average-sharper1}, \eqref{eq:binomial-exp-2} } \leq (n+1) \beta_n \bbE  \left( \log\left( \frac{n}{\eta(n+1)} - 1 \right) + \frac{1}{ \eta (n+1)}\right) \\
	& \leq n \beta_n \log\left( \frac{1}{\eta} - 1 \right) + \beta_n/\eta + \beta_n \log(1/\eta - 1) \\
	& \leq n \beta_n \log\left( \frac{1}{\eta} - 1 \right)  + 6 \beta_n,
\end{split}
\end{equation} where the last inequality is because $\eta \in [1/4,1/2]$.

	Combining with \eqref{ineq:lower-bound}, we have that 
	\begin{equation*}
		\begin{split}
			(1 - 2 \eta) \|\theta(P_0)  - \theta(P_1)\|_2 - 2\sigma_n \leq  n \beta_n \log\left( \frac{1}{\eta} - 1 \right) + 6 \beta_n \\
			 \Longrightarrow \sigma_n \geq \frac{1}{2} \left( (1- 2 \eta) \|\theta(P_0)  - \theta(P_1)\|_2  - n \beta_n \log\left( \frac{1}{\eta} - 1 \right) - 6 \beta_n \right)_{+}.
		\end{split}
	\end{equation*}  Since the above lower bound holds for any estimator satisfying the stability constraint and any $\eta \in [1/4,1/2]$, we then obtain the desired result.

\section{Proofs for Section \ref{sec:heavy-tailed-mean} }
\subsection{Proof of Theorem \ref{th:worst-heavy-tail} Part (a) }
We first consider the lower bound. We note that the lower bound for $\frac{r^2}{n^{1 \wedge 2(1-1/k) }}$ is from the unconstrained minimax risk lower bound. When $k \geq 2$, it can be proved by Le Cam's two-point argument, and when $k \in [1,2)$, it is proved in \cite[Theorem 3.1]{devroye2016sub}. Next, we show $R_{n,\infty}\left(\theta(\cP^{d}_{k}(r)),  \beta_n \right) \gtrsim   \frac{r^{2k}}{(n\beta_n)^{2(k-1)}}\wedge r^2$. Since the unconstrained minimax risk is already $c r^2$ when $k =1$, we focus on $k > 1$. We will consider $d = 1$ while for the multivariate setting, we just need to consider one coordinate. Given any $\epsilon \in [0,1]$, let $P_0$ and $P_1$ be two distributions such that 
\begin{equation*}
	P_0(X= 0) = 1 \quad \textnormal{ and }\quad P_1(X = x) = \left\{ \begin{array}{ll}
		1 - \epsilon & \text{ if } x = 0\\
		\epsilon & \text{ if } x = r \epsilon^{-1/k}. 
	\end{array} \right.
\end{equation*} It is easy to check $P_0, P_1 \in \cP^{d}_{k}(r)$. In addition, $\|\theta(P_0) - \theta(P_1)\|_2 = r \epsilon^{1-1/k}$. Notice that $\TV(P_0, P_1) = \epsilon$ and by a standard result in probability theory regarding maximal coupling, e.g., see \cite[Lemma 4.1.13]{roch2024modern}, we can find a maximum coupling between $P_0$ and $P_1$ such that $\bbE[ d_{\Ham}(\cD_n, \cD_n') ] = n \TV(P_0, P_1) = n \epsilon$, where marginally, $\cD_n$ and $\cD_n'$ consist $n$ i.i.d. samples from $P_0$ and $P_1$, respectively. Then by Theorem \ref{th:worst-lower-bound-multi}, we get
\begin{equation*}
	R_{n,\infty}\left(\theta(\cP^{d}_{k}(r)),  \beta_n \right) \gtrsim \left( r \epsilon^{1 - 1/k} - n \epsilon \beta_n \right)_+^2.
\end{equation*} When $\beta_n \leq (1 - 1/k )\frac{r}{n} $, then take $\epsilon = 1$, we get $R_{n,\infty}\left(\theta(\cP^{d}_{k}(r)),  \beta_n \right) \gtrsim r^2 \asymp  \frac{r^{2k}}{(n\beta_n)^{2(k-1)}}\wedge r^2$. When $\beta_n >  (1 - 1/k )\frac{r}{n}$, the function $f(\epsilon) = r \epsilon^{1 - 1/k} - n \epsilon \beta_n$ achieves its maximum when $\epsilon = \left( \frac{r(1-1/k)}{n \beta_n} \right)^k $ and the maximum value is $\frac{(1-1/k)^{k-1} r^k}{k (n \beta_n)^{k-1}}$. This shows $R_{n,\infty}\left(\theta(\cP^{d}_{k}(r)),  \beta_n \right) \gtrsim \frac{r^{2k}}{(n\beta_n)^{2(k-1)}} \asymp  \frac{r^{2k}}{(n\beta_n)^{2(k-1)}}\wedge r^2$ and finishes the proof for the lower bound. We note the constant in $\gtrsim$ is independent of $\beta_n$.

Next, we consider the upper bound. Notice that since $\|Y_i\|_2 \leq \frac{n \beta_n}{2}$, it is easy to check that the proposed estimators are $\beta_n$-worst-case-stable.  Next,
\begin{equation} \label{ineq:bias-variance-deco}
	\begin{split}
		\bbE\left[\| \bar{Y}- \theta\|_2^2 \right] =  \frac{\bbE\|Y - \bbE[Y]\|^2_2 }{n} + \left\| \bbE[Y ]  - \theta \right\|_2^2.
	\end{split}
\end{equation} Now we bound the bias in the above equation.
\begin{equation} \label{ineq:heavy-bias}
		\begin{split}
			\left\| \bbE[Y ]  - \theta \right\|_2&= \left\| \bbE[X \indi(\|X\|_2 \leq \rho_k)]  - \theta \right\|_2 = \| \bbE(X \indi(\|X\|_2 > \rho_k)) \|_2 \leq\bbE[\| X\|_2 \indi(\|X\|_2 > \rho_k)] \\
			&\overset{(a)}\leq (\bbE\|X\|_2^k)^{1/k} (\bbE(\indi(\|X\|_2 > \rho_k)))^{1-1/k} \\
			 & \overset{(b)}\leq r\left( \bbP(\|X\|_2 > \rho_k ) \right)^{1-1/k} \\
			 & \overset{(c)}\leq r\left( \frac{\bbE\|X\|_2^k}{\rho^k_k} \right)^{1-1/k} \leq r\left(\frac{r}{\rho_k}\right)^{k-1},
		\end{split}
	\end{equation} here (a) is by Holder's inequality, (b) is by the assumption $P \in \cP^{d}_{k}(r)$ and (c) is by Markov's inequality. At the same time,
	\begin{equation*}
		\left\| \bbE[Y ]  - \theta \right\|_2 =  \left\| \bbE[X \indi(\|X\|_2 \leq \rho_k)]  - \theta \right\|_2  \leq \rho_k + r
	\end{equation*} as $\|\theta\|_2 \leq r$. In summary,
	\begin{equation}\label{ineq:heavy-bias2}
		\left\| \bbE[Y ]  - \theta \right\|_2 \leq \left( r\left(\frac{r}{\rho_k}\right)^{k-1} \right) \wedge (\rho_k + r).
	\end{equation}

	Now we bound the variance term in \eqref{ineq:bias-variance-deco}. When $k \geq 2$, 
	\begin{equation} \label{ineq:heavy-variance0}
		\bbE\|Y - \bbE[Y]\|^2_2 \leq \bbE\|X\|_2^2\leq (\bbE\|X\|_2^k)^{2/k} \leq r^2
	\end{equation} by Holder's inequality and the fact $P \in \cP^{d}_{k}(r)$. When $k \in [1,2)$, 
	\begin{equation}\label{ineq:heavy-variance}
		\begin{split}
			\bbE\|Y - \bbE[Y]\|^2_2\leq \bbE\|Y\|_2^2 & = \bbE [\|X\|_2^{2-k} \indi(\|X\|_2 \leq \rho_k) \|X\|_2^k ] \\
			& \leq \rho_k^{2-k}\bbE [ \indi(\|X\|_2 \leq \rho_k)\|X\|_2^k ] \\
			& \leq \rho_k^{2-k} r^k.
		\end{split}
	\end{equation} When $k \geq 2$, by plugging in $\rho_k $ value into \eqref{ineq:heavy-bias2}, \eqref{ineq:heavy-variance0} and \eqref{ineq:bias-variance-deco}, we obtain $$\bbE\left[\| \bar{Y}- \theta\|_2^2 \right] \leq \frac{r^2}{n} + \left(\left( r\left(\frac{r}{n \beta_n/2}\right)^{k-1} \right) \wedge (\frac{n\beta_n}{2} + r) \right)^2 \asymp \frac{r^2}{n} + \frac{r^{2k}}{(n\beta_n)^{2(k-1)}}\wedge r^2  .$$ When $k \in [1,2)$, we have
    \begin{equation}
        \bbE\left[\| \bar{Y}- \theta\|_2^2 \right] \leq \frac{\rho_k^{2-k} r^k}{n} + \left( \frac{r^k}{\rho_k^{k-1}} \wedge (\rho_k + r) \right)^2.
    \end{equation} When $\beta_n \geq 2r n^{1/k -1}$, then $\rho_k = rn^{1/k}$ and $\bbE\left[\| \bar{Y}- \theta\|_2^2 \right] \leq \frac{2r^2}{n^{2(1-1/k)}} \asymp \frac{r^{2k}}{(n\beta_n)^{2(k-1)}}\wedge r^2 + \frac{r^2}{n^{ 2(1-1/k) }} $. When $\beta_n < 2r n^{1/k -1}$, then $\rho_k = \frac{n \beta_n}{2}$ and 
    \begin{equation*}
         \bbE\left[\| \bar{Y}- \theta\|_2^2 \right] \leq \frac{\rho_k^{2-k} r^k}{n} + \frac{r^{2k}}{\rho_{k}^{2(k-1)}} \wedge (\rho_k + r)^2 \asymp \frac{r^{2k}}{(n\beta_n)^{2(k-1)}}\wedge r^2 + \frac{r^2}{n^{ 2(1-1/k) }} .
    \end{equation*} This finishes the proof.

\subsection{Proof for Theorem \ref{th:worst-heavy-tail} Part (b) }
First, we note that the lower bound for $\beta_n \geq \frac{24r}{n}$ is from the unconstrained minimax risk lower bound. When $k \geq 2$, it can be proved by Le Cam's two point argument and when $k \in [1,2)$, it is proved in \cite[Theorem 3.1]{devroye2016sub}. When $\beta_n \leq \frac{r}{10n}$, the lower bound is obtained from Corollary \ref{coro:average-case-stability}.

Next, we take a look at the estimation error and stability guarantees of $\wh{\theta}$.

\vskip.2cm 
{\bf \noindent Estimation Error Guarantee.} Notice that by definition $\|\wb{\theta}\|_2 \leq 2r$, so when $\beta_n \leq \frac{r}{10n} $, we have $\bbE\|\wh{\theta} - \theta\|_2^2 \lesssim r^2$. Now let us focus on $\beta_n \geq \frac{24r}{n} $, in this case our estimator is $\wh{\theta} = \wb{\theta}$ and its error guarantee is given in the following lemma.

\begin{Lemma} \label{lm:average-heavy-accuracy}
	Given any $P \in \cP^{d}_{k}(r)$ with $k \geq 1$, $\bbE\|\wh{\theta} - \theta\|_2^2 \leq \frac{11r^2}{n^{1 \wedge 2(1-1/k) }}$.
\end{Lemma}
\begin{proof}[Proof of Lemma \ref{lm:average-heavy-accuracy}]
	\begin{equation} \label{ineq:heavy-tailed-estimation-error}
		\begin{split}
			&\bbE\|\wh{\theta} - \theta\|_2^2 \\
			&= \bbE\left[\|\wh{\theta} - \theta\|_2^2 \indi\left\{ \sum_{i=1}^n \|Y_i\|_2\leq 2nr \right\} \right] + \bbE\left[\|\wh{\theta} - \theta\|_2^2 \indi\left\{ \sum_{i=1}^n \|Y_i\|_2> 2nr \right\} \right] \\
			& \overset{(a)}\leq \bbE\left[\| \bar{Y}- \theta\|_2^2 \indi\left\{ \sum_{i=1}^n \|Y_i\|_2 \leq 2nr \right\} \right] + 9r^2 \bbP\left(\sum_{i=1}^n \|Y_i\|_2 > 2nr \right) \\
			& \leq \bbE\left[\| \bar{Y}- \theta\|_2^2 \right]  + 9r^2 \bbP\left(\sum_{i=1}^n \|Y_i\|_2 - n \bbE[\|Y\|_2] > 2nr - nr \bbE[\|Y\|_2] \right) \\
			& \overset{(b)}\leq \frac{\bbE\|Y - \bbE[Y]\|^2_2 }{n} + \left\| \bbE[Y ]  - \theta \right\|_2^2  + 9r^2 \bbP\left(\sum_{i=1}^n \|Y_i\|_2- n \bbE[\|Y\|_2] > nr \right) \\
			& \overset{(c)}\leq \frac{\bbE\|Y - \bbE[Y]\|^2_2 }{n} + \left\| \bbE[Y ]  - \theta \right\|_2^2   + 9r^2 \cdot \frac{\textnormal{Var}(\|Y\|_2)}{nr^2},
		\end{split}
	\end{equation} where (a) is by the definition of $\wh{\theta}$ and the fact that on the event $\left\{ \sum_{i=1}^n \|Y_i\|_2 > 2nr \right\}$, $\|\wh{\theta}\|_2 \leq 2r$ almost surely and $\|\theta\|_2 \leq r$ since $P \in \cP^{d}_{k}(r)$, (b) is because $\bbE\|Y\|_2 \leq \bbE\|X\|_2 \leq r$ and (c) is by Chebyshev's inequality.
	
	Notice that the upper bound for $\left\| \bbE[Y ]  - \theta \right\|_2^2$ has been derived in \eqref{ineq:heavy-bias2} and the upper bound for the variance term $\bbE\|Y - \bbE[Y]\|^2_2$ has been derived in \eqref{ineq:heavy-variance0} for $k \geq 2$ and \eqref{ineq:heavy-variance} for $k \in [1,2)$. In addition, following the same analysis as in \eqref{ineq:heavy-variance0} and \eqref{ineq:heavy-variance}, we can also obtain the same upper bound for $\textnormal{Var}(\|Y\|_2)$ when $k \geq 2$ and $k \in [1,2)$. In summary, by plugging \eqref{ineq:heavy-bias2}, \eqref{ineq:heavy-variance0}, \eqref{ineq:heavy-variance} and the value of $\rho_k$ into \eqref{ineq:heavy-tailed-estimation-error}, we have
	\begin{equation*}
		\begin{split}
			\bbE\|\wh{\theta} - \theta\|_2^2  \leq \frac{11r^2}{n^{1 \wedge 2(1-1/k) }}. 
		\end{split}
	\end{equation*} This finishes the proof.
\end{proof}

\vskip.2cm 
{\bf \noindent Stability Guarantee.} First, we note that by the definition of $\wh{\theta}$, we see that to show the statement, it is enough to show that
		\begin{equation*}
		\frac{1}{(n+1)^2} \sup_{\cD_{n+1}} \sum_{1 \leq i, j \leq n+1}  \|\tilde{\theta}(\cD_{n+1}^{\setminus i} )  - \tilde{\theta}(\cD_{n+1}^{\setminus j} )\|_2 \leq \frac{24r}{n+1},
	\end{equation*} where $\tilde{\theta} = \left(1 \wedge \frac{2nr}{\sum_{i=1}^n\|X_i\|_2} \right)\bar{X}$.
The following lemma shows that this is the case, and this also finishes the proof of Theorem \ref{th:worst-heavy-tail} Part (b).

\begin{Lemma} \label{lm:average-stability-bound}
If given $\cD_n = \{X_1, \ldots, X_{n}\}$, $\wh{\theta}(\cD_n) = \left(1 \wedge \frac{2nr}{\sum_{i=1}^n\|X_i\|_2} \right)\bar{X}$ for some $r > 0$. Then given any dataset $\cD_{n+1} = \{X_1, \ldots, X_{n+1}\}$,
	\begin{equation} \label{ineq:average-stability-heavy-tailed}
		\frac{1}{(n+1)^2} \sup_{\cD_{n+1}} \sum_{1 \leq i, j \leq n+1}  \|\wh{\theta}(\cD_{n+1}^{\setminus i} )  - \wh{\theta}(\cD_{n+1}^{\setminus j} )\|_2 \leq \frac{24r}{n+1}.
	\end{equation}
\end{Lemma}

\begin{proof}[Proof of Lemma \ref{lm:average-stability-bound}]
		First, without loss of generality, we assume $r = 1$. For notation convenience, we will denote $\|\cD_{n+1}^{\setminus i} \|_1 = \sum^{n+1}_{z = 1, z \neq i} \|X_z\|_2$. In addition, without loss of generality, we can assume $\|X_1\|_2 \leq \ldots \leq \|X_{n+1}\|_2$. We will divide the proof into different cases.
\begin{itemize}
		\item (Case 1: $\|\cD_{n+1}^{\setminus i} \|_1 \leq 2n$ for all $i \in [n+1]$) Then
		\begin{equation*}
			\begin{split}
				&\frac{1}{(n+1)^2} \sum_{1 \leq i, j \leq n+1} \|\wh{\theta}(\cD_{n+1}^{\setminus i} )  - \wh{\theta}(\cD_{n+1}^{\setminus j} )\|_2 = \frac{2}{(n+1)^2} \sum_{1 \leq i < j \leq n+1} \frac{\|X_j - X_i\|_2}{n} \\
				 & \leq \frac{2}{(n+1)^2} \sum_{1 \leq i < j \leq n+1} \frac{\|X_j\|_2 + \|X_i\|_2}{n} = \frac{2}{(n+1)^2} \frac{n (\sum_{i=1}^{n+1} \|X_i\|_2)}{n} \leq \frac{2 (\|\cD_{n+1}^{\setminus 1} \|_1 + \|\cD_{n+1}^{\setminus 2} \|_1 )}{(n+1)^2} \\
				 & \leq \frac{8}{n+1}.
			\end{split}
		\end{equation*}

		\item (Case 2: $\|\cD_{n+1}^{\setminus i} \|_1 > 2n$ for all $i \in [n+1]$) 
		\begin{equation} \label{ineq:case-2-bound}
			\begin{split}
				&\frac{1}{(n+1)^2} \sum_{1 \leq i, j \leq n+1} \|\wh{\theta}(\cD_{n+1}^{\setminus i} )  - \wh{\theta}(\cD_{n+1}^{\setminus j} )\|_2  = \frac{2}{(n+1)^2} \sum_{1 \leq i < j \leq n+1} \|\wh{\theta}(\cD_{n+1}^{\setminus i} )  - \wh{\theta}(\cD_{n+1}^{\setminus j} )\|_2 \\
				& =  \frac{2}{(n+1)^2}  \sum_{1 \leq i < j \leq n} \left\| \frac{2( \sum_{z \neq i} X_z )}{\|\cD_{n+1}^{\setminus i}\|_1}  - \frac{2( \sum_{z \neq j} X_z )}{\|\cD_{n+1}^{\setminus j}\|_1} \right\|_2  \\
				& \leq \frac{4}{(n+1)^2}  \sum_{1 \leq i < j \leq n} \left( \underbrace{\left\|(\sum_{z \neq i,j} X_z) \left( \frac{1}{\|\cD^{\setminus i}_{n+1}\|_1} - \frac{1}{\|\cD^{\setminus j}_{n+1}\|_1} \right) \right\|_2}_{(A)} +  \underbrace{\left\| \frac{X_j}{\|\cD_{n+1}^{\setminus i}\|_1}  - \frac{X_i}{\|\cD_{n+1}^{\setminus j}\|_1} \right\|_2}_{(B)}  \right).
			\end{split}
		\end{equation} Next, we bound (A) and (B) separately. 
		\begin{equation} \label{ineq:case2-A}
			\begin{split}
				(A) &\leq \|\cD_{n+1}^{\setminus i,j}\|_1 \frac{| \|X_i\|_2 - \|X_j\|_2 |}{\|\cD_{n+1}^{\setminus i}\|_1 \|\cD_{n+1}^{\setminus j}\|_1} \\
				& \leq \|\cD_{n+1}^{\setminus i,j}\|_1 \frac{\|X_i\|_2 + \|X_j\|_2 }{(\|\cD_{n+1}^{\setminus i,j}\|_1 + \|X_j\|_2) (\|\cD_{n+1}^{\setminus i,j}\|_1 + \|X_i\|_2 )} \\
				& \leq \frac{\|X_i\|_2 + \|X_j\|_2 }{\|\cD_{n+1}^{\setminus i,j}\|_1 + \|X_i\|_2 \vee \|X_j\|_2 } =  \frac{2(\|X_i\|_2 + \|X_j\|_2) }{2\|\cD_{n+1}^{\setminus i,j}\|_1 + 2(\|X_i\|_2 \vee \|X_j\|_2) } \\
				& \leq \frac{2(\|X_i\|_2 + \|X_j\|_2) }{\|\cD_{n+1}\|_1}.
			\end{split}
		\end{equation} At the same time,
		\begin{equation} \label{ineq:case2-B}
			\begin{split}
				(B)  &= \left\|  \frac{X_j \|\cD_{n+1}^{\setminus j}\|_1 - X_i \|\cD_{n+1}^{\setminus i}\|_1}{\|\cD_{n+1}^{\setminus i}\|_1 \|\cD_{n+1}^{\setminus j}\|_1 }  \right\|_2  =   \left\|  \frac{X_j (\|\cD_{n+1}^{\setminus i,j}\|_1 + \|X_i\|_2 ) - X_i (\|\cD_{n+1}^{\setminus i,j}\|_1 + \|X_j\|_2 )}{\|\cD_{n+1}^{\setminus i}\|_1 \|\cD_{n+1}^{\setminus j}\|_1 }  \right\|_2 \\
				& \leq \|\cD_{n+1}^{\setminus i,j}\|_1\frac{\|X_j\|_2 + \|X_i\|_2}{\|\cD_{n+1}^{\setminus i}\|_1 \|\cD_{n+1}^{\setminus j}\|_1 } + \frac{2\|X_i\|_2 \|X_j\|_2}{\|\cD_{n+1}^{\setminus i}\|_1 \|\cD_{n+1}^{\setminus j}\|_1} \\
				& \overset{ \eqref{ineq:case2-A}}\leq  \frac{2(\|X_i\|_2 + \|X_j\|_2) }{\|\cD_{n+1}\|_1} + \frac{2\|X_i\|_2}{\|\cD_{n+1}^{\setminus j}\|_1 } \\
				& \overset{(a)}\leq \frac{2(\|X_i\|_2 + \|X_j\|_2) }{\|\cD_{n+1}\|_1} + \frac{2(\|X_i\|_2 + \|X_j\|_2 )}{\|\cD_{n+1}\|_1 }  = \frac{4(\|X_i\|_2 + \|X_j\|_2) }{\|\cD_{n+1}\|_1},
			\end{split}
		\end{equation} where in (a), we use the fact for any real numbers $ 0 < a < b$, then for any $c > 0$, we have $\frac{a}{b} \leq \frac{a+c}{b+c}$. By plugging the upper bounds of (A) and (B) into \eqref{ineq:case-2-bound}, we have
		\begin{equation*}
			\begin{split}
				&\frac{1}{(n+1)^2} \sum_{1 \leq i, j \leq n+1} \|\wh{\theta}(\cD_{n+1}^{\setminus i} )  - \wh{\theta}(\cD_{n+1}^{\setminus j} )\|_2  \leq \frac{4}{(n+1)^2}  \sum_{1 \leq i < j \leq n} \frac{6(\|X_i\|_2 + \|X_j\|_2) }{\|\cD_{n+1}\|_1}\\
				& = \frac{24}{(n+1)^2} \frac{n \|\cD_{n+1}\|_1}{\|\cD_{n+1}\|_1} \leq \frac{24}{n+1}.
			\end{split}
		\end{equation*}

		\item (Case 3: there exists $k \in [n+1]$ such that $\|\cD_{n+1}^{\setminus i} \|_1 \leq 2n$ for all $i > k$ and $\|\cD_{n+1}^{\setminus i} \|_1 > 2n$ for all $i \leq k$) Notice that when $k = n+1$, it is the Case 2, we have shown the result. Next, let us consider $k \leq n$. 
	\begin{equation} \label{ineq:case3-bound}
	\begin{split}
			&\frac{1}{(n+1)^2} \sum_{1 \leq i, j \leq n+1} \|\wh{\theta}(\cD_{n+1}^{\setminus i} )  - \wh{\theta}(\cD_{n+1}^{\setminus j} )\|_2  \\
				&= \frac{2}{(n+1)^2} \sum_{1 \leq i < j \leq n+1} \|\wh{\theta}(\cD_{n+1}^{\setminus i} )  - \wh{\theta}(\cD_{n+1}^{\setminus j} )\|_2 \\
				& = \frac{2}{(n+1)^2} \Big( \sum_{1 \leq i < j \leq k} \underbrace{\|\wh{\theta}(\cD_{n+1}^{\setminus i} )  - \wh{\theta}(\cD_{n+1}^{\setminus j} )\|_2}_{(C)} + \sum_{k < i < j \leq n+1} \underbrace{\|\wh{\theta}(\cD_{n+1}^{\setminus i} )  - \wh{\theta}(\cD_{n+1}^{\setminus j} )\|_2}_{(D)} \\
				& +\sum_{1 \leq i \leq k < j \leq n+1} \underbrace{\|\wh{\theta}(\cD_{n+1}^{\setminus i} )  - \wh{\theta}(\cD_{n+1}^{\setminus j} )\|_2}_{(E)} \Big).
	\end{split}
\end{equation}	Again, we bound (C), (D) and (E) separately below.
\begin{equation*}
	\begin{split}
		(C) =  \left\| \frac{2( \sum_{z \neq i} X_z )}{\|\cD_{n+1}^{\setminus i}\|_1}  - \frac{2( \sum_{z \neq j} X_z )}{\|\cD_{n+1}^{\setminus j}\|_1} \right\|_2 \leq \frac{12(\|X_i\|_2 + \|X_j\|_2) }{\|\cD_{n+1}\|_1},
	\end{split}
\end{equation*} where the inequality follows the same analysis as in \eqref{ineq:case2-A} and \eqref{ineq:case2-B}. In addition, 
\begin{equation*}
	\begin{split}
		(D) = \frac{\|X_j - X_i\|_2}{n} \leq \frac{\|X_j\|_2 + \|X_i\|_2}{n} \leq \frac{4(\|X_j\|_2 + \|X_i\|_2)}{\|\cD_{n+1}\|_1},
	\end{split}
\end{equation*} where the last inequality is because when $i,j > k$, $\|\cD_{n+1}^{\setminus i}\|_1 \leq 2n, \|\cD_{n+1}^{\setminus j}\|_1 \leq 2n$, so $\|\cD_{n+1}\|_1 \leq \|\cD_{n+1}^{\setminus i}\|_1 + \|\cD_{n+1}^{\setminus j}\|_1 \leq 4n$.
		
		Finally, in the (E) term, we have $\|\cD_{n+1}^{\setminus i}\|_1 > 2n$ and $\|\cD_{n+1}^{\setminus j}\|_1 \leq 2n$. Now
		\begin{equation*}
			\begin{split}
				(E) &= \left\|  \frac{2( \sum_{z \neq i} X_z )}{\|\cD_{n+1}^{\setminus i}\|_1} - \frac{\sum_{z \neq j} X_z}{n}  \right\|_2 = \left\|  \frac{2( \sum_{z \neq i} X_z )n - (\sum_{z \neq j} X_z) \|\cD_{n+1}^{\setminus i}\|_1 }{\|\cD_{n+1}^{\setminus i}\|_1 n}  \right\|_2 \\
				& = \left\|  \frac{( \sum_{z \neq i,j} X_z )(2n - \|\cD_{n+1}^{\setminus i}\|_1) + 2nX_j - X_i \|\cD_{n+1}^{\setminus i}\|_1 }{\|\cD_{n+1}^{\setminus i}\|_1 n}  \right\|_2\\
				& \leq \left\|  \frac{( \sum_{z \neq i,j} X_z )(2n - \|\cD_{n+1}^{\setminus i}\|_1) }{\|\cD_{n+1}^{\setminus i}\|_1 n} \right\|_2 + \left\| \frac{ 2nX_j - X_i \|\cD_{n+1}^{\setminus i}\|_1 }{\|\cD_{n+1}^{\setminus i}\|_1 n} \right\|_2 \\
				& \leq \|\cD_{n+1}^{\setminus i,j}\|_1   \frac{\|\cD_{n+1}^{\setminus i}\|_1 - 2n}{\|\cD_{n+1}^{\setminus i}\|_1 n}  + \frac{2\|X_j\|_2}{\|\cD_{n+1}^{\setminus i}\|_1} + \frac{\|X_i\|_2}{n} \\
				& \overset{(a)}\leq 2  \frac{\|\cD_{n+1}^{\setminus i}\|_1 - 2n}{\|\cD_{n+1}^{\setminus i}\|_1} + \frac{2\|X_j\|_2}{\|\cD_{n+1}^{\setminus i}\|_1} + \frac{2\|X_i\|_2}{\|\cD_{n+1}^{\setminus j}\|_1}  \\
				& \overset{(a)}\leq 2  \frac{\|\cD_{n+1}^{\setminus i}\|_1 - \|\cD_{n+1}^{\setminus j}\|_1}{\|\cD_{n+1}^{\setminus i}\|_1} + \frac{2\|X_j\|_2}{\|\cD_{n+1}^{\setminus i}\|_1} + \frac{2\|X_i\|_2}{\|\cD_{n+1}^{\setminus j}\|_1}  \\
				& \leq 2  \frac{\|X_j\|_2 - \|X_i\|_2}{\|\cD_{n+1}^{\setminus i}\|_1} + \frac{2\|X_j\|_2}{\|\cD_{n+1}^{\setminus i}\|_1} + \frac{2\|X_i\|_2}{\|\cD_{n+1}^{\setminus j}\|_1} = \frac{4\|X_j\|_2}{\|\cD_{n+1}^{\setminus i}\|_1} \\
				& \overset{(b)}\leq 4\frac{\|X_j\|_2 + \|X_i\|_2}{\|\cD_{n+1}^{\setminus i}\|_1 + \|X_i\|_2} = 4\frac{\|X_j\|_2 + \|X_i\|_2}{\|\cD_{n+1}\|_1}
			\end{split}
		\end{equation*} where in (a) we use the fact $\|\cD_{n+1}^{\setminus j}\|_1 \leq 2n$ and in (b), we use the fact for any real numbers $ 0 < a < b$, then for any $c > 0$, we have $\frac{a}{b} \leq \frac{a+c}{b+c}$.
		
		In summary, based on the bounds for (C), (D) and (E), \eqref{ineq:case3-bound} can be bounded further as
		\begin{equation*}
			\begin{split}
				&\frac{1}{(n+1)^2} \sum_{1 \leq i, j \leq n+1} \|\wh{\theta}(\cD_{n+1}^{\setminus i} )  - \wh{\theta}(\cD_{n+1}^{\setminus j} )\|_2  \leq \frac{2}{(n+1)^2} \sum_{1 \leq i < j \leq n+1} 12\frac{\|X_j\|_2 + \|X_i\|_2}{\|\cD_{n+1}\|_1}\\
				& = \frac{24}{(n+1)^2} \frac{n \|\cD_{n+1}\|_1}{\|\cD_{n+1}\|_1} = \frac{24}{n+1}.
			\end{split}
		\end{equation*}
This finishes the proof of this lemma.
\end{itemize}

\end{proof}

\section{Proof for Section \ref{sec:sparse-mean}: Theorem \ref{th:sparse-mean}}

\begin{proof}[Proof of Theorem \ref{th:sparse-mean}]
We first argue the upper bound, followed by the lower bound.
	\paragraph{Upper bound.}
	The sparsity-ensured soft-thresholding estimator
	is $s$-sparse by construction and has worst-case stability of  $\beta_n$ by Lemma~\ref{lem:sparse-thresh-stability}, which implies average-case ($\ell_p,\, p\geq 1$) stability  as well. When $\beta_n \geq 4\sqrt{2s}r/n$, its high-probability and in-expectation $\ell_2$ risks are both of order $r^2 s\log d/n$ by Theorem~\ref{thm:sparse-thresh-risk} and Corollary~\ref{cor:sparse-thresh-expectation}. Hence $\forall p \, \geq 1$,
	\[
	R_{n,p}\bigl(\theta(\cP^d(r,s)),\beta_n\bigr)\ \lesssim\ \frac{r^2 s\log d}{n}
	\qquad\text{whenever }\beta_n \geq \frac{4\sqrt 2r\sqrt{s}}{n}.
	\]
When $\beta_n \leq \tfrac{4\sqrt{2} \cdot r\sqrt{s}}{n}$, the risk admits a trivial upper bound. Indeed, each coordinate of the estimator is bounded by $r$ and the estimator itself is $s$-sparse, due to the boundedness of the data. Consequently, the overall risk is at most on the order of $r^2 s$, yielding
	\[
	R_{n,p}\bigl(\theta(\cP^d(r,s)),\beta_n\bigr)\ \lesssim\ r^2 s.
	\]
	
	\paragraph{Lower bound (average case, $p=1$).}
	Applying Corollary~\ref{coro:average-case-stability} with the pair
	$\theta_0=0$ and $\theta_1=r\sum_{j=1}^s e_j$ (notice that we can easily find $P_0$ and $P_1$ such that $(1-t)P_0 + tP_1 \in \cP^d(r,s)$ for all $t \in [0,1]$) yields
	\[
	R_{n,1}\bigl(\theta(\cP^d(r,s)),\beta_n\bigr)
	\ \ge\
	\Bigl(\tfrac{r\sqrt{s}}{4}-2(n+1)\beta_n\Bigr)_+^{\,2}
	\ \vee\ 
	R\bigl(\theta(\cP^d(r,s))\bigr).
	\]
	Since $R_n(\theta(\cP^d(r,s)))\asymp r^2 s\log d/n$, this gives the lower branch
	$R_{n,1}(\cdot)\gtrsim r^2 s$ when $\beta_n\leq \frac{r\sqrt{s}}{8(n+1)}$, and it matches the optimal unconstrained rate above this threshold up to constants.
	
	\paragraph{Lower bound (for $p > 1$).}
	Since $R_{n,p}(\cdot,\beta_n)$ is a nondecreasing function of $p$, the same lower bounds hold for $p \geq 1$.
	
	Combining the upper and lower bounds establishes the sharp phase transition at $\beta_n^\star\asymp r\sqrt{s}/n$ for $\ell_p$ stability for all $p \geq 1$.
\end{proof}

\begin{Lemma}[Worst-case stability of the sparsity-ensured soft-thresholding estimator]\label{lem:sparse-thresh-stability}
Fix $\beta_n >0$, let  $X\in\bbB_\infty(r)\subseteq\bbR^d$, and let $\cD_n,\cD_n'$ be neighboring datasets (they differ in exactly one sample). 
Consider the data-driven soft-thresholding estimator in~\eqref{eq:data-driven-soft} 
evaluated on $\cD_n,\cD_n'$ given by $\wh\theta_j(\cD_n)$ and $\wh\theta_j(\cD_n')$. Then
\[
\bigl\|\wh\theta(\cD_n)-\wh\theta(\cD_n')\bigr\|_2\ \le\ \beta_n.
\]
Hence $\wh\theta$ is $\beta_n$-worst-case-stable (Definition~\ref{def:worst-cases-stability}) .
\end{Lemma}
\begin{proof}
Let 
\[
c\ :=\ \left(1 \wedge \frac{n\beta_n}{4\sqrt{2s}\,r}\right)\in(0,1].
\]
By definition,
\[
\wh\theta_j(\cD_n)=c\cdot \mathrm{sign}(\wb X_j)\big(|\wb X_j|-\tau_s\big)_+,\qquad
\wh\theta_j(\cD_n')=c\cdot \mathrm{sign}(\wb X'_j)\big(|\wb X'_j|-\tau_s'\big)_+,
\] where $\tau_s = |\bar{X}|_{(s+1)}$ and $\tau'_s = |\bar{X}'|_{(s+1)}$

\paragraph{Bound the change in sample means.}
Since $\cD_n$ and $\cD_n'$ differ in one observation and $X\in\bbB_\infty(r)$,
\[
\Delta\ :=\ \|\wb X-\wb X'\|_\infty\ \le\ \frac{2r}{n}.
\]

\paragraph{Order-statistic stability.}
Let $s_j:=|\wb X_j|$ and $s'_j:=|\wb X'_j|$, with order statistics $s_{(1)}\ge\cdots\ge s_{(d)}$ and $s'_{(1)}\ge\cdots\ge s'_{(d)}$.
Since $x\mapsto|x|$ is $1$-Lipschitz, $|s_j-s'_j|\le \Delta$ for all $j$. Hence, for each $k$,
\[
\big|\,s_{(k)}-s'_{(k)}\,\big|\ \le\ \Delta,
\qquad\text{in particular}\qquad
|\tau_s-\tau_s'|\ =\ \big|\,s_{(s+1)}-s'_{(s+1)}\,\big|\ \le\ \Delta.
\]

\paragraph{Coordinatewise change of the thresholded magnitudes.}
Define $a_j:=(|\wb X_j|-\tau_s)_+$ and $a'_j:=(|\wb X'_j|-\tau_s')_+$. Since $x\mapsto (x)_+$ is $1$-Lipschitz,
\[
|a_j-a'_j|
\ \le\ \big||\wb X_j|-|\wb X'_j|\big| + |\tau_s-\tau_s'|
\ \le\ \Delta+\Delta
\ =\ 2\Delta.
\]

\paragraph{Translate to the signed (and scaled) estimator; handle sign flips.}
We have $\wh\theta_j=c\cdot \mathrm{sign}(\wb X_j)\,a_j$ and $\wh\theta'_j=c\cdot \mathrm{sign}(\wb X'_j)\,a'_j$.

\emph{(i) No sign flip:} If $\mathrm{sign}(\wb X_j)=\mathrm{sign}(\wb X'_j)$, then
\[
|\wh\theta_j-\wh\theta'_j|
\ =\ c\,|a_j-a'_j|
\ \le\ c\cdot 2\Delta.
\]

\emph{(ii) Sign flip:} If $\mathrm{sign}(\wb X_j)\neq \mathrm{sign}(\wb X'_j)$, then
$|\wb X_j|+|\wb X'_j|\le |\wb X_j-\wb X'_j|\le \Delta$, so $|\wb X_j|\le\Delta$ and $|\wb X'_j|\le\Delta$.
Therefore $a_j\le (\Delta-\tau_s)_+\le \Delta$ and $a'_j\le (\Delta-\tau_s')_+\le \Delta$, yielding
\[
|\wh\theta_j-\wh\theta'_j|
\ \le\ c\,(a_j+a'_j)
\ \le\ 2c\Delta.
\]

Combining (i)-(ii), for all $j$,
\[
|\wh\theta_j-\wh\theta'_j|\ \le\ 2c\Delta\ \le\ 2c\cdot \frac{2r}{n}\ =\ \frac{4cr}{n}.
\]

\paragraph{Sparsity of the outputs and the $\ell_2$ bound.}
By construction $\tau_s=s_{(s+1)}$, hence $a_j>0$ only on at most $s$ coordinates; the same holds for $a'_j$.
Thus the union of supports of $\wh\theta(\cD_n)$ and $\wh\theta(\cD_n')$ has size at most $2s$, and by Cauchy--Schwarz,
\[
\|\wh\theta(\cD_n)-\wh\theta(\cD_n')\|_2
\ \le\ \sqrt{2s}\cdot \max_{1\le j\le d}|\wh\theta_j-\wh\theta'_j|
\ \le\ \sqrt{2s}\cdot \frac{4cr}{n}
\ =\ c\cdot \frac{4\sqrt{2s}\,r}{n}.
\]
Finally, by the definition of $c$,
\[
c\cdot \frac{4\sqrt{2s}\,r}{n}\ \le\ \beta_n.
\]
Hence $\|\wh\theta(\cD_n)-\wh\theta(\cD_n')\|_2\le \beta_n$, proving $\beta_n$-worst-case stability.
\end{proof}

\begin{Theorem}[Tail bound with purely data--driven cutoff $\tau_s=|\wb X|_{(s+1)}$]\label{thm:sparse-thresh-risk}
Consider data from the sparse mean model  $X_1,\dots,X_n \stackrel{\text{i.i.d.}}{\sim} P \in \cP^d(r,s)$ with  mean parameter $\theta(P)\in\bbR^d$. Let $\wh \theta$ be the estimator defined in~\eqref{eq:data-driven-soft} and assume that $\beta_n \geq \frac{4\sqrt 2 r\sqrt s}{n} $, then for every $u\ge 0$,
\[
\bbP\!\left(\ \|\wh\theta-\theta(P)\|_2 \;>\; 8\,r\,\sqrt{\frac{s\,(\log(4d)+u)}{n}}\ \right)
\;\le\; e^{-u}.
\]
\end{Theorem}

\begin{proof}
When $\beta_n \geq \frac{4\sqrt 2 r\sqrt s}{n}$, then the estimator reduces to
\[
\wh\theta_j \;=\; \mathrm{sign}(\wb X_j)\,\bigl(|\wb X_j|-\tau_s\bigr)_+,\qquad
\tau_s \;:=\; |\wb X|_{(s+1)},
\]
where $|\wb X|_{(s+1)}$ is the $(s{+}1)$-st order statistic of $\{|\wb X_j|\}_{j=1}^d$.
Write $\theta:=\theta(P)$ and fix $u\ge 0$. Set
\[
t_d(u)\ :=\ r\sqrt{\frac{2(\log(4d)+u)}{n}}.
\]

\paragraph{Coordinate concentration.}
For any $j\in[d]$, since $X_{ij}\in[-r,r]$, Hoeffding's inequality gives
\[
\bbP\big(|\wb X_j-\theta_j|>t_d(u)\big)\ \le\ 2\exp\!\Big(-\frac{2n\,t_d(u)^2}{(2r)^2}\Big)\ =\ 2e^{-(\log(4d)+u)}.
\]
A union bound over $j=1,\dots,d$ yields
\[
\bbP\!\Big(\ \|\wb X-\theta\|_\infty \le t_d(u)\ \Big)\ \ge\ 1-4d\,e^{-(\log(4d)+u)}\ =\ 1-e^{-u}.
\]
Define the high-probability event
\[
\mathcal E_u\ :=\ \bigl\{\ \|\wb X-\theta\|_\infty \le t_d(u)\ \bigr\}.
\]

\paragraph{$\tau_s$ is small on $\mathcal E_u$.}
Let $S^\star:=\{j:\theta_j\neq 0\}$ so $|S^\star|\le s$. On $\mathcal E_u$, if $j\notin S^\star$ then $\theta_j=0$ and hence $|\wb X_j|\le t_d(u)$. Thus \emph{at most} the $s$ indices in $S^\star$ can exceed $t_d(u)$, which forces
\[
|\wb X|_{(s+1)}\ \le\ t_d(u)\qquad\text{on }\mathcal E_u.
\]
Since $\tau_s=|\wb X|_{(s+1)}$, we conclude
\[
0\ \le\ \tau_s\ \le\ t_d(u)\qquad\text{on }\mathcal E_u.
\]

\paragraph{Selected support has size $\le s$.}
By definition,
\[
\widehat S\ :=\ \{j:\wh\theta_j\neq 0\}\ =\ \{j:|\wb X_j|>\tau_s\}.
\]
Because $\tau_s$ is the $(s{+}1)$-st largest magnitude, at most $s$ coordinates can exceed it; hence
\[
|\widehat S|\ \le\ s\qquad\text{(deterministically)}.
\]

\paragraph{Coordinatewise error on $\mathcal E_u$.}
We bound $|\wh\theta_j-\theta_j|$ in the three disjoint cases below, always using $\tau_s\le t_d(u)$ and $||\wb X_j|-|\theta_j||\le t_d(u)$ on $\mathcal E_u$.

\emph{(a) $j \in \widehat S$.} Then $|\bar X_j| > \tau_s$.  
We analyze the coordinate-wise error:
\[
\begin{aligned}
|\widehat\theta_j - \theta_j|
&= \bigl|\,\mathrm{sign}(\bar X_j)\,(|\bar X_j| - \tau_s) - \theta_j\,\bigr| \\[2pt]
&\le \bigl|\,|\bar X_j| - |\theta_j|\,\bigr| + \tau_s
  + 2|\theta_j|\,\mathbf{1}\!\left\{\mathrm{sign}(\bar X_j) \ne \mathrm{sign}(\theta_j)\right\}.
\end{aligned}
\]
On the event $\mathcal E_u$, we have 
$|\,|\bar X_j| - |\theta_j|\,| \le t_d(u)$ and $\tau_s \le t_d(u)$.  
Moreover, a sign flip implies $|\theta_j| \leq|\bar X_j - \theta_j| \le t_d(u)$, hence $|\theta_j| \le t_d(u)$.  
Combining these, we obtain
\[
|\widehat\theta_j - \theta_j| \;\le\; 4\,t_d(u)
\quad\text{for all } j \in \widehat S \text{ on } \mathcal E_u.
\]

\emph{(b) $j\in S^\star\setminus \widehat S$}. Then $|\wb X_j|\le \tau_s$ and $\wh\theta_j=0$, so
\[
|\wh\theta_j-\theta_j|
=|\theta_j|
\le |\theta_j|-|\wb X_j|+|\wb X_j|
\le t_d(u)+\tau_s
\le 2\,t_d(u).
\]

\emph{(c) $j\notin S^\star\cup \widehat S$}. Then $\theta_j=0$ and $|\wb X_j|\le \tau_s$, hence $\wh\theta_j=0$ and the error is $0$.

Thus, on $\mathcal E_u$, for all $j\in \widehat S\cup S^\star$,
\[
|\wh\theta_j-\theta_j|\ \le\ 4\,t_d(u),
\qquad\text{and}\qquad
\wh\theta_j-\theta_j=0\ \text{if }j\notin \widehat S\cup S^\star.
\]

\paragraph{$\ell_2$ error on $\mathcal E_u$.}
Since $|\widehat S|\le s$ and $|S^\star|\le s$, we have $|\widehat S\cup S^\star|\le 2s$. Therefore,
\[
\|\wh\theta-\theta\|_2
=\Big(\sum_{j\in \widehat S\cup S^\star}|\wh\theta_j-\theta_j|^2\Big)^{1/2}
\le \sqrt{2s}\cdot (4\,t_d(u))
= 8\,r\,\sqrt{\frac{s\,(\log(4d)+u)}{n}},
\]
on $\mathcal E_u$.

\paragraph{Tail.}
Finally,
\[
\bbP\!\left(\ \|\wh\theta-\theta\|_2 \;>\; 8\,r\,\sqrt{\frac{s\,(\log(4d)+u)}{n}}\ \right)
\ \le\ \bbP(\mathcal E_u^c)
\ \le\ e^{-u}.
\]
\end{proof}

\begin{Corollary}[In-expectation risk]\label{cor:sparse-thresh-expectation}
Under the assumptions of Theorem~\ref{thm:sparse-thresh-risk},
\[
\bbE\big\|\wh\theta-\theta(P)\big\|^2_2 \ \le\ C\,\frac{sr^2\log(4d)}{n},
\]
for a universal constant $C>0$.
\end{Corollary}

\begin{proof}
Let's proceed in a step-by-step fashion.
\paragraph{Tail bound to a usable form.}
Theorem~\ref{thm:sparse-thresh-risk} says that for every $u\ge 0$,
\[
\bbP\!\left(\ \|\wh\theta-\theta\|_2 \;>\; 8r\sqrt{\frac{s(\log(4d)+u)}{n}}\ \right)\ \le\ e^{-u}.
\]
Let $Z:=\|\wh\theta-\theta\|_2$. Setting $t=8r\sqrt{s(\log(4d)+u)/n}$ and solving for $u$ gives
$u=\frac{n t^2}{64 r^2 s}-\log(4d)$. Hence, for every $t>0$,
\[
\bbP(Z>t)\ \le\ \min\!\left\{\,1,\ 4d\exp\!\Big(-\frac{n t^2}{64 r^2 s}\Big)\right\}.
\]
Define $a:=\frac{n}{64 r^2 s}$ for brevity.

\smallskip
\paragraph{Pick the splitting point where the two pieces match.}
Choose
\[
t_0\ :=\ 8r\sqrt{\frac{s\log(4d)}{n}},
\]
so that $4d\,e^{-a t_0^2}=1$. Then
\[
\bbP(Z>t)\ \le\ \begin{cases}
1, & 0<t\le t_0,\\[2pt]
4d\,e^{-a t^2}, & t>t_0.
\end{cases}
\]

\smallskip
\paragraph{Tail-integration identity.}
For any nonnegative random variable $Z$,
$\ \bbE Z^2=2\displaystyle\int_0^\infty t\bbP(Z>t)\,dt$. Therefore,
\[
\bbE Z^2\ \le\ \underbrace{\int_0^{t_0} t_0\,dt}_{=\,t_0}\ +\ \int_{t_0}^{\infty} 8dt\,e^{-a t^2}\,dt.
\]

\smallskip
\paragraph{Bound the Gaussian tail integral.}
Using the standard bound $\displaystyle\int_x^\infty te^{-a t^2}dt \le \frac{e^{-a x^2}}{2 a }$ for $a,x>0$,
\[
\int_{t_0}^{\infty} 8dt\,e^{-a t^2}\,dt
\ \le\ 8d\cdot \frac{e^{-a t_0^2}}{2 a }
\ =\ \frac{1}{ a}\qquad\text{since }\ 4d\,e^{-a t_0^2}=1.
\]
With $a=\frac{n}{64 r^2 s}$ and $t_0=8r\sqrt{\frac{s\log(4d)}{n}}$, $\frac{1}{ a}
=\frac{64 r^2 s}{n}
$

\smallskip
\paragraph{Combine and simplify.}
Thus
\[
\bbE Z^2\ \le\ t^2_0\ +\ \frac{64 r^2 s}{n}
\ \lesssim \ \frac{sr^2\log(4d)}{n}.
\] This finishes the proof.
\end{proof}

\section{Proofs for Section \ref{sec:nonparametric}: Theorem \ref{th:worst-nonparametric}}

\subsection{Additional Properties of Besov Space} \label{app:besov-intro}
{ Here we collect the standard properties of the wavelet basis that are invoked throughout the proofs. As specified in Section~\ref{sec:nonparametric}, let $\phi$ and $\psi$ be the father (scaling) function and mother wavelet from a compactly supported Daubechies $\mathrm{db}N$ family, with $N$ sufficiently large that the regularity $A$ exceeds $s$; see \citet{daubechies1992ten}. The mother wavelet $\psi$ has $N$ vanishing moments. We use the associated Cohen--Daubechies--Vial construction \citep{cohen1993wavelets}, which gives a compactly supported, orthonormal, $A$-regular wavelet basis of $L_2[0,1]$ of the form
\[
\big\{\phi_{l_0+1,m}:0\le m<2^{l_0+1}\big\}
\ \cup\
\big\{\psi_{lk}:l\ge l_0+1,\ 0\le k<2^l\big\}.
\]
For interior indices,
\begin{equation} \label{eq:basis-construction}
	\phi_{l_0+1,m}(x)=2^{(l_0+1)/2}\phi(2^{l_0+1}x-m),
\qquad
\psi_{lk}(x)=2^{l/2}\psi(2^lx-k).
\end{equation}
For indices whose supports intersect the endpoints of $[0,1]$, the CDV construction replaces these translates by boundary-corrected functions having the same orthonormality, regularity, and localization properties.}

\paragraph{Wavelet basis properties.}
Let $\{\psi_{lk}\}$ denote a compactly supported, $A$-regular wavelet basis (with $A > s$) on $[0,1]$.  
{ By compact support and the standard localization properties of the CDV construction, both the interior and boundary-corrected functions satisfy the following bounds due to \eqref{eq:basis-construction}: there exist constants $C>0$ and $K_{\mathrm{loc}}>0$, independent of $x_0$, such that}
\begin{equation}\label{eq:wavelet-localization}
\begin{split}
	&\max(\|\psi_{lk}\|_\infty,\|\phi_{lk}\|_\infty) \le C\,2^{l/2},
\qquad
\#\{k\in\mathcal K_l:\ \psi_{lk}(x_0) \neq 0\} \le K_{\mathrm{loc}},\\
&
\max(\mathrm{Leb}(\mathrm{supp}(\psi_{lk})), \mathrm{Leb}(\mathrm{supp}(\phi_{lk}))) = O(2^{-l})
\quad \text{for all } l \ge l_0.
\end{split}
\end{equation}

The following lemma collects a few convenient properties we are going to use for the projection kervel $K_L(u,v) = \sum_{l=l_0}^{L} \sum_{k\in\mathcal K_l} \psi_{lk}(u) \psi_{lk}(v)$. 

\begin{Lemma}\label{lm:projection-kervel}
	For any $L \geq \ell_0$ and $x_0 \in (0,1)$, there exists a constant $C_0>0$ depending only on the choice of the wavelet basis such that
\begin{itemize}
	\item (i) $\|K_L(\cdot,x_0)\|_\infty \leq C_0 2^L$
	\item (ii) $\int_0^1 |K_L(x,x_0)| dx \leq C_0$
	\item (iii) $\int_0^1K_L(x,x_0)^2\,dx=K_L(x_0,x_0)\le C_02^L$.
\end{itemize}
\end{Lemma}
\begin{proof}
	To prove (i), use~\eqref{eq:wavelet-localization}, for fixed $l$ and any $x_0\in(0,1)$,
\[
\sum_{k\in\mathcal K_l}\big|\psi_{lk}(x)\psi_{lk}(x_0)\big|
\le \Big(\max_{k\in\mathcal K_l}|\psi_{lk}(x)|\Big)
\sum_{\substack{k\in\mathcal K_l:\\ \psi_{lk}(x_0)\neq 0}}|\psi_{lk}(x_0)|
\ \overset{\eqref{eq:wavelet-localization}}\le\ C\,2^{l/2}\cdot (K_{\rm loc}\cdot C\,2^{l/2})
= K_{\rm loc}C^2\,2^{l}.
\]
Summing over $l=l_0,\dots,L$ yields $\|K_L(\cdot,x_0)\|_\infty \ \le\ \sum_{l=l_0}^{L} C^2K_{\rm loc}\,2^{l}\ \lesssim 2^{L}$.

Conclusion (ii) is proved by using the multiresolution structure of subspaces formed by a wavelet basis. Specifically, the multiresolution analysis of the wavelet basis \citep{cohen1993multiresolution} says that for any $L \geq l_0$, the following two sets of basis functions form the same space:
\begin{equation*}
	\{ \psi_{lk}: l_0 \leq l \leq L, 0 \leq k \leq 2^{l}-1  \} \quad \text{ and } \quad \{ \phi_{L+1,k}: 0 \leq k \leq 2^{L+1}-1\}.
\end{equation*} Therefore, $K_L(x,x_0) = \sum_{l=l_0}^{L} \sum_{k\in\mathcal K_l} \psi_{lk}(x) \psi_{lk}(x_0) = \sum_{k \in \mathcal K_{L+1}} \phi_{L+1,k}(x) \phi_{L+1,k}(x_0)$. Then
\begin{equation*}
	\int_0^1 |K_L(x,x_0)| dx  = \int_0^1 \left|\sum_{k \in \mathcal K_{L+1}} \phi_{L+1,k}(x) \phi_{L+1,k}(x_0)\right| dx \overset{ \eqref{eq:wavelet-localization} }\leq  K_{\mathrm{loc}} C^2 2^{L+1} \cdot   \mathrm{Leb}(\mathrm{supp}(\phi_{L+1,k})) \leq C_0.
\end{equation*}

The final conclusion (iii) follows from the orthonormal property of the wavelet basis and (i).
\end{proof}

\subsection{Proofs for Theorem \ref{th:worst-nonparametric} Part (a)}
We claim that
\begin{equation*}
  R_{n,\infty}\left(\theta(\cB^s_{q_0,q_1}(1)), \beta_n\right) \asymp (n \beta_n)^{-2\nu} \wedge 1 + n^{- \frac{2\nu}{2\nu + 1}}.
\end{equation*}
The lower bound follows directly from Lemma~\ref{lem:np-lb-worst-case-stability}.  
For the upper bound, we invoke Lemma~\ref{lem:wavelet-stability}, which guarantees that the estimator~\eqref{eq:wavelet-est2} is worst-case $\beta_n$-stable, and then apply Lemma~\ref{lem:np-worst-case-risk} to control its risk.

\begin{Lemma}[Worst-case stability lower bound for nonparametric pointwise estimation]\label{lem:np-lb-worst-case-stability}
Consider the nonparametric regression model defined in~\eqref{eq:dgp-non-para}. Let $\cP=\{P_f:\ f\in\cB^s_{q_0,q_1}(1)\}$ and the target functional be $\theta(P_f)=f(x_0)$ for a fixed $x_0\in(0,1)$. Then for any $\beta_n\ge 0$,
\[
R_{n,\infty}\!\big(\theta(\cB^s_{q_0,q_1}(1)),\beta_n\big)
\;\;\ge\;\;
c\,((n\beta_n)^{-2\nu}\wedge 1)\ \ \vee\ \ n^{-2\nu/(2\nu+1)},
\]
where $c>0$ depends only on $(\sigma,s,q_0,q_1)$ and the chosen wavelet system, and the second term is the classical (unconstrained) minimax pointwise risk.
\end{Lemma}

\begin{proof}
We invoke Theorem~\ref{th:worst-lower-bound-multi}, which states that for any $f_0,f_1\in \cB^s_{q_0,q_1}(1)$ and any coupling $(\cD_n,\cD_n')$ of $P_{f_0}^{\otimes n},P_{f_1}^{\otimes n}$,
\[
R_{n,\infty}(\theta(\cB^s_{q_0,q_1}(1)),\beta_n)\ \ge\
\left[\frac{\big(|f_1(x_0)-f_0(x_0)|-\bbE d_{\Ham}(\cD_n,\cD_n')\,\beta_n\big)_{+}}{2}\right]^2
\ \vee\ R_n(\theta(\cB^s_{q_0,q_1}(1))).
\]

\paragraph{Two-point construction.} We borrow the following construction from proof of Theorem 4 in \cite{cai2003rates}.
Choose a compactly supported, function $g\ge 0$ , with compact support, $g(0)>0$, and $g\in \cB^s_{q_0,q_1}(1)$. For scale $\alpha>0$ and amplitude $\eta>0$, set
\[
h_{\alpha,\eta}(x):=\eta\,g\!\big(\alpha(x-x_0)\big),\qquad f_0\equiv 0,\quad f_1:=h_{\alpha,\eta}.
\]
For $\alpha$ large (so that the support of $h_{\alpha,\eta}$ lies in $(0,1)$) i.e $\exists C$ such that $\alpha >C$, and by Besov scaling (via wavelet characterization, Lemma 1 in \cite{cai2003rates}),
\[
\|h_{\alpha,\eta}\|_{\cB^s_{q_0,q_1}}\ \lesssim\ \eta\,\alpha^{\nu},\qquad \nu=s-1/q_0.
\]
Hence choosing $\eta=c_\star\,\alpha^{-\nu}$ with $c_\star>0$ small (depending only on $(s,q_0,q_1)$ and the basis) ensures $f_0,f_1\in\cB^s_{q_0,q_1}(1)$. { Let $A:=c_\star g(0)>0$. The separation at the target point is then
\[
|f_1(x_0)-f_0(x_0)|=A\alpha^{-\nu}.
\]
}

\paragraph{Single-sample total variation.}
Let $P_f$ denote the law of $(X,Y)$ under $f$.
Conditioning on $X=x$, $Y\sim N(f(x),\sigma^2)$. Using
\[
\mathrm{TV}\big(N(\mu_1,\sigma^2),N(\mu_2,\sigma^2)\big)
=2\Phi\!\Big(\tfrac{|\mu_1-\mu_2|}{2\sigma}\Big)-1
\ \le\ \frac{|\mu_1-\mu_2|}{\sigma\sqrt{2\pi}},
\]
and $X\sim\mathrm{Unif}[0,1]$,
\[
\mathrm{TV}(P_{f_1},P_{f_0})
=\bbE_X\ \mathrm{TV}\big(N(f_1(X),\sigma^2),N(0,\sigma^2)\big)
\le \frac{1}{\sigma\sqrt{2\pi}}\ \bbE|f_1(X)|
= C_\sigma \int_0^1 |h_{\alpha,\eta}(x)|\,dx.
\]
{ Since the support of $h_{\alpha,\eta}$ lies in $(0,1)$, a change of variables gives
\[
\int_0^1 |h_{\alpha,\eta}(x)|\,dx
=c_\star\|g\|_1\alpha^{-(\nu+1)}.
\]
Therefore, defining $B:=c_\star\|g\|_1/(\sigma\sqrt{2\pi})$, we have
\[
\mathrm{TV}(P_{f_1},P_{f_0})\ \le\ B\,\alpha^{-(\nu+1)},
\]
where $A,B>0$ depend only on $(\sigma,s,q_0,q_1)$, the wavelet system, and $g$.
}

\paragraph{Product coupling and Hamming cost.}
By maximal i.i.d.\ coupling of the $n$ draws, see \cite[Lemma 4.1.13]{roch2024modern},
\[
\bbE d_{\Ham}(\cD_n,\cD_n')\ =\ n\,\mathrm{TV}(P_{f_1},P_{f_0})
\ { \le\ B\,n\,\alpha^{-(\nu+1)}.}
\]

\paragraph{Plug and optimize.}
Theorem~\ref{th:worst-lower-bound-multi} specialized to $\theta(P_f)=f(x_0)$ yields
\[
R_{n,\infty}\!\big(\theta(\cB^s_{q_0,q_1}(1)),\beta_n\big)\ \ge\
\left[\frac{\Big(\underbrace{|f_1(x_0)-f_0(x_0)|}_{ = A\alpha^{-\nu}}
\ -\ \underbrace{\bbE d_{\Ham}(\cD_n,\cD_n')\,\beta_n}_{ \ \ \le\ B n\beta_n\,\alpha^{-(\nu+1)}}\Big)_{+}}{2}\right]^2
\ \vee\ R_n\!\big(\theta(\cB^s_{q_0,q_1}(1))\big).
\]
{ Let $\alpha_0\ge 1$ be a fixed constant sufficiently large that the support condition and the Besov scaling argument above hold for every $\alpha\ge\alpha_0$. Choose a fixed $K\ge 2B/A$ and set
\[
\alpha=\alpha_0+Kn\beta_n.
\]
Then $n\beta_n/\alpha\le 1/K$, and hence
\begin{align*}
|f_1(x_0)-f_0(x_0)|
-\bbE d_{\Ham}(\cD_n,\cD_n')\,\beta_n
&\ge A\alpha^{-\nu}-Bn\beta_n\alpha^{-(\nu+1)}\\
&=\alpha^{-\nu}\left(A-B\frac{n\beta_n}{\alpha}\right)\\
&\ge \alpha^{-\nu}\left(A-\frac{B}{K}\right)
\ge \frac{A}{2}\alpha^{-\nu}.
\end{align*}
Finally, $\alpha=\alpha_0+Kn\beta_n\asymp 1+n\beta_n$, so
\[
\alpha^{-\nu}\asymp (n\beta_n)^{-\nu}\wedge 1.
\]
Thus the positive part is uniformly bounded below by a constant multiple of $(n\beta_n)^{-\nu}\wedge 1$.}
Squaring yields
\[
R_{n,\infty}\!\big(\theta(\cB^s_{q_0,q_1}(1)),\beta_n\big)\ \ge\ c\,((n\beta_n)^{-2\nu}\wedge 1)
\ \vee\ R_n\!\big(\theta(\cB^s_{q_0,q_1}(1))\big),
\]
for some constant $c>0$ depending only on $(\sigma,s,q_0,q_1)$ and $g$. Using $R_n\!\big(\theta(\cB^s_{q_0,q_1}(1))\big)\asymp n^{-2\nu/(2\nu+1)}$ (e.g., \cite{cai2003rates}) completes the proof.
\end{proof}

\begin{Lemma}[Worst-case stability of the estimator]\label{lem:wavelet-stability}
Fix $x_0\in(0,1)$ and consider the estimator
	 \begin{equation*}
	 \wh{f}_{\cD_n}(x_0) = h^{-1} \left( \pi_h \left( \frac{1}{n} \sum_{i=1}^n K_{L}(X_i, x_0) [Y_i]_{M}\right) \right) , 
	 \end{equation*} where $h(t) = \bbE_{\xi \sim N(0, \sigma^2)}([\xi + t]_M)$ denotes a smooth and monotone transformation under Gaussian noise, $\pi_{h}(\cdot)$ denotes the Euclidean projection onto the interval $h([-M, M])$ and $L = \left\lfloor\log_2\left(c_{\psi}n \beta_n\right) \vee (\ell_0-1)\right \rfloor \wedge L_{\text{opt}}$. Here $M$ denotes the upper bound for $\|f\|_{\infty}$, $[\cdot]_M$ is the clipping operator at level $M$ (see its definition in Section \ref{sec:organization}). The wavelet basis is chosen to satisfy~\eqref{eq:wavelet-localization}. Then we can choose $c_\psi > 0$ depends only on $(\sigma,s,q_0,q_1)$ and the wavelet system such that $\widehat f(x_0)$ is $\beta_n$-worst-case-stable.
\end{Lemma}

\begin{proof}

If $  \log_2\!(c_\psi\,n\beta_n) < \ell_0 $, $L =\ell_0 -1$ which implies $\hat{f} = 0$ and hence it is trivially $\beta_n$ stable. For the rest of the proof, we focus on the other case $ \log_2\!(c_\psi\,n\beta_n)  \geq \ell_0$. It is easy to check that $h'(t) = \bbP_{\xi \sim N(0, \sigma^2)}(-M \leq t + \xi \leq M)$, so $h(t)$ is a smooth and strictly increasing function. Thus, $h^{-1}$ is well-defined. Moreover, $\inf_{|t| \leq M} h'(t)  = \bbP_{\xi \sim N(0, \sigma^2)}(-2M \leq \xi \leq 0) = \Phi(2M/\sigma) - 1/2 > 0$, where $\Phi(\cdot)$ denotes the cdf of the standard Gaussian. Therefore, $h^{-1}$ is Lipschitz and for any $u, v \in h([-M,M])$,
\begin{equation} \label{ineq:h-inver-lip}
	|h^{-1}(u) - h^{-1}(v)| \leq \frac{1}{\Phi(2M/\sigma) - 1/2} |u-v|
\end{equation} by the mean-value theorem.

Then for any neighboring datasets $\cD_n,\cD_n'$ differing only at index $j$,
\begin{equation*}
	\begin{split}
		\bigl|\widehat f_{\cD_n}(x_0)-\widehat f_{\cD_n'}(x_0)\bigr|
&\overset{(a)}\leq  \frac{1}{n(\Phi(2M/\sigma) - 1/2)}\,\bigl|\,[Y_j]_M K_L(X_j,x_0)-[Y'_j]_M K_L(X'_j,x_0)\,\bigr|
\\
&\le \frac{M}{n(\Phi(2M/\sigma) - 1/2)}\big(|K_L(X_j,x_0)|+|K_L(X'_j,x_0)|\big).
	\end{split}
\end{equation*}
 where in (a), we use \eqref{ineq:h-inver-lip} and the fact that $\pi_h(\cdot)$ is $1$-Lipschitz.
Taking the suprema over the neighboring pair and inputs gives
\begin{equation}\label{eq:stab_core_updated}
\sup_{\cD_n\sim \cD_n'}\bigl|\widehat f_{\cD_n}(x_0)-\widehat f_{\cD_n'}(x_0)\bigr|
\ \le\ \frac{2M}{n(\Phi(2M/\sigma) - 1/2)}\,\|K_L(\cdot,x_0)\|_\infty \overset{\text{Lemma  \ref{lm:projection-kervel}(i) }}\leq \frac{2M C_0}{n(\Phi(2M/\sigma) - 1/2)} \cdot 2^L.
\end{equation}

Recall that $L \leq  \left\lfloor\log_2\left(c_{\psi}n \beta_n\right)\right \rfloor$, then the conclusion follows by taking $c_\psi =(\Phi(2M/\sigma) - 1/2)/(2M C_0)$. 
\end{proof}

\begin{Lemma}[Risk of the transformed-clipped--truncated wavelet estimator]\label{lem:np-worst-case-risk}
	Let $f\in\mathcal{B}^s_{q_0,q_1}(1)$ with $\nu=s-1/q_0 > 0, q_0 \geq 1, q_1 \geq 1$, and fix $x_0\in(0,1)$. Consider the same estimator as in Lemma \ref{lem:wavelet-stability}, then
	\[
	\sup_{f\in\mathcal{B}^s_{q_0,q_1}(1)}\,
	\bbE\bigl(\widehat f(x_0)-f(x_0)\bigr)^2
	\;\;\lesssim\;\;
	(n\beta_n)^{-2\nu}\wedge 1
	\;+\;
	n^{-\frac{2\nu}{2\nu+1}}.
	\]

\end{Lemma}

\begin{proof} We proceed in a step-by-step fashion. First, if $\log_2\!(c_\psi\,n\beta_n) < \ell_0 $, then $\hat{f}_{\cD_n}(x_0) = 0$ and $\sup_{f \in \cB^s_{q_0,q_1}(1) } \bbE|\widehat{f}_{\cD_n}(x_0) - f(x_0)|^2 \lesssim 1$ since $f(x_0) \in [-M, M]$. For the rest of the proof, we focus on $\log_2\!(c_\psi\,n\beta_n) \geq \ell_0 $. Given $f(x_0) \in [-M, M]$, we have $h^{-1}(\pi_h(h\circ f(x_0))) = h^{-1}((h\circ f(x_0))) = f(x_0)$. Thus,
\begin{equation} \label{ineq:risk-nonpara-first-bound}
	\begin{split}
		|\widehat{f}_{\cD_n}(x_0) - f(x_0)| &= \left|  h^{-1} \left( \pi_h \left( \frac{1}{n} \sum_{i=1}^n K_{L}(X_i, x_0) [Y_i]_{M}\right) \right) - h^{-1}(\pi_h(h\circ f(x_0)))\right| \\
		& \leq \frac{1}{\Phi(2M/\sigma) - 1/2} \left| \frac{1}{n} \sum_{i=1}^n K_{L}(X_i, x_0) [Y_i]_{M} - h(f(x_0)) \right|,
	\end{split}
\end{equation} where the inequality is due to $h^{-1}(\cdot)$ is $ \frac{1}{\Phi(2M/\sigma) - 1/2}$-Lipschize as we have shown in the proof of Lemma \ref{lem:wavelet-stability} and the projection operator $\pi_h$ is $1$-Lipschitz.

Thus, to control $\bbE|\widehat{f}_{\cD_n}(x_0) - f(x_0)|^2$, we just need to bound the bias and variance of $\frac{1}{n} \sum_{i=1}^n K_{L}(X_i, x_0) [Y_i]_{M}$ as an estimator to $ h(f(x_0))$.

\vskip.2cm
{\bf \noindent Bias Control.} Since $h(t)$ is smooth and $\|f\|_{\infty} \leq M$, by the smooth composition of Besov function property, specifically Theorem 2.87 of \cite{bahouri2011fourier}, we have
\begin{equation*}
\begin{split}
	\|h\circ f - h(0) \|_{\cB^s_{q_0,q_1}} \leq C \|f\|_{\cB^s_{q_0,q_1}} .
\end{split}
\end{equation*} where the constant $C>0$ above only depends on $M$, $s$ and the derivative bound of $h$ on the interval $[-M,M]$. This implies that $\|h\circ f \|_{\cB^s_{q_0,q_1}} \leq \|h(0)\|_{\cB^s_{q_0,q_1}} + C \lesssim 1$ and
\begin{equation} \label{ineq:hf-besov-coeff}
	\|\left((h\circ f)_{lk}\right)_{k \in \mathcal{K}_l}\|_{\infty} \leq \|\left((h\circ f)_{lk}\right)_{k \in \mathcal{K}_l}\|_{q_0} \lesssim 2^{-l(s+1/2-1/q_0)}.
\end{equation}

 Denote $\Pi_L f = \sum_{l = l_0}^L \sum_{k \in \mathcal{K}_l } f_{lk} \psi_{lk}$ as the projection operator of $f$ at resolution level $L$. Then
\begin{equation*}
	\begin{split}
		\left| \bbE(K_L(X, x_0) [Y]_M) - h(f(x_0)) \right| &= \left| \bbE(K_L(X, x_0) [f(X) + \xi]_M) - h(f(x_0)) \right| \\
		& = \left| \bbE(K_L(X, x_0) h(f(X))) - h(f(x_0))\right| \\
		& = | \Pi_L (h \circ f)(x_0) -  h(f(x_0)) | \\
		& = \left| \sum_{l > L} \sum_{k \in \mathcal{K}_l} (h\circ f)_{lk} \psi_{lk}(x_0) \right| \\
		&  \overset{ \eqref{eq:wavelet-localization} }\leq \left|\sum_{l > L} K_{\mathrm{loc}}  (\max_{k \in \mathcal{K}_{l}}(h\circ f)_{lk})  C 2^{l/2} \right|  \overset{ \eqref{ineq:hf-besov-coeff} } \lesssim \sum_{l > L}  2^{-l \nu} \\
		& \lesssim 2^{-L \nu}.
	\end{split}
\end{equation*}

\vskip.2cm
{\bf \noindent Variance Control.}
\begin{equation*}
	\begin{split}
		\mathrm{Var}\left(\frac{1}{n} \sum_{i=1}^n K_L(X_i, x_0) [Y_i]_M\right) & = \frac{\mathrm{Var}\left(K_L(X, x_0) [Y]_M\right)}{n} \leq \frac{\bbE[K^2_L(X, x_0) [Y]^2_M]}{n} \\
		&\leq \frac{M^2\bbE[K^2_L(X, x_0)]}{n} \overset{\text{Lemma } \ref{lm:projection-kervel}(iii) }\lesssim \frac{2^L}{n}
	\end{split}
\end{equation*}

By combining the bias and variance above and \eqref{ineq:risk-nonpara-first-bound}, we get
\begin{equation*}
	\begin{split}
		\sup_{f \in \cB^s_{q_0,q_1}(1) } \bbE|\widehat{f}_{\cD_n}(x_0) - f(x_0)|^2 \lesssim 2^{-2 L \nu} + \frac{2^L}{n} \lesssim (n\beta_n)^{-2\nu}
	\;+\;
	n^{-\frac{2\nu}{2\nu+1}},
	\end{split}
\end{equation*} where in the last inequality, we plug in the value of $L = \left\lfloor\log_2\left(c_{\psi}n \beta_n\right)\right \rfloor \wedge L_{\text{opt}}$. 

Finally, combining with the result when $\log_2\!(c_\psi\,n\beta_n) < \ell_0 $, we get the unified upper bound
\begin{equation*}
	\begin{split}
		\sup_{f \in \cB^s_{q_0,q_1}(1) } \bbE|\widehat{f}_{\cD_n}(x_0) - f(x_0)|^2 \lesssim 2^{-2 L \nu} + \frac{2^L}{n} \lesssim (n\beta_n)^{-2\nu} \wedge 1
	\;+\;
	n^{-\frac{2\nu}{2\nu+1}}.
	\end{split}
\end{equation*}
\end{proof}

\subsection{Proofs for Theorem \ref{th:worst-nonparametric} Part (b) }
Under the average-case stability, we will show
\begin{equation*}
  R_{n,1}\left(\theta(\cB^s_{q_0,q_1}(1)), \beta_n\right)
  \;\asymp\;
  \begin{cases}
    n^{- \frac{2\nu}{2\nu + 1}}, & \text{if } \beta_n \ge \,\dfrac{24C_{\psi}}{n}, \\[0.8em]
    1, & \text{if } \beta_n \le \,\dfrac{c_{\star}}{n \log n}.
  \end{cases}
\end{equation*}
The upper bounds follow from Lemma \ref{lm:aver-stability-nonpara} and Lemma~\ref{lem:np-avg-case-rate}.  The corresponding lower bound is an immediate consequence of Lemma~\ref{lem:avg-lb-two-regime}.
\begin{Lemma}[Average-case stability lower bound: two-regime phase transition]\label{lem:avg-lb-two-regime}
Consider the nonparametric regression model~\eqref{eq:dgp-non-para}.
{Then there exist positive constants $c_0,c_1,c_{\star}>0$, depending only on $(\sigma,s,q_0,q_1)$ and the wavelet system, such that, the following holds under average-case ($p=1$) stability:}
\[
R_{n,1}\left(\theta(\cB^s_{q_0,q_1}(1)), \beta_n\right)\ \geq \begin{cases}
\ \frac{c_0^2}{16}\quad & \beta_n \le \dfrac{c_{\star}}{n\log n},\\
c_1n^{-\frac{2\nu}{2\nu+1}} \quad &\text{otherwise}
.   
\end{cases}
\]
\end{Lemma}

\begin{proof}
We invoke Theorem~\ref{th:average-lower-bound-multi} with $p=1$, which states that for any $\beta_n\ge 0$,
\[
R_{n,1}\left(\theta(\cB^s_{q_0,q_1}(1)), \beta_n\right)\ \ge\ 
\sup_{f_0,f_1\in \cB^s_{q_0,q_1}(1)}\left[\frac{\big( f_1(x_0) - f_0(x_0)-(n{+}1)(\log n{+}1)\beta_n\big)_+}{2}\right]^2
\ \vee\ R_n(\theta(\cB^s_{q_0,q_1}(1))).
\]
\textbf{Choice of the pair.} Let $f_0\equiv 0$ and $f_1\equiv c_0$ be constants. Constants belong to $\cB^s_{q_0,q_1}$ with norm proportional to $|c_0|$; taking $c_0>0$ small enough (depending only on $s,p,q$) ensures $f_0,f_1\in\cB^s_{q_0,q_1}(1)$. Therefore, the lower bound yields
\begin{equation}\label{eq:avg-lb-core}
R_{n,1}\left(\theta(\cB^s_{q_0,q_1}(1)), \beta_n\right)
\ \ge\
\left[\frac{\big(c_0-(n{+}1)(\log n{+}1)\beta_n\big)_+}{2}\right]^2
\ \vee\ R_n\!\big(\theta(\cB^s_{q_0,q_1}(1))\big).
\end{equation}

\paragraph{Regime (i).} Suppose $\beta_n \le \dfrac{c}{2(n+1)(\log n+1)}$, where we take $c$ to be sufficiently small depending only on $c_0$ such that
$c_0-(n{+}1)(\log n{+}1)\beta_n \ge c_0/2$. Then the first term in \eqref{eq:avg-lb-core} is at least $(c_0/4)^2$.

\paragraph{Regime (ii).}
\[
R_{n,1}\left(\theta(\cB^s_{q_0,q_1}(1)), \beta_n\right)\ \ge\ R_n\!\big(\theta(\cB^s_{q_0,q_1}(1))\big)\ \asymp\ n^{-\frac{2\nu}{2\nu+1}}.
\]
 This completes the proof.
\end{proof}

\begin{Lemma}[Average-case Stability bound] \label{lm:aver-stability-nonpara}
    The estimator defined in~\eqref{eq:wavelet-est-ave} is $\beta_n$-average-case-stable.
\end{Lemma}

\begin{proof} By applying Lemma~\ref{lm:average-stability-bound} with $X_i$ being replaced by $Z_i$, we have
\begin{equation*}
	\begin{split}
		\frac{1}{(n{+}1)^2}\sup_{\cD_{n+1}}
\sum_{i,j=1}^{n+1}
\big|\wh f_{\cD_{n+1}^{\setminus i}}(x_0)-\wh f_{\cD_{n+1}^{\setminus j}}(x_0)\big| \leq  \left(\frac{n\beta_n }{24C_{\psi}} \wedge 1\right) \cdot \frac{24 C_{\psi}}{n} \leq \beta_n.
	\end{split}
\end{equation*}
\end{proof}

\begin{Lemma}[Average-case estimator attains the unconstrained minimax rate]\label{lem:np-avg-case-rate}
Let $(X_i,Y_i)_{i=1}^n$ be i.i.d. data generated\ from~\eqref{eq:dgp-non-para}. Fix $x_0\in(0,1)$  and let $\{\psi_{lk}\}$ be a wavelet basis satisfying~\eqref{eq:wavelet-localization}.
Then for $\hat f(x_0)$ as defined in~\eqref{eq:wavelet-est-ave} we have
\[
\sup_{f\in\cB^s_{q_0,q_1}(1)}\ \mathbb{E}\bigl(\widehat f(x_0)-f(x_0)\bigr)^2\lesssim\begin{cases}
 \ \,n^{-\frac{2\nu}{2\nu+1}} \quad & \beta_n \geq \frac{24C_{\psi}}{n}\\
 \ 1 &\beta_n \leq \frac{1}{n}
    
\end{cases}
\] where $C_{\psi} >0$ depends only on $(\sigma,s,q_0,q_1)$
\end{Lemma}

\begin{proof} Let's proceed in a step-by-step fashion. First, when $\beta_n \le \frac{1 }{n}$. By design, we have $|\wh{f}_{\cD_n}(x_0)| \leq \frac{n \beta_n}{12} \leq \frac{1}{12}$. Thus, $|\wh{f}_{\cD_n}(x_0) - f(x_0)| \lesssim 1$ since $\|f\|_{\infty} \leq M$.

Now, we consider the regime $ \beta_n \geq \frac{24C_{\psi}}{n}$, where $C_{\psi} >0$ is some constant depending only on $(\sigma,s,q_0,q_1)$ to be specified later. In this regime, $ \wh{f}_{\cD_n}(x_0) = \left( \frac{2n C_{\psi}}{S} \wedge 1 \right) \bar{Z} $. We note that 
\begin{equation} \label{ineq:bias-nonpara}
	 \bbE|Z_i| = \bbE|(f(X_i) + \xi_i) K_L(X_i, x_0)|
 \le(M+\sigma\sqrt{2/\pi})
 \int_{0}^1 |K_L(x,x_0)|\,dx \overset{\text{Lemma  \ref{lm:projection-kervel} (ii) }}\leq (M+\sigma\sqrt{2/\pi})C_0.
\end{equation} In addition,
\begin{equation}\label{ineq:variance-nonpara}
	 \bbE Z_i^2
 \le(M^2+\sigma^2)\int_0^1 K_L(x,x_0)^2\,dx
\overset{\text{Lemma  \ref{lm:projection-kervel} (iii) }} \le (M^2+\sigma^2)C_0 2^L.
\end{equation} Take $C_{\psi} = (M+\sigma\sqrt{2/\pi})C_0$, then
\begin{equation} \label{ineq:prob-shrink}
	\begin{split}
		\bbP(S > 2n C_{\psi})\overset{ \eqref{ineq:bias-nonpara} } \leq \bbP(S - n \bbE|Z| > n C_{\psi} )\leq \frac{n \mathrm{Var}(|Z|) }{n^2 C^2_{\psi}} \overset{ \eqref{ineq:variance-nonpara} }\leq \frac{ (M^2+\sigma^2) 2^L}{4n(M+\sigma\sqrt{2/\pi})^2C_0},
	\end{split}
\end{equation} where the second inequality is by Chebyshev's inequality. Finally,
\begin{equation*}
	\begin{split}
		\mathbb{E}\bigl(\widehat f(x_0)-f(x_0)\bigr)^2 & = \mathbb{E}[\bigl(\widehat f(x_0)-f(x_0)\bigr)^2 \indi(S \leq 2nC_{\psi}) ] + \mathbb{E}[\bigl(\widehat f(x_0)-f(x_0)\bigr)^2 \indi(S > 2nC_{\psi}) ] \\
		&\overset{(a)} \lesssim 2^{-2\nu L} + \frac{2^L}{n} + \bbP(S > 2n C_{\psi}) \\
		& \overset{ \eqref{ineq:prob-shrink} }\lesssim 2^{-2\nu L} + \frac{2^L}{n} \lesssim n^{-\frac{2\nu}{2\nu+1}}.
	\end{split}
\end{equation*} where (a) is because when $S \leq 2nC_{\psi}$, $\widehat f(x_0)$ reduces to the classical truncated wavelet estimator in \eqref{eq:wavelet-est} and its bias and variance at resolution level $L$ follows from the usual bias-variance tradeoff for the unconstrained estimator (see the proof of Theorem 2 of \cite{cai2003rates}); moreover, on the event $S > 2nC_{\psi}$, $|\wh{f}_{\cD_n}(x_0)| \leq 2 C_{\psi}$, therefore, $|\wh{f}_{\cD_n}(x_0) - f(x_0)| \lesssim 1$ since $\|f\|_{\infty} \leq M$; the final inequality follows as we choose $L = L_{\mathrm{opt}}$. This finishes the proof of this lemma.
\end{proof}

\section{Additional Proofs for Section \ref{sec:discussion} }

\subsection{Proof of Lemma \ref{lm:privacy-to-stability} }
	 For any two neighboring datasets $\cD_n, \cD_n'$, by the DP guarantee of $M$, we know that given $\cD_n, \cD_n'$, $ e^{-\epsilon_n} \leq dM(\cD_n)/dM(\cD_n') \leq e^{\epsilon_n}$. Then 
\begin{equation*}
	\begin{split}
		|\wh{\theta}(\cD_n) - \wh{\theta}(\cD'_n)| &= \left| \int_{-r}^r m dM(\cD_n) - \int_{-r}^r m  dM(\cD'_n)  \right| = \left|\int_{-r}^r m\left(\frac{dM(\cD_n)}{dM(\cD_n')} - 1\right) d M(\cD_n') \right| \\
		& \leq  r \int_{-r}^r (e^{\epsilon_n }-1) d M(\cD_n')  = r (e^{\epsilon_n}-1).
	\end{split}
\end{equation*} For $\epsilon_n < 1$, we have $\exp(\epsilon_n) -1 \leq 2 \epsilon_n$.

\subsection{Proof of Proposition \ref{prop:convert-privacy-stability-curve} }
First, we note \eqref{eq:dp-stability} can be directly obtained from the proof of \eqref{eq:stability-dp}, so we omit it here. We just show \eqref{eq:stability-dp}. The upper bound in \eqref{eq:stability-dp} is a direct consequence of Lemma \ref{lm:stability-to-privacy}. Next, we show the lower bound in \eqref{eq:stability-dp}. Denote $\wh{\theta}(\cD_n) = \bbE[M(\cD_n)|\cD_n] $, then
\begin{equation*}
	\begin{split}
		 R_{\textnormal{DP}}( \theta(\cP) , \epsilon_n ) &= \inf_{M \textnormal{ satisfies } \epsilon_n-\DP } \sup_{P \in \cP} \bbE [ M(\cD_n) - \theta(P) ]^2 \\
		 & \geq \inf_{M \textnormal{ satisfies } \epsilon_n-\DP } \sup_{P \in \cP} \bbE \left[ \wh{\theta}(\cD_n) - \theta(P) \right]^2 \\
		 & = \inf_{\substack{M \textnormal{ satisfies } \epsilon_n-\DP \\ |M(\cD_n)| \leq r \textnormal{ a.s. } } } \sup_{P \in \cP} \bbE \left[ \wh{\theta}(\cD_n) - \theta(P) \right]^2 \quad(\textnormal{as } |\theta| \leq r)\\
		 & \overset{\textnormal{Lemma } \ref{lm:privacy-to-stability} }\geq \inf_{ \wh{\theta}(\cD_n) \textnormal{ is }r(e^{\epsilon_n}-1)-\textnormal{worst-case-stable} } \sup_{P \in \cP} \bbE \left[ \wh{\theta}(\cD_n) - \theta(P) \right]^2 \\
		 & = R_{n,\infty}(\theta(\cP), r(e^{\epsilon_n} - 1)).
	\end{split}
\end{equation*} This finishes the proof of this proposition.

{
\subsection{Proof of Theorem \ref{th:dis-risk}}
We note that the first part of the lower bound is the trivial unconstrained minimax lower bound. Thus, we focus on proving the second one.
Following the same proof as in Eq. \eqref{ineq:symmetrization} of Theorem \ref{th:average-lower-bound-multi}, it is without loss of generality to consider $\wh{\theta}$ is symmetric in the data points, as given any estimator, we can symmetrize it while maintaining the same stability level and have risk no bigger than the original one. For any estimator that is $\beta_n$-distributional-stable, let 
	\begin{equation*} 
		\sigma_n^2 := \sup_{P \in \cP} \bbE_{P} \|  \wh{\theta}(\cD_n;\xi) - \theta(P) \|_2^2,
	\end{equation*} where the expectation is taken with respect to $\cD_n \overset{i.i.d.}\sim P$ and $\xi$. Given any $P_0, P_1 \in \cP$ such that $P_t:=(1-t) P_0 + tP_1 \in \cP$ for all $t \in [0,1]$, for any $0<\eta<0.5$,
		\begin{equation} \label{ineq:mean-difference-lower-bound}
		\begin{split}
			&\|\bbE_{P_{\eta}}[ \wh{\theta}(\cD_n;\xi)] - \bbE_{P_{1-\eta}}[\wh{\theta}(\cD'_n; \xi) ] \|_2 \\
			&= \|\bbE_{P_{\eta}}[ \wh{\theta}(\cD_n;\xi) - \theta(P_{\eta})] + \theta(P_{\eta}) - \theta(P_{1-\eta}) + \bbE_{P_{1-\eta}}[\theta(P_{1-\eta}) - \wh{\theta}(\cD'_n; \xi) ] \|_2 \\
			& \geq \|\theta(P_{\eta}) - \theta(P_{1-\eta}) \|_2 - 2 \sigma_n.
		\end{split}
	\end{equation} On the other hand, we expect $\|\bbE_{P_{\eta}}[ \wh{\theta}(\cD_n;\xi)] - \bbE_{P_{1-\eta}}[\wh{\theta}(\cD'_n; \xi) ] \|_2$ would not be too large given that $\wh{\theta}$ is $\beta_n$-distributional-average-case stable. To show that, for any $t \in (0,1)$ and $x \in \cX$, let us define 
	\begin{equation*}
		g_t(x) = \bbE_{X_2,\ldots, X_n \overset{i.i.d}\sim P_t, \xi}[\wh{\theta}((x, X_2, \ldots, X_n); \xi)]\quad \text{and}\quad m(t) = \bbE_{P_t}[ \wh{\theta}(\cD_n;\xi)]. 
	\end{equation*} By definition,
	\begin{equation*}
		m(t) = \sum_{i=0}^n {n \choose i} (1-t)^i t^{n-i} \int_{x_1, \ldots, x_n} \bbE_{\xi}[\wh{\theta}((x_1, x_2, \ldots, x_n); \xi)] d P_0(x_1, \ldots, x_i) dP_1(x_{i+1}, \ldots, x_n).
	\end{equation*} Notice that $m(t)$ is a polynomial function of $t$. Therefore, 
	\begin{equation} \label{ineq:d-m-lower}
	\begin{split}
		m'(t) & = \left( \sum_{i=1}^n i {n \choose i} (1-t)^{i-1}t^{n-i} -  \sum_{i=0}^{n-1} (n-i) {n \choose i} (1-t)^i t^{n-i-1}\right) \\
		& \times  \int_{x_1, \ldots, x_n} \bbE_{\xi}[\wh{\theta}((x_1, x_2, \ldots, x_n); \xi)] d P_0(x_1, \ldots, x_i) dP_1(x_{i+1}, \ldots, x_n) \\
		& = n (\bbE_{P_1}[g_t(X)] - \bbE_{P_0}[g_t(X)]).
	\end{split}
	\end{equation} At the same time, distributional-average-case-stability implies that for any $t \in [0,1]$,
	\begin{equation} \label{ineq:mixture-lower}
		\begin{split}
			\beta_n \geq& \bbE_{X_1,\ldots, X_n,X_1' \overset{i.i.d.}\sim P_t; \xi} \|\wh{\theta}(X_1,X_2, \ldots, X_n; \xi) - \wh{\theta}(X_1',X_2, \ldots, X_n; \xi)\|_2 \\
			\geq & \bbE_{X_1, X_1' \overset{i.i.d.}\sim P_t} \|g_t(X_1) - g_t(X_1')\|_2 \geq 2t(1-t) \bbE_{X_1 \sim P_1, X_1'\sim P_0, X_1 \perp X_1'}\| g_t(X_1) - g_t(X_1') \|_2 \\
			\geq & 2t(1-t) \| \bbE_{P_1}[g_t(X)] - \bbE_{P_0}[g_t(X)]\|_2.
		\end{split}
	\end{equation} A combination of \eqref{ineq:d-m-lower} and \eqref{ineq:mixture-lower} yields that 
	\begin{equation} \label{ineq:mean-difference-upper-bound}
		\begin{split}
			\|\bbE_{P_{\eta}}[ \wh{\theta}(\cD_n;\xi)] - \bbE_{P_{1-\eta}}[\wh{\theta}(\cD'_n; \xi) ] \|_2 & = \|m(\eta) - m(1-\eta)\|_2 = \left \|\int_{\eta}^{1-\eta} m'(t)dt\right \|_2 \\
			& \leq\int_{\eta}^{1-\eta} \left \| m'(t)\right \|_2 dt \overset{ \eqref{ineq:d-m-lower} }\leq n \int_{\eta}^{1-\eta} \| \bbE_{P_1}[g_t(X)] - \bbE_{P_0}[g_t(X)]\|_2 dt\\
			& \overset{ \eqref{ineq:mixture-lower} }\leq \frac{n \beta_n}{2} \int_{\eta}^{1-\eta} \frac{1}{t(1-t)} dt \\
			& = \frac{n \beta_n}{2} \log\left( \frac{t}{1-t} \right)\Big|^{t=1-\eta}_{t=\eta} =  n \beta_n  \log \left( \frac{1-\eta}{\eta} \right).
		\end{split}
	\end{equation} Finally, a combination of \eqref{ineq:mean-difference-lower-bound} and \eqref{ineq:mean-difference-upper-bound} yields the second part of the lower bound by observing that \eqref{ineq:mean-difference-lower-bound} and \eqref{ineq:mean-difference-upper-bound} hold for any $0<\eta < 1/2$.
	
	When $\theta(P)$ is linear in $P$, then the second conclusion follows by taking $\eta = 1/4$ in the first conclusion.
}

\end{sloppypar}

\end{document}